\documentclass[12pt]{elsarticle}

\usepackage[utf8]{inputenc}
\usepackage{amsmath}
\usepackage{amsfonts}
\usepackage{amssymb}
\usepackage{amsthm}
\usepackage{graphicx}
\usepackage{xcolor}
\usepackage{arydshln}
\usepackage{enumitem}
\usepackage{hhline}
\usepackage{tikz}
\usetikzlibrary{automata,positioning,shapes}
\usepackage[hypcap=false]{caption}
\usepackage{mathtools}
\usepackage{multicol}
\usepackage{multirow}
\usepackage{bbm} 
\usepackage{bm} 
\usepackage{cancel}
\usepackage[colorlinks=true, linkcolor=blue, citecolor=blue, pdfauthor=author]{hyperref}
\usepackage[nameinlink]{cleveref}
\usepackage{ifthen}
\usetikzlibrary{calc}

\usepackage[left=2cm,right=2cm,top=1.5cm,bottom=1.5cm]{geometry}
\def\lbl#1{\label{#1}}

\renewcommand{\leq}{\leqslant}
\renewcommand{\geq}{\geqslant}
\newcommand{\eps}{\varepsilon}

\newcommand{\N}{\mathbbm N}
\newcommand{\Z}{\mathbbm Z}

\newcommand{\R}{\mathbbm R}
\newcommand{\C}{\mathbbm C}

\newtheorem{theorem}{Theorem}
\newtheorem{proposition}[theorem]{Proposition}
\newtheorem{lemma}[theorem]{Lemma}
\newtheorem{corollary}[theorem]{Corollary}

\newtheorem{definition}[theorem]{Definition}

\newtheorem*{remark}{Remark}

\newcommand{\thm}[1]{\Cref{thm:#1}}
\newcommand{\prop}[1]{\Cref{prop:#1}}
\newcommand{\lem}[1]{\Cref{lem:#1}}
\newcommand{\cor}[1]{\Cref{cor:#1}}
\newcommand{\defn}[1]{\Cref{def:#1}}

\newcommand{\fig}[1]{\Cref{fig:#1}}
\newcommand{\sect}[1]{\Cref{sec:#1}}

\makeatletter
\def\namedlabel#1#2{\begingroup
    #2%
    \def\@currentlabel{#2}%
    \phantomsection\label{#1}\endgroup
}
\makeatother

\newcommand{\lgray}{\color{lightgray}}

\newcommand{\mute}{}

\newcommand{\matnot}{\mathcal}
\newcommand{\Mat}[5]{\matnot{#1}_{#2,#3}^{#4}(#5)}
\newcommand{\Matminor}[3]{\matnot{#1}_{#2}^{#3}}
\newcommand{\Mattilde}[5]{\tilde{\matnot{#1}}_{#2,#3}^{#4}(#5)}

\newcommand{\diag}{\operatorname{diag}}
\newcommand{\detnew}[5]{#1_{#2,#3}^{#4}(#5)}
\newcommand{\dettilde}[5]{\tilde{#1}_{#2,#3}^{#4}(#5)}
\newcommand{\B}[4]{B_{#1,#2}^{#3}({#4})}
\newcommand{\D}[4]{D_{#1,#2}^{#3}({#4})}
\newcommand{\E}[4]{E_{#1,#2}^{#3}({#4})}

\newcommand{\Cof}[3]{\text{Cof}_{#1,#2}(#3-1)}

\begin{document}

\begin{frontmatter}

\title{Binomial Determinants for Tiling Problems\\ Yield to the Holonomic Ansatz}

\author[bupt]{Hao Du}
\author[ricam]{Christoph Koutschan}
\author[muic]{Thotsaporn Thanatipanonda}
\author[ricam]{Elaine Wong\corref{mycorrespondingauthor}}\ead{elaine.wong@ricam.oeaw.ac.at}

\cortext[mycorrespondingauthor]{Corresponding author}

\address[ricam]{Johann Radon Institute for Computational and Applied Mathematics (RICAM),\\ Austrian Academy of Sciences,
  Altenberger Stra\ss e 69, 4040 Linz, Austria}
\address[bupt]{School of Sciences, Beijing University of Posts and Telecommunications (BUPT),
  Beijing 100876, China}
\address[muic]{Science Division, Mahidol University International College (MUIC),
  Nakhonpathom, 73170, Thailand}

\nonumnote{This paper is published in the \textit{European Journal of Combinatorics}. DOI:10.1016/j.ejc.2021.103437}

\begin{abstract}
We present and prove closed form expressions for some families of binomial
determinants with signed Kronecker deltas that are located along an arbitrary
diagonal in the corresponding matrix. They count cyclically symmetric rhombus
tilings of hexagonal regions with triangular holes. We extend a previous
systematic study of these families, where the locations of the Kronecker
deltas depended on an additional parameter, to families with negative
Kronecker deltas. By adapting Zeilberger's holonomic ansatz to make it work
for our problems, we can take full advantage of computer algebra tools for
symbolic summation. This, together with the combinatorial interpretation,
allows us to realize some new determinantal relationships. From there, we are
able to resolve all remaining open conjectures related to these determinants,
including one from 2005 due to Lascoux and Krattenthaler.
\end{abstract}


\end{frontmatter}


\section{Introduction and History}
\lbl{sec:intro}

We tell a tale of two matrix families, whose determinants we want to calculate.
Suppose that $\mu$ is an indeterminate, $n\in\N$, and $s,t\in\Z$.
We define the matrices
\begin{align*}
  \Mat Dst{\mu}n &:= \begin{pmatrix}
    \binom{\mu+i+j+s+t-4}{j+t-1} + \delta_{i+s,j+t}
  \end{pmatrix}_{1\leq i, j \leq n},\\
  \Mat Est{\mu}n &:= \begin{pmatrix}
     \binom{\mu+i+j+s+t-4}{j+t-1} - \delta_{i+s,j+t}
  \end{pmatrix}_{1\leq i, j \leq n},
\end{align*}
and denote $\D st{\mu}n, \E st{\mu}n$ to be their corresponding determinants.
At a first glance, these two families appear almost the same: their entries
have the same binomial coefficient formula, with some entries (along a
diagonal) differing only by a genetic mutation of~$\pm 1$. The genealogy of
these families extends back to 1979 in a classic paper by
Andrews~\cite{Andrews79}, where we encounter the first result of the kind that
we will see in this paper, namely, that the determinant of a matrix from one
of the families has a closed form and counts certain combinatorial objects
(more precisely: descending plane partitions). 
For more background on plane partitions and their connections to
determinants up to the year 1999, see~\cite{Bressoud99}.

Then in 2005, Krattenthaler
published a rich collection of results and open
problems about determinants~\cite{Krattenthaler05}, containing four conjectures of a similar flavor
with various levels of difficulty: Problem 34 goes back to George Andrews, and
Conjectures 35--37 were formulated by Krattenthaler, Xin, and Lascoux. Two of
them were resolved by the second and third authors in 2013 \cite[Theorems~2
  and~5]{KoutschanThanatipanonda13} using Zeilberger's holonomic
ansatz~\cite{Zeilberger07} and automated tools for dealing with symbolic
sums~\cite{Koutschan09}. We briefly describe these techniques in
\sect{holonomicansatz} and \sect{computeralgebra}.

Despite the elegance and simplicity of the method, Problem 34 was only
partially resolved \cite[Theorem~1]{KoutschanThanatipanonda13}, and with the
introduction of an additional parameter, Conjecture 37 remained elusive even
with the available machinery. In their attempt to complete the work, the same
authors observed that the few determinants from their previous paper could be
generalized to infinite families that count cyclically symmetric rhombus
tilings of a hexagonal-shaped region with triangular holes. This is discussed
in much further detail in \sect{comb}, but the main idea originates from the
connection between counting rhombus tilings of a lozenge-shaped region and
counting non-intersecting lattice paths in the integer lattice, and using the
(known) fact that the latter are counted by determinants of binomial
coefficients.  The addition of the Kronecker deltas to the matrix complicates
the counting, as we have to consider all tuples of paths with certain selected
start and end points. This corresponds to adding up the number of rhombus
tilings of many different lozenge-shaped regions. Instead, we can construct a
single hexagonal region from three rotated copies of the original
lozenge-shaped regions, where the additional variations due to the Kronecker
deltas correspond to the presence or absence of rhombi crossing borders. To make this work, one has to enforce cyclic symmetry on the
rhombus tilings.

Armed with this interpretation,
a slew of new results was achieved in 2019, and more conjectures for these binomial
determinants were posed. Some of the results were proven using algebraic
manipulations and the computer as was done in~\cite{KoutschanThanatipanonda13}, but
also the combinatorial interpretation turned out to be crucial in a few of the
proofs.  Nevertheless, Conjecture 37 still resisted, as well as
newly introduced conjectures. We can summarize the exposition so far in
\Cref{tab:olddet}.

\begin{table}[ht]
\centering
\begin{tabular}{l|c|c|c}
  Determinant & First Proposed & Resolved & Year\\
  \hline
  $\D00{\mu}n$ & \cite[Theorem 8]{Andrews79} &
  \cite[Theorem 8]{Andrews79} & 1979 \\
  $\E11{\mu}n$ & \cite[Conjecture 35]{Krattenthaler05} &
  \cite[Theorem 2]{KoutschanThanatipanonda13} & 2013 \\
  $\E22{\mu}n$ & \cite[Conjecture 36]{Krattenthaler05} &
  \cite[Theorem 5]{KoutschanThanatipanonda13} & 2013 \\
  $\D11{\mu}n$ & \cite[Problem 34]{Krattenthaler05} &
  \cite[Theorem 13]{KoutschanThanatipanonda19} & 2019 \\
  $\D{2r}0{\mu}n$ & \cite[Theorem 18]{KoutschanThanatipanonda19} &
  \cite[Theorem 18]{KoutschanThanatipanonda19} & 2019 \\
  $\D{2r-1}0{\mu}n$ & \cite[Theorem 19]{KoutschanThanatipanonda19} &
  \cite[Theorem 19]{KoutschanThanatipanonda19} & 2019
\end{tabular}
\caption{Previous work on the determinants $\D st\mu{n}$ and $\E
  st\mu{n}$. What is remarkable about the theorems in the ``resolved'' column
  is that they give reasonably nice closed forms for the corresponding
  determinant. In all cases, the results are valid for $n,r$ being positive
  integers. It will be eventually revealed in this manuscript, that some
  of the $D$ and $E$ families exhibit an interesting symmetrical and
  combinatorial relationship with each other.}
\lbl{tab:olddet}
\end{table}
\vspace{0.5cm}

We now take on the ambitious goal of not only confirming that all previously
unproven conjectures are true, but also highlighting the relationships that we
found between the families that enabled us to accomplish that goal, as well as
the discovery of some new relationships. This work culminates in
\fig{3Dfamilies}. In particular, we give the closed forms of determinants
for four different families. Some of these were ``to do'' from previous
papers, one is simply an easy ``switch'' of the other (see \sect{switch}) and
one has been proposed in this paper as an analog of an old conjecture. We link
to their resolution in \Cref{tab:newdet}.

\begin{table}[ht]
\centering
\begin{tabular}{l|c|c|c}
  Determinant & Condition & First Proposed & Resolved \\
  \hline
  $\E1{2r-1}{\mu}{2m-1}$ & $m\geq r$& \cite[Conjecture 37]{Krattenthaler05} &
  \thm{Krat37ugly} \\
  $\E{2r-1}1{\mu}{2m-1}$ & $m\geq r$& This paper &
  \thm{Krat37nice} \\
  $\D {2r}1{\mu}{2m}$ & $m \geq r$ & \cite[Conjecture 20]{KoutschanThanatipanonda19} &
  \thm{KTConj20} \\
  $\E{-1}{2r-1}{\mu}{2m-1}$ & $m>r$ & This paper &
  \thm{Eneg1CF}\\
  $\D{-1}{2r}{\mu}{2m}$ & $m>r$ & \cite[Conjecture 21]{KoutschanThanatipanonda19} &
  \thm{ktconj21}
\end{tabular}
\caption{Main results of the present paper. The references in the
  ``resolved'' column give a closed form for the corresponding determinant.
  These results are valid for $m,r$ being positive integers under the
  given condition.}
\lbl{tab:newdet}
\end{table}

We remark that much of the ground work to prove the conjectures has already
been laid out
in~\cite{KoutschanThanatipanonda13} and~\cite{KoutschanThanatipanonda19}. Similar
to those papers, we make heavy use of Zeilberger's holonomic
ansatz~\cite{Zeilberger07} (see \sect{holonomicansatz}) and then creative
telescoping~\cite{Zeilberger91}
for proving identities containing symbolic sums that result from the
method. There were three key challenges that we had to overcome in order to be
successful:
{\parskip=0pt
\begin{itemize}
\item The holonomic ansatz could not be applied directly. Certain
  algebraic manipulations had to be invoked to sufficiently simplify our
  matrices before we could apply the ansatz to deduce certain relationships
  between the $E$ and $D$ determinants. Then we still had to use an induction
  argument to arrive at the desired conclusions.
\item In trying to find formulas for ratios of determinants, we sometimes
  encountered the indeterminate form $\frac{0}{0}$. In order to prevent a
  determinant from evaluating to zero, we chose to perturb our parameters $s$
  and~$t$.  Hence, it was not possible to use the classical definition of the
  binomial coefficient over the integers, but we needed to extend the
  definition of the binomial coefficient to the real numbers (see
  \sect{prelim}).
\item Automated symbolic computation was not entirely automatic. We ran into
  many computational bottlenecks, partly due to the extra parameter $r$ in our
  determinants. This is briefly described in the proof of \lem{biglemma1} and
  shown in full detail in the online supplemental material \cite{EM2}. We
  believe that one of the major contributions of this paper is the fact that
  it demonstrates the amazing power of computer algebra to solve combinatorial
  problems, while at the same time reveals limitations in the software.
\end{itemize}

The rest of this paper is mostly organized around the resolution of the
conjectures, but we also include some additional motivation and several new
results. In \sect{prelim}, we introduce all of the important vocabulary,
notations, definitions and properties that we use throughout, and briefly describe the main technique and computational tools. We explain the
combinatorial interpretation for the $E$ determinant in~\sect{comb}. Its
relationship to $D$ (whose interpretation was already described
in~\cite{KoutschanThanatipanonda19}) is shown in~\lem{famA}. The proof
of~\thm{KTConj24} relies heavily on this result. \sect{krat37ktconj20}
highlights the first main event: the proofs of \cite[Conjecture~37]{Krattenthaler05}
and \cite[Conjecture~20]{KoutschanThanatipanonda19}. \sect{ktconj21} highlights
the second main event: the proofs of
\cite[Conjecture~21]{KoutschanThanatipanonda19} and its $E$-analog
(introduced here, not conjectured anywhere
else). In \sect{misc}, we use Andrews' famous
determinant~\cite[Theorem~8]{Andrews79} together with \lem{famA} to prove
\cite[Conjecture~24]{KoutschanThanatipanonda19}.
We are also able to identify two more nice determinant ratios. Finally, in \sect{triangle}, we
conclude with a few relationships that we discovered between certain $E$
determinants (and similarly: $D$ determinants) that do not admit a ``nice''
(i.e., fully factored) closed-form evaluation.
}


\def\ckrat{blue}
\def\cEDnegone{red}
\def\cEDzero{green!75!black!75!white}
\def\cEDzeroref{green}
\def\cDtriA{lime!90!black}
\def\cDtriAref{lime}
\def\cDtriB{brown}
\def\cEtriA{cyan}
\def\cEtriB{yellow!90!black}
\def\cEtriBref{yellow}
\def\cEDcor{magenta}

\newcommand{\nice}[1]{\draw[nice] (#1) ellipse (0.2 and 0.1);}
\newcommand{\ugly}[1]{\draw ($(#1)+(0.14,0.093)$) to ++(-0.346,-0.09) to ++(0.274,-0.097) to ++(0.073,0.188);}

\begin{tikzpicture}[scale=0.7]

\foreach \DEoff in {0,12} {
  \ifthenelse{\DEoff>0}{\newcommand{\DE}{D}}{\newcommand{\DE}{E}}
  \node at (\DEoff+4,15) {$\DE_{s,t}^{\mu}(n)$ Family};
  \foreach \noff in {0,10} {
    \ifthenelse{\DEoff>0}{
      \ifthenelse{\noff>0}{\newcommand{\pn}{o}}{\newcommand{\pn}{e}}
    }{
      \ifthenelse{\noff>0}{\newcommand{\pn}{e}}{\newcommand{\pn}{o}}
    }
    \foreach \s in {-2,...,6} {
      \foreach \t in {-2,...,6} {
        \coordinate (\DE\s\t\pn) at (\DEoff+1.053*\s+0.932*\t,{\noff-0.374*\s+0.618*\t)});
      };
    };
    \foreach \i in {-1,...,5} {
      \draw[-, lightgray!30] (\DE\i-2\pn) to (\DE\i6\pn);
      \draw[-, lightgray!30] (\DE-2\i\pn) to (\DE6\i\pn);
    }
    \filldraw[draw=white,fill=white]
      ($(\DE61\pn)!0.5!(\DE62\pn)$) to ($(\DE16\pn)!0.5!(\DE26\pn)$) to (\DE66\pn);
    \draw[->] (\DE-20\pn) to (\DE60\pn);
    \draw[->] (\DE0-2\pn) to (\DE06\pn);
    \node at ($(\DE60\pn)+(0.4,-0.2)$) {$s$};
    \node at ($(\DE06\pn)+(0.4,0.2)$) {$t$};
  };
};
\node at (E-1-1e) {$n$ even\phantom{m}};
\node at (E-1-1o) {$n$ odd\phantom{m}};
\node at (D44e) {$n$ even};
\node at (D44o) {$n$ odd};

\tikzstyle{nice}=[draw=black,fill=white]
\foreach \s in {0,...,5} {\nice{E\s-1o}\nice{E\s-1e}\nice{D\s-1o}\nice{D\s-1e}};
\foreach \i in {0,2,4} {\nice{E\i0o}\nice{E0\i o}};
\foreach \i in {1,3,5} {\nice{D\i0e}\nice{D0\i e}};
\nice{D6-1e}

\tikzstyle{nice}=[draw=black,fill=black]
\foreach \i in {0,...,5} {\nice{D\i0o}\nice{D0\i o}};
\foreach \i in {0,2,4} {\nice{D\i0e}\nice{D0\i e}};
\nice{E11e}\nice{E11o} 
\nice{E22o} 

\foreach \i in {1,...,6} {\ugly{D-1\i o}\ugly{D\i1o}\ugly{E-1\i e}};
\foreach \i in {2,...,5} {\ugly{D1\i o}\ugly{E1\i e}\ugly{E\i1e}};
\foreach \i in {1,3,5}   {\ugly{D-1\i e}\ugly{D1\i e}\ugly{D\i1e}};
\foreach \i in {2,4}     {\ugly{E1\i o}\ugly{E\i1o}\ugly{E-1\i o}};
\ugly{E-16o}
\ugly{E-10e}\ugly{E-10o}\ugly{D-10o}


\foreach \i in {1,3,5} {
  \pgfmathtruncatemacro{\j}{\i+1}
  \coordinate (TD\i1) at ($(D\i1o)!0.25!(D\i1e)$);
  \draw[-, thick, \cDtriA]
    (D\i1e) to (TD\i1) to (D\i1o) to (D\j1o) to[out=210,in=93] (TD\i1);
};
\node[label={right:\hypersetup{linkcolor=\cDtriA}\cor{triangleD1}}] at ($(D51e)!0.9!(TD51)$) {};

\foreach \i in {1,3,5} {
  \pgfmathtruncatemacro{\j}{\i+1}
  \coordinate (TD-1\i) at ($(D-1\i o)!0.25!(D-1\i e)$);
  \draw[-, thick, \cDtriB]
    (D-1\i e) to (TD-1\i) to (D-1\i o) to (D-1\j o) to[out=240,in=85] (TD-1\i);
};
\node[label={left:\hypersetup{linkcolor=\cDtriB}\cor{triangleD-1}}] at ($(D-11e)!0.75!(TD-11)$) {};


\foreach \i in {2,4} {
  \pgfmathtruncatemacro{\j}{\i+1}
  \coordinate (TE\i1) at ($(E\i1e)!0.25!(E\i1o)$);
  \draw[-, thick, \cEtriA]
    (E\i1o) to (TE\i1) to (E\i1e) to (E\j1e) to[out=210,in=93] (TE\i1);
};
\node[label={right:\hypersetup{linkcolor=\cEtriA}\cor{triangleE1}}] at ($(E41o)!0.8!(TE41)$) {};

\foreach \i in {1,3,5} {
  \pgfmathtruncatemacro{\j}{\i+1}
  \coordinate (TE-1\i) at ($(E-1\j e)!0.25!(E-1\j o)$);
  \draw[-, thick, \cEtriB]
    (E-1\j o) to (TE-1\i) to[out=93,in=330] (E-1\i e) to (E-1\j e) to (TE-1\i);
};
\node[label={left:\hypersetup{linkcolor=\cEtriB}\cor{triangleEneg1}}] at ($(E-12o)!0.6!(TE-11)$) {};

\tikzstyle{angle1}=[out=-15,in=-160]
\tikzstyle{angle2}=[out=-160,in=-15]
\tikzstyle{angle1}=[out=0,in=180]
\tikzstyle{angle2}=[out=185,in=5]
\tikzstyle{connection}=[-,thick, \cEDzero]
\tikzstyle{nice}=[draw=\cEDzero,fill=\cEDzero]
\foreach \i in {0,...,5} {\nice{E\i0e}\nice{E0\i e}};
\foreach \i in {1,3,5}   {\nice{E\i0o}\nice{E0\i o}};
\draw[connection] (E00e) to[angle1] (D10o) to[angle2] (E20e) to[angle1] (D30o) to[angle2] (E40e) to[angle1] node[midway,below left] {\hypersetup{linkcolor=\cEDzero}\lem{famA}} (D50o);
\draw[connection] (D00o) to[angle2] (E10e) to[angle1] (D20o) to[angle2] (E30e) to[angle1] (D40o) to[angle2] (E50e);

\tikzstyle{connection}=[-,thick,\cEDnegone]
\tikzstyle{nice}=[draw=\cEDnegone,fill=\cEDnegone]
\draw[connection]
  (D10o) to[out=130,in=60]
  node[pos=0.15, above left] {\hypersetup{linkcolor=\cEDnegone}\lem{quoED1}}
  (E1-1o) to[out=15,in=220]
  node[pos=0.74,above] {\hypersetup{linkcolor=\cEDnegone}\lem{biglemma2}}
  (D2-1e) to[out=230,in=10]
  (E3-1o) to[out=0,in=215]
  (D4-1e) to[out=225,in=-3]
  (E5-1o) to[out=348,in=215] (D6-1e);
\draw[connection] (E1-1o) to[out=25,in=225] (D0-1e);
\foreach \i in {1,3,5} {
  \nice{E-1\i o}
  \node[above] at (E-1\i o) {\hypersetup{linkcolor=\cEDnegone}\labelcref{thm:Eneg1CF}};
};
\foreach \i in {0,2,4,6} {
  \nice{D-1\i e}
  \node[above] at (D-1\i e) {\hypersetup{linkcolor=\cEDnegone}\labelcref{thm:ktconj21}};
};

\tikzstyle{connection}=[-,thick,\ckrat]
\tikzstyle{nice}=[draw=\ckrat,fill=\ckrat]
\draw[connection]
  (E11o) to[out=25,in=170]
  node[midway,above] {\hypersetup{linkcolor=\ckrat}\lem{biglemma1}}
  (D21e) to[out=180,in=40]
  (E31o) to[out=30,in=180]
  (D41e) to[out=190,in=40]
  (E51o) to[out=30,in=190] (D61e);
\foreach \i in {3,5} {
  \nice{E1\i o}
  \node[above] at (E\i1o) {\hypersetup{linkcolor=\ckrat}\labelcref{thm:Krat37nice}};
  \nice{E\i1o}
  \node[above] at (E1\i o) {\hypersetup{linkcolor=\ckrat}\labelcref{thm:Krat37ugly}};
};
\foreach \i in {2,4,6} {
  \nice{D\i1e}
  \node[above] at (D\i1e) {\hypersetup{linkcolor=\ckrat}\labelcref{thm:KTConj20}};
  \nice{D1\i e}
};

\draw[-,thick,\cEDcor] (E11e) to[out=50,in=130]
  node[pos=0.55,above] {\hypersetup{linkcolor=\cEDcor}\cor{EDCorollary1}} (D01o);
\draw[-,thick,\cEDcor] (E22o) to[out=50,in=130]
  node[midway,above] {\hypersetup{linkcolor=\cEDcor}\cor{EDCorollary2}} (D12e);
  
\end{tikzpicture}\kern-42pt

{
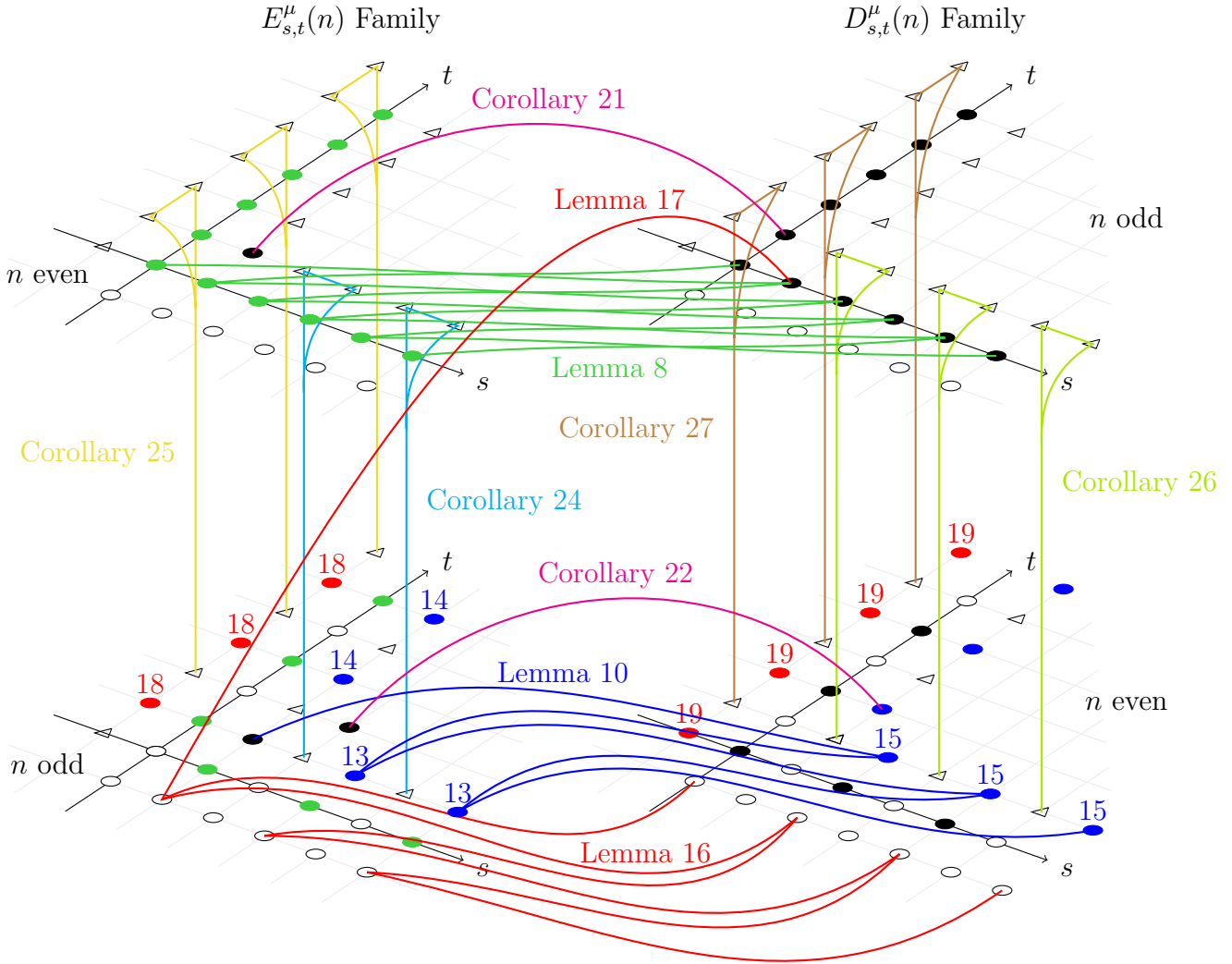
\captionof{figure}
{A $4$-dimensional graphical outline of our contributions: the four
  $(s,t)$-coordinate systems represent the $D$- (resp.~$E$-) determinants, for
  even (resp.\ odd)~$n$. Empty circles refer to zero determinants, filled
  circles to determinants which admit a closed-form product formula, and
  triangles to those which do not. Black circles indicate previously known
  results, while colored circles (together with their corresponding theorem
  number) stand for new results.  Each connection indicates a nice ratio of
  determinants (or limit of a ratio in the case of \lem{biglemma2}). Note that
  \lem{switch} always allows us to derive similar identities with the indices
  $s$ and~$t$ switched. For the sake of clarity, all connections that
  emanate from this symmetry are omitted. For the same reason, the analogous
  connections of \lem{famA} have been omitted in the bottom part of the
  figure.}
  
\lbl{fig:3Dfamilies}
}


\section{Preliminaries}
\lbl{sec:prelim}

This section introduces definitions, notations, vocabulary and properties that
will be used throughout the paper. We also briefly describe the main technique and computational tools that we use to prove some key lemmas. We include them here for the reader's convenience.

\subsection{Pochhammer Symbols and Generalized Binomial Coefficients}
\lbl{sec:gamma}

Many of the formulas in this paper contain rising factorials, which we
represent by the Pochhammer symbol, defined for an indeterminate~$a$, and
$b\in\Z$\,:
\begin{equation*}
(a)_b:=
\begin{cases}
a(a+1)\cdots(a+b-1), & b>0,\\
1, & b=0,\\
\frac{1}{(a+b)_{-b}}, & b<0.
\end{cases}
\end{equation*}
This symbol can also be written as a quotient of gamma functions. We list some
(well-known) properties of the Pochhammer symbol that we use most often
throughout our proofs.
\begin{multicols}{2}
\begin{enumerate}
\item[\namedlabel{itm:pochdef}{(P1)}]\ $(a)_{b}=\frac{\Gamma(a+b)}{\Gamma(a)}$,
\item[\namedlabel{itm:pochinverse}{(P2)}]\ $(a)_{-b} = \frac{1}{(a-b)_b}$,
\item[\namedlabel{itm:pochinterlace}{(P3)}]\ $2^{2b}\cdot(a)_b \cdot \left(a+\frac{1}{2}\right)_b=(2a)_{2b}$,
\item[\namedlabel{itm:pochconnect}{(P4)}]\ $(a)_{b}\cdot (a+b)_{c}=(a)_{b+c}$,
\item[\namedlabel{itm:pochshift}{(P5)}]\ $\frac{(a)_{b}}{(a)_k}=(a+k)_{b-k}$,
\item[\namedlabel{itm:negpoch}{(P6)}]\ $(-a)_b = (-1)^b (a-b+1)_b$,
\item[\namedlabel{itm:pochprod}{(P7)}]\ $\prod_{i=0}^{b-1} (a+i)_k = \prod_{i=0}^{k-1} (a+i)_b$,
\item[\namedlabel{itm:pochprodstep}{(P8)}]\ $\prod_{i=0}^{k-1} (a+ib)_b = (a)_{k b}$.
\end{enumerate}
\end{multicols}

In our work, we find that there is a need to use a more generalized definition
of the binomial coefficient in order to be able to realize our proofs. To be
more specific, for the case $t=-1$, one can see that all entries in the first
column of the matrices $\Mat Ds{-1}{\mu}n$ and $\Mat Es{-1}{\mu}n$ will be
zero, giving us a zero determinant. Since we want to consider ratios of such
determinants, this would result in an indeterminate form $\frac{0}{0}$ rather
than some potentially useful expression. We move away from the offending form
by applying a small perturbation to the parameters and then observing the
ratio's behavior in the limit (see \sect{ktconj21}). Hence, the binomial
coefficients would need to make sense at these perturbations, and for this
purpose, we make great use of the gamma function, which is defined for all
$z\in\C\setminus\lbrace 0, -1, -2, \ldots \rbrace$ in such a way that
\begin{equation}\lbl{eq:gamma}
\Gamma(z+1)=z\Gamma(z).
\end{equation}

\begin{definition}\lbl{def:newbc}
For an indeterminate $x$ and $y\in\C\setminus\{-1,-2,\dots\}$,
we define
\[
  \binom{x}{y} :=  \frac{\Gamma(x+1)}{\Gamma(x-y+1) \, \Gamma(y+1)}.
\]
\end{definition}

Using this definition, we can easily derive a generalization of Pascal's
identity, as well as a useful summation identity.

\begin{lemma}\lbl{lem:pascal}
Let $x$ be an indeterminate and $y\in\C\setminus\{-1,-2,\dots\}$ and
$j\in\N$. Then the following identities hold:
\begin{equation}\lbl{eq:pascal}
\binom{x+1}{y}-\binom{x}{y}=\binom{x}{y-1},
\end{equation}
\begin{equation}\lbl{eq:pascalsum}
\sum\limits_{\ell=0}^{j-1}\binom{x+\ell}{y+\ell}=\binom{x+j}{y+j-1}-\binom{x}{y-1}.
\end{equation}
\end{lemma}

\begin{proof}
The first identity is derived using a direct application of \defn{newbc} to
each binomial coefficient and suitable usages of \eqref{eq:gamma}:
\begin{align*}
  \frac{\Gamma(x+2)}{\Gamma(y+1)\,\Gamma(x-y+2)}-
  \frac{\Gamma(x+1)}{\Gamma(y+1)\,\Gamma(x-y+1)}
  &=\frac{\Gamma(x+2)-\Gamma(x+1)\cdot (x-y+1)}{\Gamma(y+1)\,\Gamma(x-y+2)}\\
  &=\frac{\Gamma(x+1)((x+1)-(x-y+1))}{\Gamma(y+1)\,\Gamma(x-y+2)}\\
  &=\frac{\Gamma(x+1)\cdot y}{\Gamma(y+1)\,\Gamma(x-y+2)}\\
  &=\frac{\Gamma(x+1)}{\Gamma(y)\,\Gamma(x-y+2)}.
\end{align*}
The second identity follows directly by applying \eqref{eq:pascal} $j$ times.
\end{proof}

\subsection{Useful Properties of Determinants}
\lbl{sec:switch}

There are three aspects of our determinants from
\cite{KoutschanThanatipanonda19} that deserve special mention because it will
explain why certain assumptions are made in the statements of our results, and
also why we may choose to omit parts of proofs that are repetitive.

\paragraph{Desnanot--Jacobi--Dodgson Identity (DJD)}
This identity is very useful in some determinant evaluations, particularly
whenever there is a need to establish a link between determinants with
parameters $s$ and $t$ that are closely related (see \sect{triangle}). The
proof of this identity can be found in~\cite{Bressoud99}. We
refer the reader to~\cite{AmdeberhanZeilberger01} for an entertaining
discussion and excellent explanation of its use. To be more precise,
if we let $(m_{i,j})_{i,j\in\Z}$ be a doubly
infinite sequence and $M_{s,t}(n)$ to be the determinant of the $n\times
n$-matrix $(m_{i,j})_{s\leq i < s+n, t\leq j < t+n}$, then
\begin{equation}\lbl{eq:DJD}
M_{s,t}(n)M_{s+1,t+1}(n-2)=M_{s,t}(n-1)M_{s+1,t+1}(n-1)-M_{s+1,t}(n-1)M_{s,t+1}(n-1).
\end{equation}
Visually, one can imagine the corresponding matrices (in gray) like this:
\newcommand{\drawbg}{{\color{magenta}\rlap{\rule[-10pt]{30pt}{30pt}}}}
\begin{center}
  {\lgray \rule[-10pt]{30pt}{30pt}}%
  $\;\times\;$%
  {\lgray \drawbg\rule{4pt}{0pt}\rule[-6pt]{22pt}{22pt}\rule{4pt}{0pt}}%
  $\quad=\quad$%
  {\lgray \drawbg\rule[-6pt]{26pt}{26pt}\rule{4pt}{0pt}}%
  $\;\times\;$%
  {\lgray \drawbg\rule{4pt}{0pt}\rule[-10pt]{26pt}{26pt}}%
  $\quad-\quad$%
  {\lgray \drawbg\rule{4pt}{0pt}\rule[-6pt]{26pt}{26pt}}%
  $\;\times\;$%
  {\lgray \drawbg\rule[-10pt]{26pt}{26pt}}\ \ .
\end{center}

\paragraph{Binomial Determinants without Kronecker Deltas}
For sufficiently small~$n$, the Kronecker deltas will not be present in
our matrices $\Mat Dst{\mu}{n}$ and $\Mat Est{\mu}{n}$ (unless $s=t$).
This simplifies the determinant computations greatly and we state
here without proof, the well-known result.

\begin{proposition}[{\cite[Section~2.3]{Krattenthaler99}},
    {\cite[Proposition~14]{KoutschanThanatipanonda19}}]\lbl{prop:detwithnoKD}
For an indeterminate $\mu$ and $n,s,t\in \Z$ with $t\geq 0$ and $n\geq 1$, we have
\[
  \det\begin{pmatrix}\binom{\mu+i+j+s+t-4}{j+t-1}\end{pmatrix}_{1\leq i, j \leq n} =
  \prod\limits_{i=0}^{t-1}\frac{(\mu+s+i-1)_n}{(i+1)_n}.
\]
\end{proposition}

In the statements of all of our lemmas, theorems, and corollaries we will
henceforth assume that $n$ is sufficiently large so that at least one
Kronecker delta is present in the matrix. For smaller~$n$,
\prop{detwithnoKD} can be used.

\paragraph{The Switching Lemma}
Lastly, we present a generalized version of
\cite[Theorem~17]{KoutschanThanatipanonda19}, where we deduce a relationship
between the determinants that have their indices $s$ and $t$ switched. Therefore,
we usually omit analogous cases in the statements of our results
because it is understood that the ``switching lemma" (\lem{switch}) can be
used to obtain them. To prove this lemma, we need a definition and two smaller
lemmas.

\begin{definition}\lbl{def:vec.uv}
  For two real numbers $s,t \notin \lbrace-1,-2,\ldots\rbrace$ and~$n\in\Z^+$,
  we define two vectors $u_{t,n}:=(u_{t,n,i})_{1\leq i \leq n}$
  and $v_{s,n}:=(v_{s,n,j})_{1\leq j \leq n}$ where
  \begin{align*}
    u_{t,n,i} &:= \frac{\Gamma(\mu+t+i-2)}{\Gamma(\mu+n-3)\,\Gamma(i+t)}, \\
    v_{s,n,j} &:= \frac{\Gamma(\mu+n-3)\,\Gamma(j+s)}{\Gamma(\mu+s+j-2)}.
  \end{align*}
\end{definition}

\begin{lemma}\lbl{lem:uv1}
  Let real numbers $s,t \notin \lbrace-1,-2,\ldots\rbrace$ with $t-s\in \N_0$
  and~$n \in \Z^+$. Then for each integer~$i$ with $1\leq i\leq n+s-t$
  we have
  \[
    u_{t,n,i} \cdot v_{s,n,i+t-s} = 1.
  \]
\end{lemma}

\begin{proof}
  By \defn{vec.uv} and the simple substitution $j\rightarrow i+t-s$, the
  result is immediate.
\end{proof}

\begin{lemma}\lbl{lem:uv2}
  For real numbers $s,t \notin \lbrace-1,-2,\ldots\rbrace$ with $t-s\in \N$
  and~$n \in \Z^+$, we have
\[
\prod_{i=1}^n \bigl(u_{t,n,i} \cdot v_{s,n,i}\bigr) =
\prod_{i=0}^{t-s-1} \frac{\bigl(\mu+i+s-1\bigr)_{\!n}}%
         {\bigl(i+s+1\bigr)_{\!n}}.
\]
\end{lemma}

\begin{proof}
By \defn{vec.uv} and the properties of the Pochhammer symbols, 
\begin{align*}
  \prod_{i=1}^n \bigl(u_{t,n,i} \cdot v_{s,n,i}\bigr)
  &= \prod_{i=1}^n \frac{\Gamma(\mu+t+i-2)\,\Gamma(i+s)}{\Gamma(\mu+s+i-2)\,\Gamma(i+t)} \\
  &\stackrel{\ref{itm:pochdef}}{=}
  \prod_{i=1}^n \frac{\bigl(\mu+s+i-2\bigr)_{t-s}}{\bigl(i+s\bigr)_{t-s}} \\
  &\stackrel{\ref{itm:pochprod}}{=}
  \prod_{i=0}^{t-s-1} \frac{\bigl(\mu+s+i-1\bigr)_{\!n}}{\bigl(i+s+1\bigr)_{\!n}}.
  \qedhere
\end{align*}
\end{proof}

\begin{lemma}\lbl{lem:switch} (Switching)
Let $\Mat Ast{\mu}n$ be either $\Mat Dst{\mu}n$ or $\Mat Est{\mu}n$, and
$\detnew Ast{\mu}n$ its corresponding determinant. For $\mu$ indeterminate,
real numbers $s,t \notin \lbrace-1,-2,\ldots\rbrace$ with $t-s\in\N$ and
$n \in \Z^+$,
\begin{equation}\lbl{eq:switch}
\detnew Ast{\mu}n = \prod_{i=0}^{t-s-1} \frac{\bigl(\mu+s+i-1\bigr)_{\!n}}{\bigl(i+s+1\bigr)_{\!n}}\cdot \detnew Ats{\mu}n.
\end{equation}
\end{lemma}

\begin{proof}
The claimed equality of determinants is a direct consequence (using \lem{uv2})
of the following identity of matrices:
\begin{equation}\lbl{eq:FactorDiag}
\bigl(\Mat Ast{\mu}n \bigr)^{\mathrm{T}} =
\diag(u_{t,n}) \cdot \Mat Ats{\mu}n \cdot \diag(v_{s,n}).
\end{equation}
Hence, the rest of the proof will be dedicated to show it.
The $(i,j)$-entry of the right-hand side of~\eqref{eq:FactorDiag} is equal to
\[
u_{t,n,i} \cdot \binom{\mu+i+j+s+t-4}{j+s-1} \cdot v_{s,n,j} \pm u_{t,n,i}\cdot v_{s,n,j}\cdot\delta_{i+t-s,j}.
\]
Since the $(i,j)$-entry of the left-hand side is equal to
$\binom{\mu+i+j+s+t-4}{i+t-1} \pm \delta_{j+s-t,i}$, \lem{uv1} and the fact
that $\delta_{j+s-t,i}=\delta_{i+t-s,j}(=\delta_{i+t,j+s})$ imply that the
Kronecker delta parts of the both sides are equal. As for the binomial
coefficient parts, by \defn{newbc} and \defn{vec.uv} we have that
\begin{align*}
  & \frac{\Gamma(\mu+t+i-2)}{\Gamma(\mu+n-3)\,\Gamma(i+t)} \cdot
  \frac{\Gamma(\mu+s+t+i+j-3)}{\Gamma(j+s)\,\Gamma(\mu+t+i-2)}\cdot
  \frac{\Gamma(\mu+n-3)\,\Gamma(j+s)}{\Gamma(\mu+s+j-2)} \\
  &= \frac{\Gamma(\mu+s+t+i+j-3)}{\Gamma(i+t)\,\Gamma(\mu+s+j-2)}
  = \binom{\mu+s+t+i+j-4}{t+i-1},
\end{align*}
which implies that~\eqref{eq:FactorDiag} holds and so does the lemma.
\end{proof}

\subsection{The Holonomic Ansatz}
\lbl{sec:holonomicansatz}

We recall here the original formulation of the
\textit{holonomic ansatz}, due to Zeilberger~\cite{Zeilberger07}. The method
is a way to deal with a potentially difficult-to-compute family of
determinants $A(n):=\det(a_{i,j})_{1\leq i,j \leq n}$, where the
dimension~$n\geq1$ is a parameter, and the $a_{i,j}$ form a bivariate
holonomic sequence not depending on~$n$. The method requires $A(n)\neq0$ for
all~$n$, but this fact (provided it is the case) can be established by an
induction argument.  By exploiting the Laplace expansion with respect to the
last row, the determinant can be expressed as
\[
  A(n)=\sum\limits_{k=1}^n a_{n,k}\cdot\Cof nkn,
\]
where $a_{n,k}$ is the $k$-th term in the expansion row and $\Cof nkn$ is the
corresponding cofactor. While $\Cof nkn$ might also be difficult to compute,
the induction hypothesis implies that $\Cof nnn=A(n-1)\neq0$, and hence we
can define
\begin{equation}\lbl{eq:cnk}
  c_{n,k} := \frac{\Cof nkn}{\Cof nnn}.
\end{equation}
For each fixed~$n$, the quantities $(c_{n,1},\ldots, c_{n,n})$ satisfy
the following system of equations:
\begin{equation}\lbl{eq:sysansatz}
\begin{cases}
  c_{n,n}=1,& n\geq1, \\[1ex]
  \sum\limits_{k=1}^{n}a_{\ell,k}\cdot c_{n,k}=0,& 1\leq \ell\leq n-1.
\end{cases}
\end{equation}
The first equation is trivially satisfied by the definition of $c_{n,k}$,
while the second equation corresponds to computing determinants with the row
of expansion replaced by a different one from the same matrix, resulting in
the matrix having two equal rows, giving a zero determinant. By the
induction hypothesis, the system in \eqref{eq:sysansatz} has
full rank, and therefore it has a unique solution.

This view is useful in the sense that, for some fixed $n$ and fixed $k$, we
can compute $c_{n,k}$ by~\eqref{eq:sysansatz}. Then we can use the result of
these computations to make a guess for bivariate recurrences with polynomial
coefficients satisfied by the $c_{n,k}$ (a so-called \emph{holonomic
  description}).  Such recurrences may or may not exist, and in the latter
case, the whole method fails.  In other words, we do not try to work with an
explicit form of the $c_{n,k}$ (which may be hard to find) but instead with an
implicit recursive definition. It remains to prove that the guessed
recurrences define the same bivariate sequence as~\eqref{eq:cnk}, for which it
is sufficient to show that this sequence satisfies~\eqref{eq:sysansatz}
\emph{for all} $n$ and~$\ell$. The machinery that can be employed to do such
confirmations is briefly described in \sect{computeralgebra}.

Now, if we have a conjectured formula~$F(n)$ for the determinant~$A(n)$, then
it suffices to prove
\begin{equation}\label{eq:eq3ansatz}
  \sum\limits_{k=1}^n a_{n,k}\cdot c_{n,k} = \frac{F(n)}{F(n-1)}
\end{equation}
for all~$n\geq2$, again using the machinery described in
\sect{computeralgebra}, to conclude that $A(n)=F(n)$. At the same time, we
complete the induction step by checking that $F(n)\neq0$. If no closed form
$F(n)$ is conjectured, then the method yields a holonomic recurrence for the
quotient $A(n)/A(n-1)$, which can be used for \emph{finding} a closed form (by
solving the recurrence), or for efficiently evaluating $A(n)$, or for studying
its asymptotics as $n\to\infty$.

In \lem{biglemma1}, \lem{biglemma2}, and \lem{quoED1}, the reader will see how
we adapt this elegant idea to give us some of our results, after a few
algebraic manipulations.

\subsection{Computational Machinery for Proving Identities}
\lbl{sec:computeralgebra}

In some proofs we will employ the holonomic machinery, which means that in order to
show that a certain identity is true, we will show that both sides satisfy the same
set of recurrences and have the same (finite number of) initial values. If a
function satisfies a sufficient number of linear recurrences with polynomial
coefficients, we will refer to it as \textit{holonomic}.

It is sometimes easier to translate such notions into an appropriate algebraic
framework, so that we can access computer packages (for our purposes,
we use \cite{Koutschan10b}) that automate the computation of these recurrences
within that framework: we will view linear recurrences as
operators in some (non-commutative) algebra, and we say that a
function \emph{satisfies a recurrence} if the corresponding operator
\emph{annihilates} it. The (infinite) set of all recurrences that a function
satisfies translates to a \emph{left ideal} of annihilating operators in the
algebra. Such an \textit{annihilating ideal} can be finitely presented by some
\emph{generators}, for example in the form of a left Gr\"obner basis.
An identity is correct if we can show that (1)
the annihilating ideals for both sides are equal, or one is a subideal of
the other, and (2) both sides agree on sufficiently many initial values
(their number being determined by the ideals).
We refer the reader to some resources (see for example \cite{Ore33, BeckerWeispfenningKredel93, Koutschan09}) if they are interested
in the algebraic theory behind these computations.

We can remark here that many of the identities that need to be proven with the
computer (see the last part of the proofs of \lem{biglemma1}, \lem{biglemma2}
and \lem{quoED1}) contain sums and products of objects that are holonomic
functions in the parameters. One feature in the theory is that we can start by
computing annihilators/recurrences for single terms or factors and then use
closure properties \cite{Zeilberger90} to deduce a grand recurrence for the
whole expression (i.e., sums and products of holonomic functions are still
holonomic and will therefore satisfy a (different) recurrence with polynomial
coefficients). For symbolic sums, as long as we can confirm that we have
``natural boundaries" (in the sense that the summands evaluate to zero beyond
the stated limits), the method of creative telescoping \cite{Zeilberger91} can
be used with minimal effort via packages that have automated these
computations \cite{Koutschan10b}. For more information on the holonomic systems
approach, we highlight the books and survey papers
\cite{PetkovsekWilfZeilberger96,KauersPaule11,Koutschan13a,Chyzak14}.


\section{Combinatorial Interpretation}
\lbl{sec:comb}

In this section, we would like to provide some additional motivation for
studying these determinants. For an early exposition on the connection between
counting plane partitions and determinants, see the work of Gessel and
Viennot~\cite{GesselViennot85} in the early 1980s.
Krattenthaler~\cite{Krattenthaler06} related the determinant $\D00{\mu}n$
to the enumeration of cyclically symmetric rhombus tilings of a hexagon
with a triangular hole whose size depends on~$\mu$. While many of the ideas
in this section have already been covered in \cite{KoutschanThanatipanonda19},
we show that we can apply them to the determinants with the negative Kronecker
delta. Moreover, we are able to present a combinatorial connection between the
$D$ and $E$ determinants (see \lem{famA} at the end of this section). For our
convenience, we will use the same naming conventions described in \sect{intro}
with letters in plain math text ($B$) being the determinant corresponding to
the matrix written in calligraphic text ($\matnot{B}$).

First, we rewrite the determinant $\E st{\mu}n$ as a sum of minors from
expanding along the $i$-th row and pulling out the single cofactor containing
$-1$ from $-\delta_{i+s-t,j}$ to get
\[
E_{s,t}^{\mu}(n)=(-1)^{s-t+1}\cdot M_{i+s-t}^i+\sum\limits_{j=1}^n(-1)^{i+j}\cdot b_{i,j}\cdot M_j^i,
\]
where
\[
b_{i,j}:=\binom{\mu+i+j+s+t-4}{j+t-1},
\]
and $(-1)^{i+j}\cdot M_j^i$ denotes the $(i,j)$-cofactor of $\Mat Est{\mu}n$. We can apply the removal of the Kronecker delta recursively so that what remains are determinants that do not contain any Kronecker deltas, but that are minors of $\Mat Bst{\mu}n :=\left(b_{i,j}\right)_{1\leq i,j\leq n}.$
This results in another formulation of our $E$ determinant (assuming $s\geq t$), that is,
\begin{align}\lbl{eq:sumofminors}
E_{s,t}^{\mu}(n)=\sum\limits_{I\subseteq \lbrace 1,\ldots, n-(s-t)\rbrace} (-1)^{(s-t+1)\cdot|I|}\cdot B_{I+s-t}^I,
\end{align}
where we are summing over all subsets of rows with a nonzero Kronecker delta (producing additional factors of $-1$ each time) and $\Matminor B{I+s-t}I$ is the submatrix obtained by deleting all rows with indices in $I$ and all columns with indices in $I+s-t=\lbrace i+s-t \mid i\in I\rbrace$ from $\Mat Bst{\mu}n$. The formulation for $s\leq t$ is analogous in that we first switch the subsets $I, I+s-t$ and then switch $s,t$ throughout on the right side of \eqref{eq:sumofminors}. 

Using the Lindstr\"om–Gessel–Viennot lemma~\cite{Andrews79, GesselViennot85, Lindstroem73}, we can deduce that $\B st{\mu}n$ counts $n$-tuples of non-intersecting paths in the integer lattice $\N^2$.  Each $b_{i,j}$ counts the number of paths that start at $(\mu+s+i-3,0)$ and end at $(0,t+j-1)$ with step set $\lbrace \leftarrow, \uparrow \rbrace$ in the first quadrant of the $(i,j)$-plane under the assumption that $\mu+s\geq 2$. These non-intersecting lattice paths are in bijection with rhombus tilings of a lozenge-shaped region, where two of the tile orientations (for example, \rotatebox[origin=c]{90}{$\lozenge$} and \rotatebox[origin=c]{330}{$\lozenge$})
correspond to the paths and the third tile orientation to empty locations (for example, \rotatebox[origin=c]{30}{$\lozenge$}). The start and end points are represented by half-rhombi (i.e., triangles) along the southern (bottom) and western (left) boundaries. We illustrate this with a simple example in \fig{tilings}.

\begin{figure}[ht]
\centering
\begin{tabular}{c|cc}
\multirow{4}{*}{$B^{\varnothing}_{\varnothing}=6$}& \includegraphics[scale=0.25]{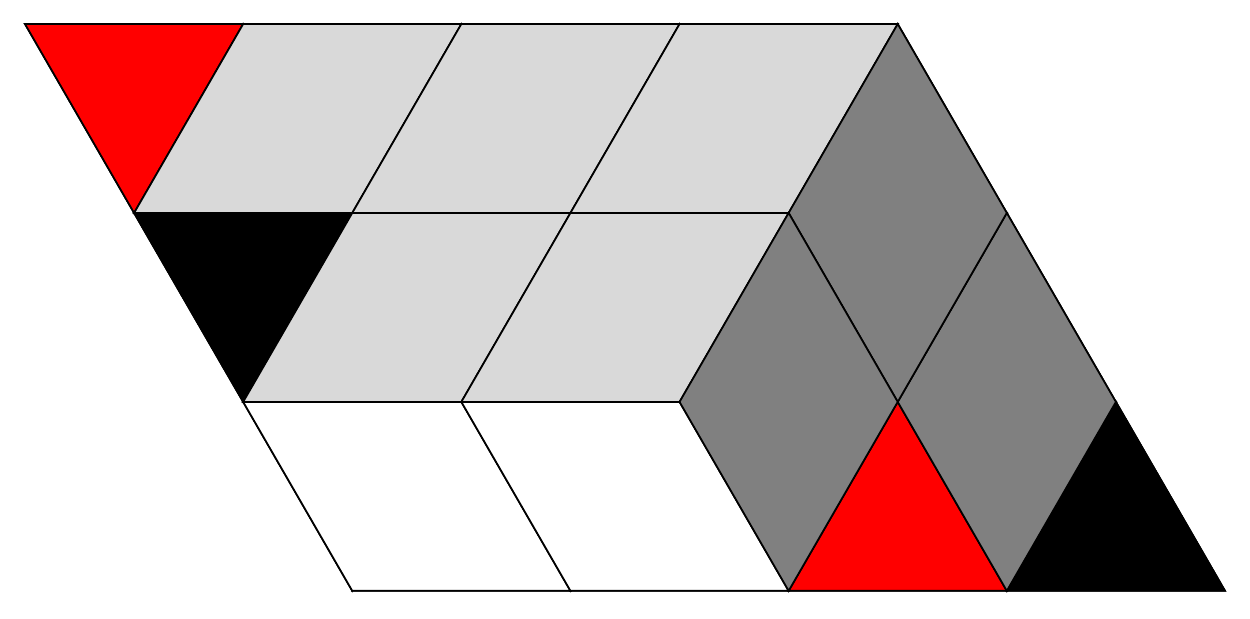}  & \includegraphics[scale=0.25]{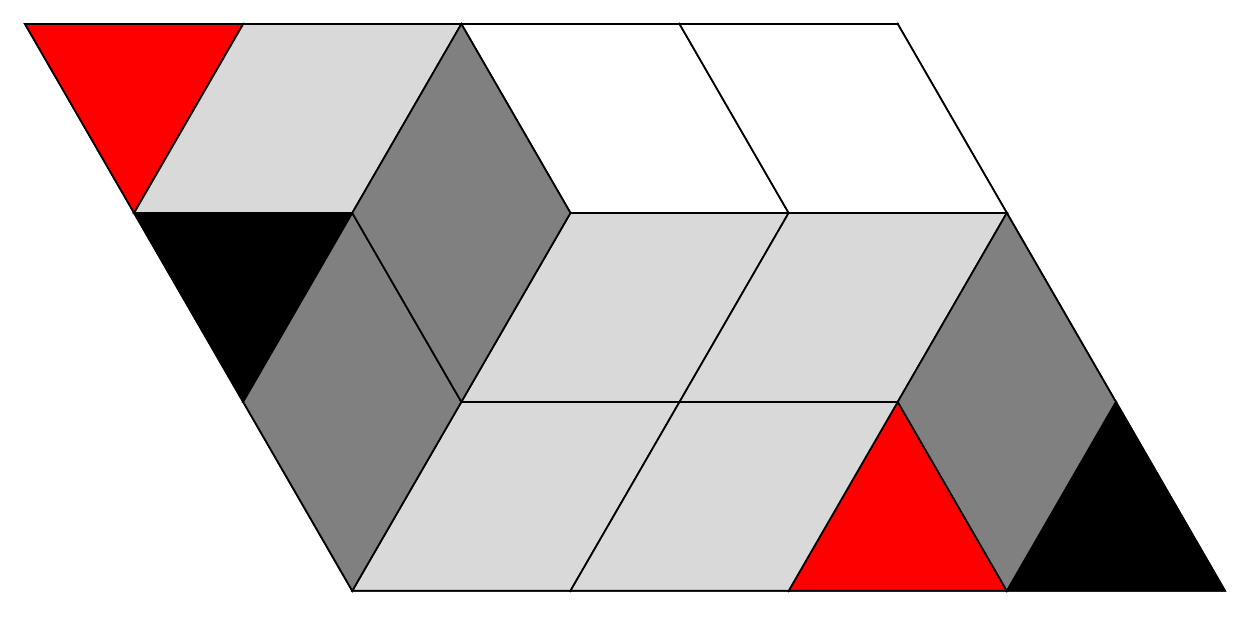} \\
& \includegraphics[scale=0.25]{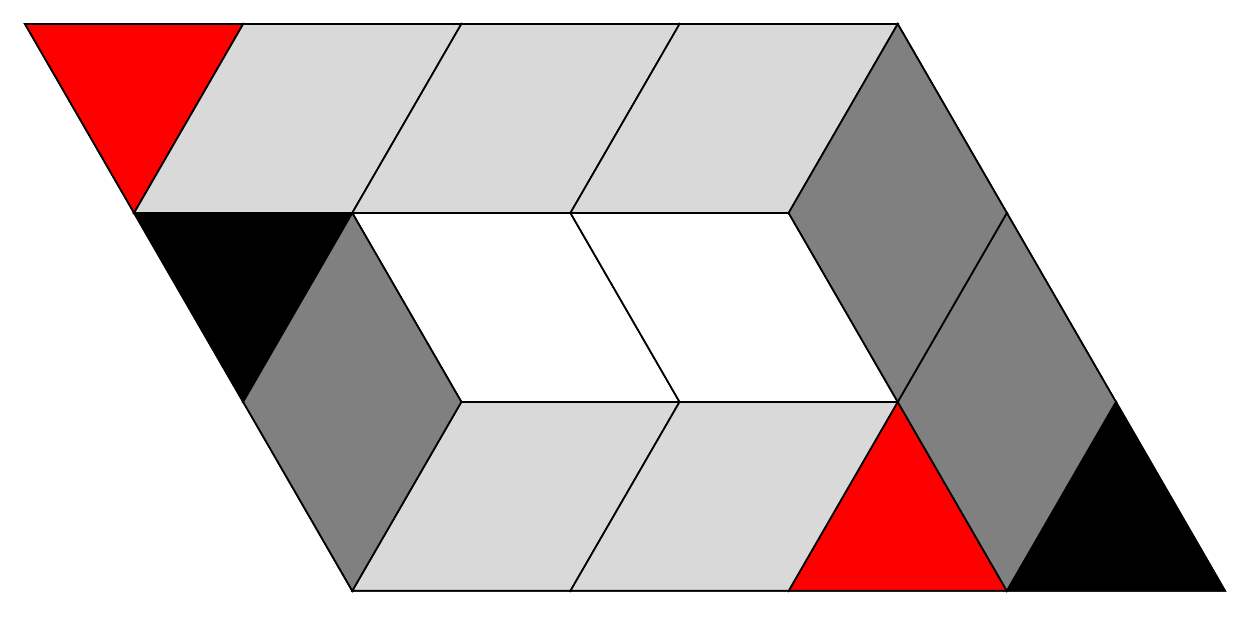} & \includegraphics[scale=0.25]{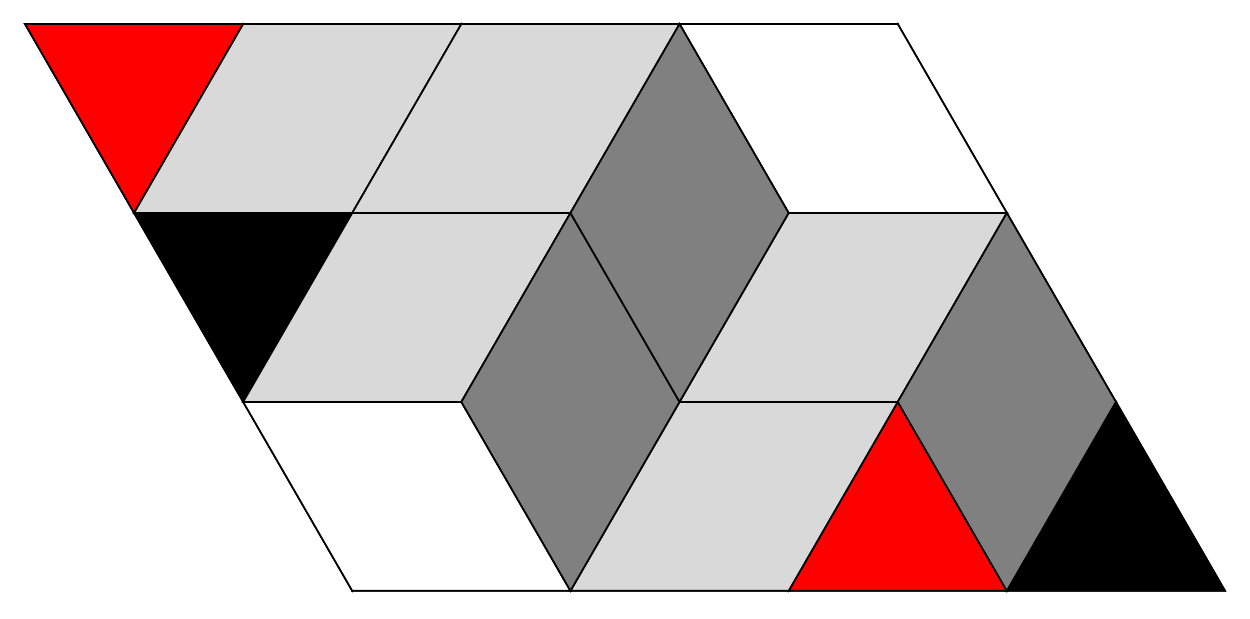} \\
& \includegraphics[scale=0.25]{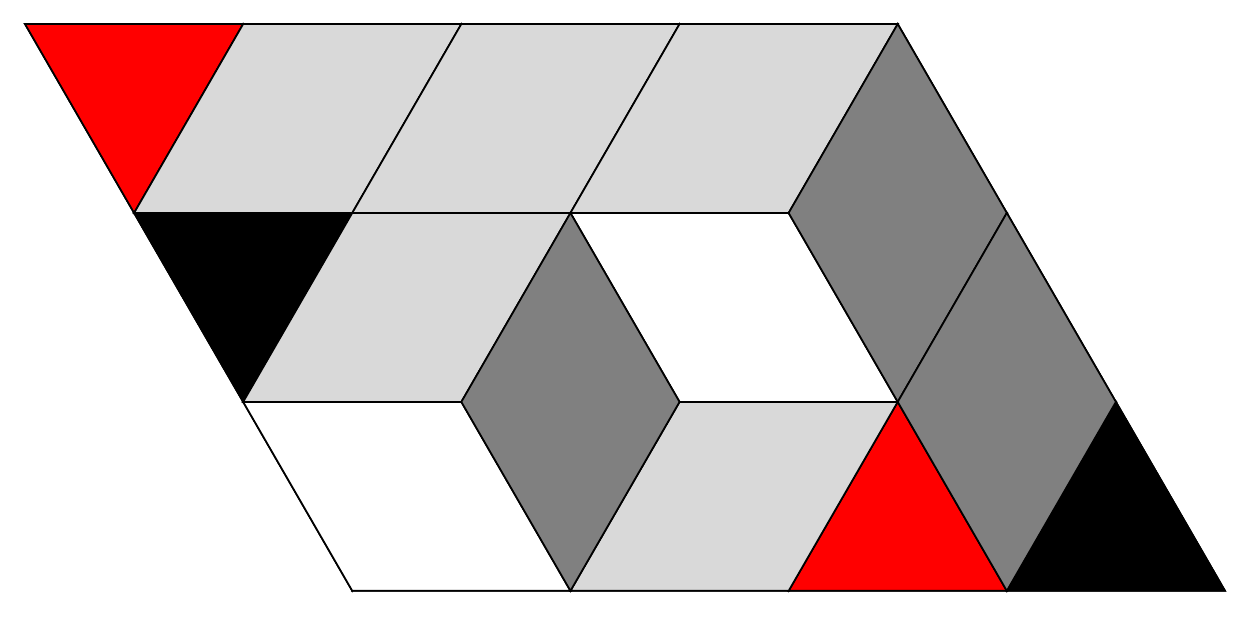} & \includegraphics[scale=0.25]{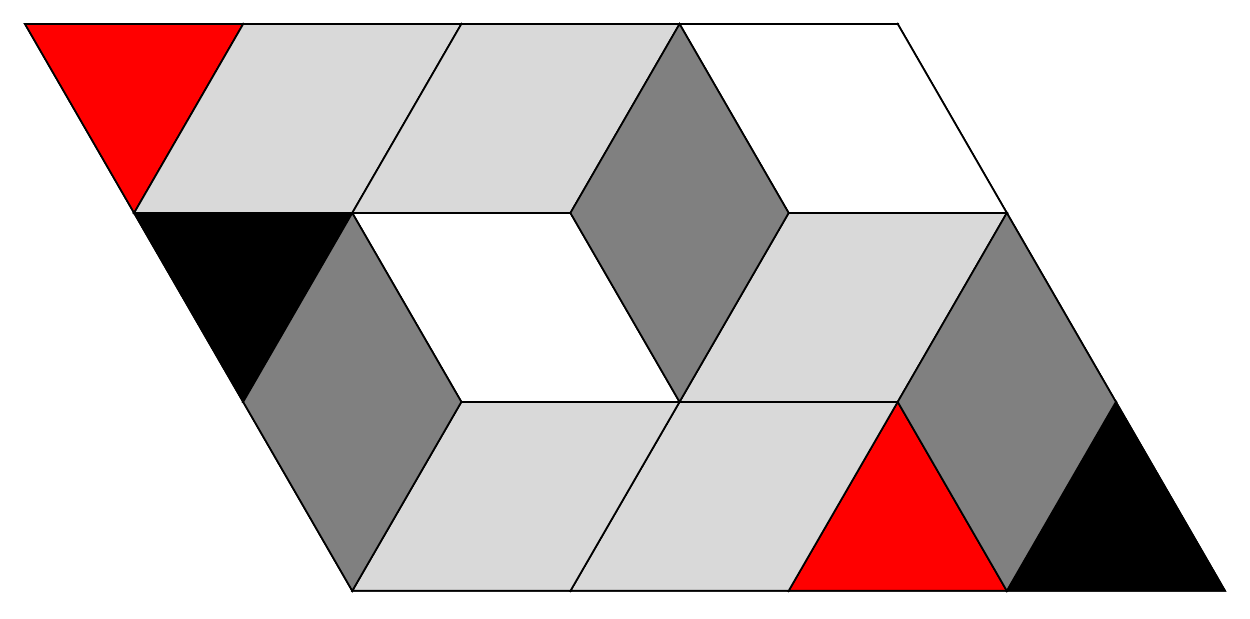} \\[2pt]
\hline\\[-10pt]
\multirow{2}{*}{$B^{\lbrace 1\rbrace}_{\lbrace 2 \rbrace}=4$}& \includegraphics[scale=0.25]{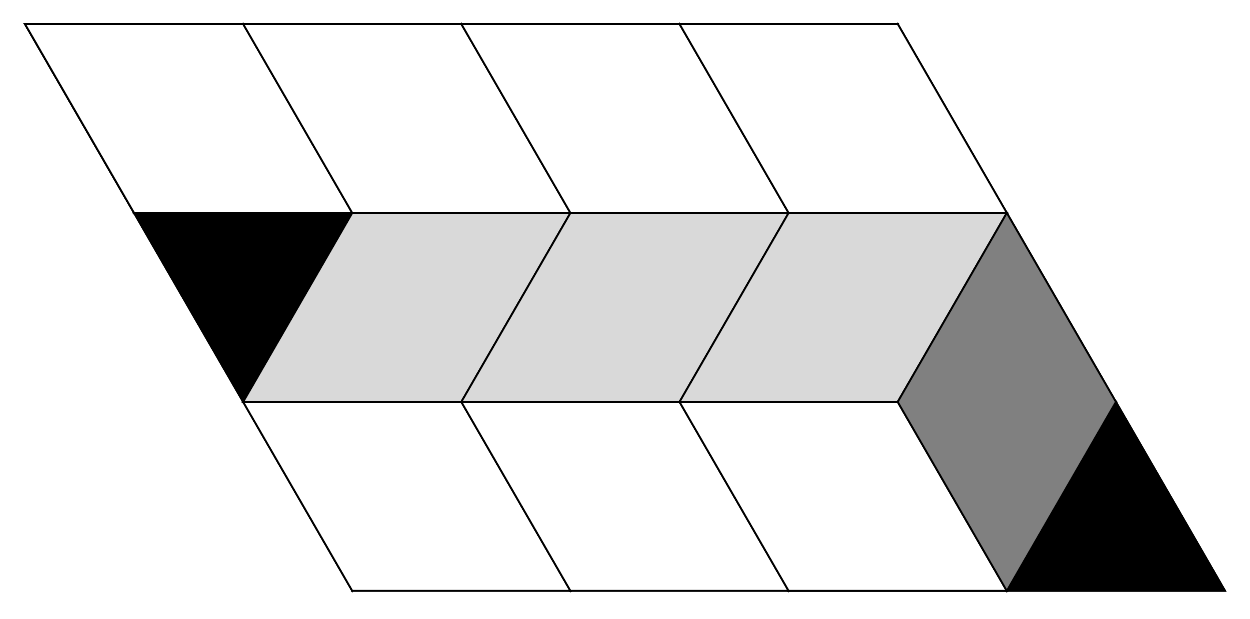} & \includegraphics[scale=0.25]{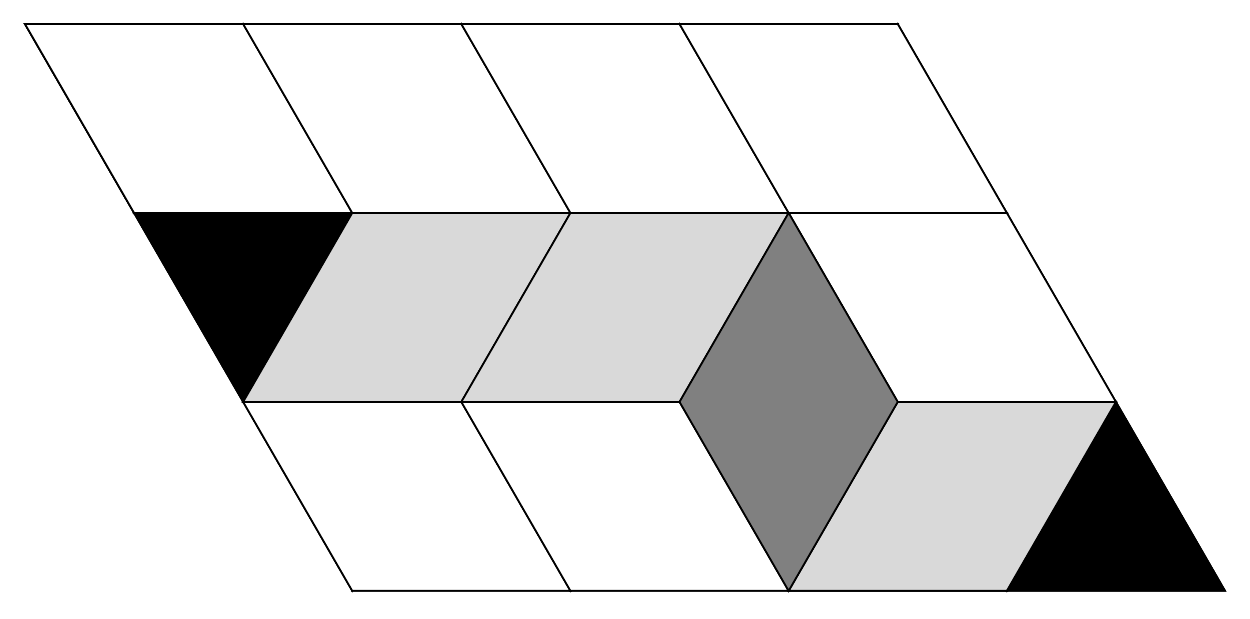} \\
& \includegraphics[scale=0.25]{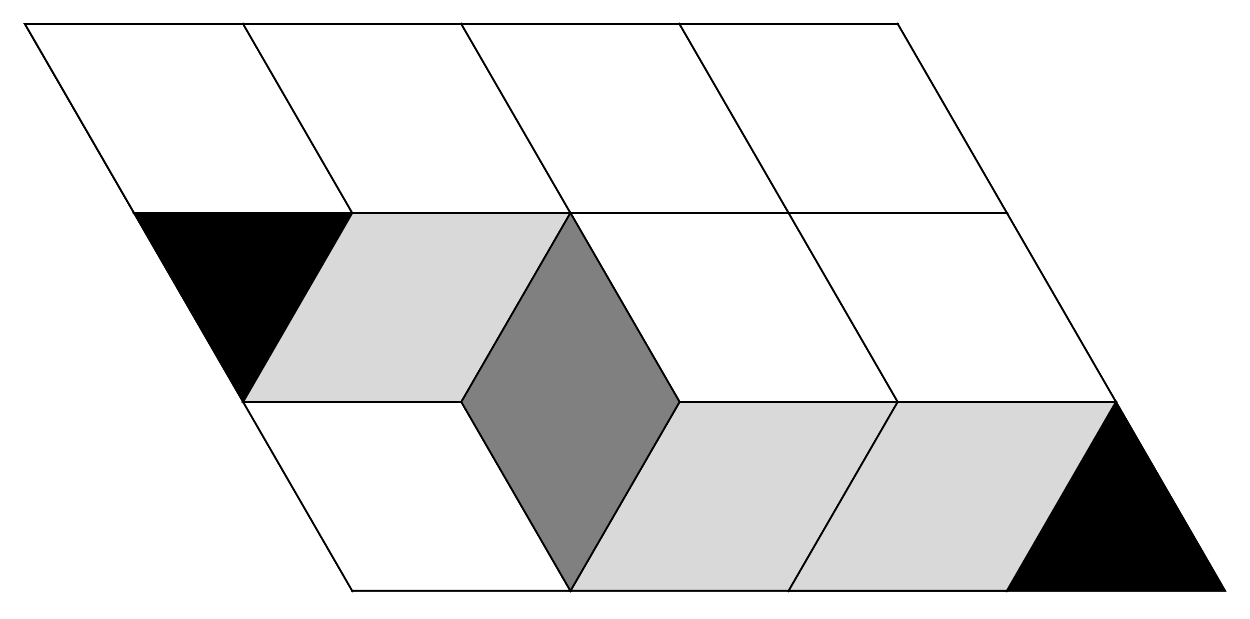} & \includegraphics[scale=0.25]{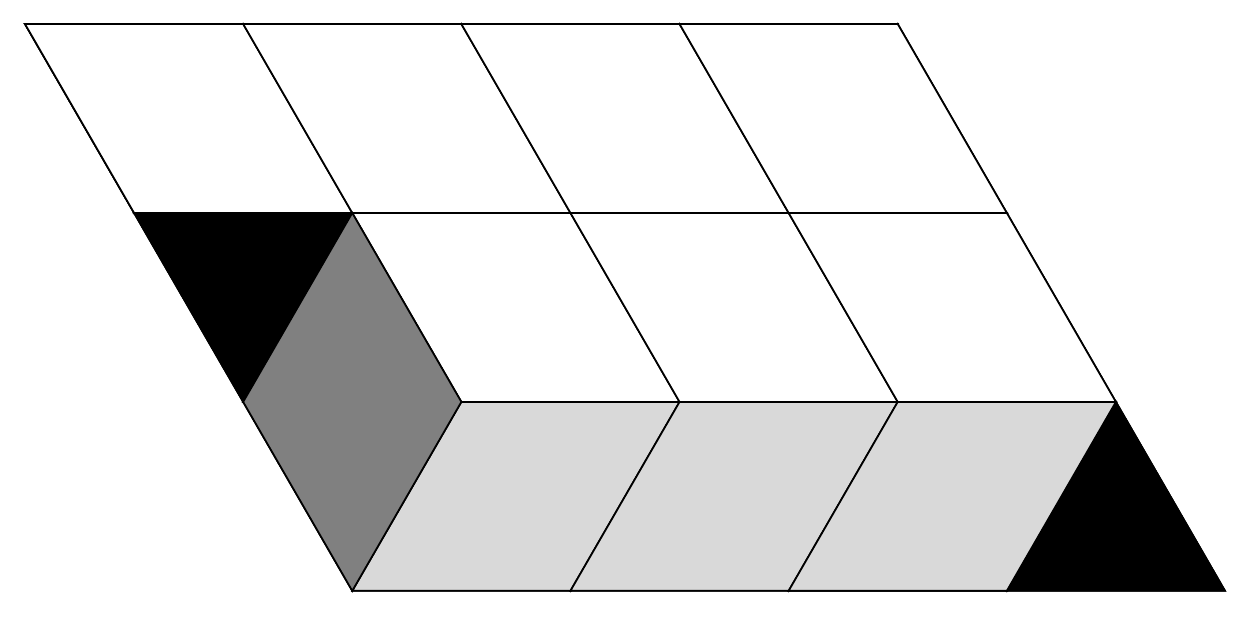} \\
\end{tabular}
\caption{Here is an explicit graphical computation to show that
  $E_{2,1}^{2}(2)=10$. The non-intersecting lattice paths are indicated by the
  light gray tiles (left steps) and dark gray tiles (up steps) with the
  triangles indicating start points at the bottom edge (ordered from left to right
  according to the rows they correspond to) and end points on the left edge
  (ordered from bottom to top according to the columns they correspond
  to). For $I=\varnothing$ and $I+s-t=\varnothing$, all points are
  present. For $I=\lbrace 1\rbrace$ and $I+s-t=I+1=\lbrace 2\rbrace$, we consider the tiling problem of the lozenge-shaped region without the left-most and top-most triangles corresponding to the Kronecker delta that was present in row~1 and column~2 of $\Mat E2122$ (this is indicated by the absence of red triangles in the bottom part of
  the table). The white tiles correspond to locations that are not visited by
  a path.}  \lbl{fig:tilings}
\end{figure}

If $|I|=|J|$, then $B_J^I$ is the determinant of the submatrix after removing
all rows with indices in $I$ and all columns with indices in $J$ and counts
the $(n-|I|)$-tuples of non-intersecting paths where the start points with
indices in $I$ and end points with indices in $J$ are omitted. Since these
points can be omitted or not, depending on the subset $I$ of the set $\lbrace
1,\ldots,n-(s-t)\rbrace$, the elements of the latter (the superset) are called
\textit{optional} points. On the other hand, we always keep certain rows and
columns (respectively, certain starting points and certain ending points of
paths). In particular, we keep the ones that do not contain the Kronecker
delta. We will call the corresponding points \textit{mandatory}. In
\fig{tilings}, the red triangles indicate the activation of some optional points while the black
triangles indicate mandatory points. We can see that these black triangles are
present in every tiling that we count, while the red triangles appear only in
the computation of $B_{\varnothing}^{\varnothing}$, since all rows and columns
are present in the computation. Thus, we can see that each determinant
$B^I_{I+s-t}$ in~\eqref{eq:sumofminors} counts the number of paths with start
and end points that are controlled by the set $I$, and the sign simply acts as
a weight.

We now have enough information to give a combinatorial interpretation of $E^{\mu}_{s,t}(n)$, which is a combination of the determinants~$B^I_{I+s-t}$ and signs. Put together, what does this combination count, and furthermore, what role does the sign play? 

Let us imagine a new counting problem that involves counting the tilings of
not one, but three copies of the same lozenge arranged in a cyclic fashion
(i.e., two are rotations of the first by 120 and 240 degrees,
respectively). For illustration purposes, we will make the following
assumptions: $\mu,s,t\in\Z$ such that $\mu+s\geq 2$, $s\geq t \geq 0$
and $n\geq s$, with the understanding that the case $t\geq s\geq 0$ and
$n\geq t$ is analogous. We remark here that also some cases $t<0$ and $s<0$
have a combinatorial interpretation which will be shown and used in
\sect{triangle}, but to simplify our explanations we exclude such cases
here. The arrangement is such that the optional starting points of one lozenge
are paired with the optional ending points of the other. The mandatory points
remain unpaired. The resulting region is a hexagon (if $s=t$) or a pinwheel
(if $s\neq t$) with a triangular hole of length $\mu-2$ in its center.

The pinwheel-shaped region can actually be viewed as a hexagon if we remove the three triangular regions that emanate from the half-rhombi corresponding to those mandatory points (i.e., the parts that are ``sticking out'' in the pinwheel): these regions can be tiled in only one way (see \fig{forcedtiling}) and removing them does not affect the final count. In this sense, we can almost always achieve a hexagonal region, with the exception being a big triangular region if $1\leq n\leq s, t=0$ (or $1\leq n\leq t, s=0$ in the analogous case). Triangular regions corresponding to the interior mandatory points can be similarly removed, resulting in three additional triangular holes surrounding the original one, further justifying the name ``holey'' (see \fig{regionexamples}).

\begin{figure}[ht]
\centering
\begin{tabular}{c@{\qquad}c@{\qquad}c}
  \includegraphics[height=80pt]{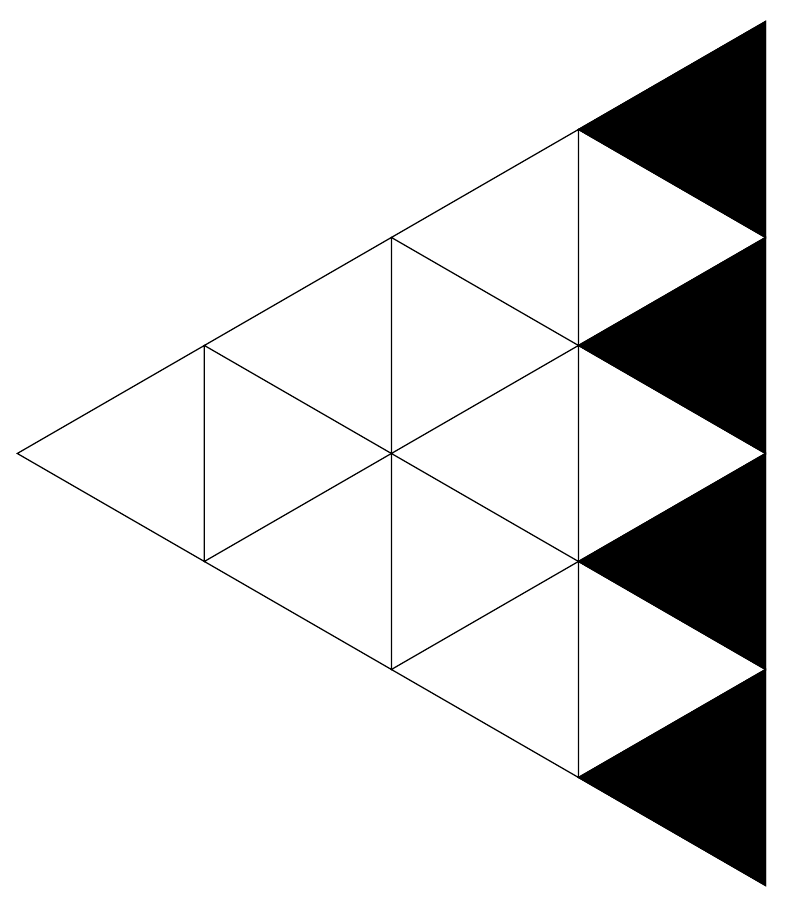}
  & \raisebox{40pt}{$\xrightarrow{\displaystyle\text{\ forced tiling\ }}$} &
  \includegraphics[height=80pt]{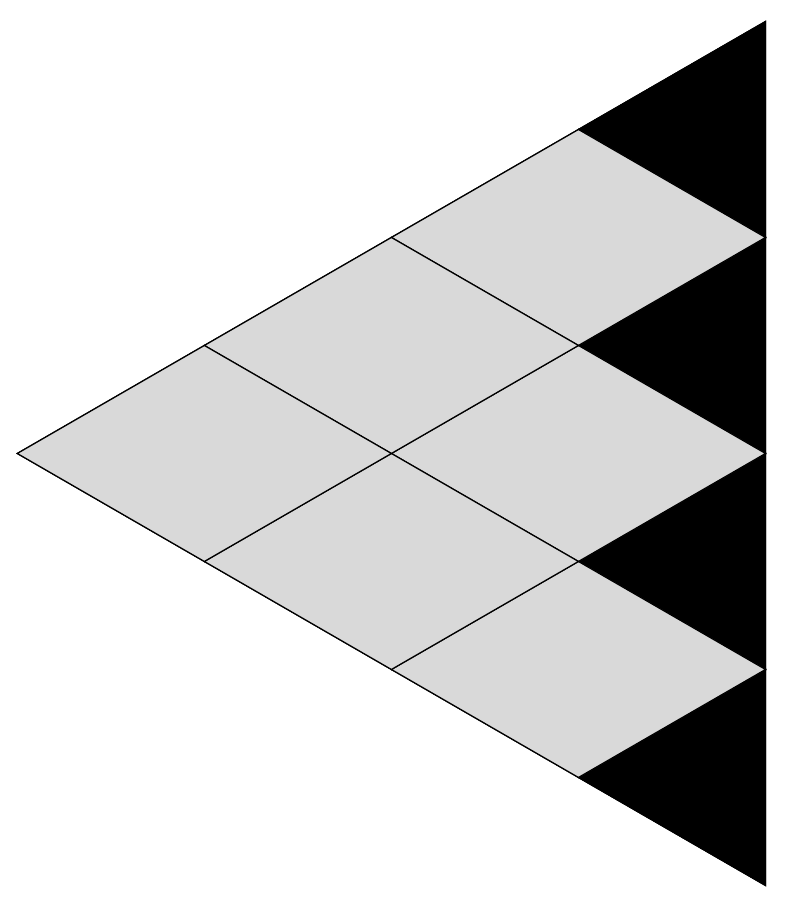}
\end{tabular}
\caption{Example of a forced tiling of a triangular region emanating from mandatory points (represented by the smaller black triangles). Since there is only one way to tile such a region, removing it will not influence the final count.
}
\lbl{fig:forcedtiling}
\end{figure}

Next, we apply a very important rule: we say that we only want to count cyclically symmetric tilings of this holey hexagon. This ensures that we only count tilings that match the tilings from one of the triplicated lozenges. We can make a few observations by imposing this new condition. First, a full tile is allowed at the optional point connection (and it will appear depending on the tiling that we are considering). Second, the mandatory points will not have a counterpart on the other side of the border and this provides a natural perimeter to prevent the counting of tilings that do not fit with our problem. Third, the space between the vertices of the central triangular hole and the mandatory points exists if $t>0$ (or $s>0$ in the analogous case). A full tile will never cross the border here. \fig{examplewithpaths} illustrates a region to be tiled and depicts a tiling with this rule applied.

\begin{figure}[ht]
\centering
\includegraphics[width=0.45\textwidth]{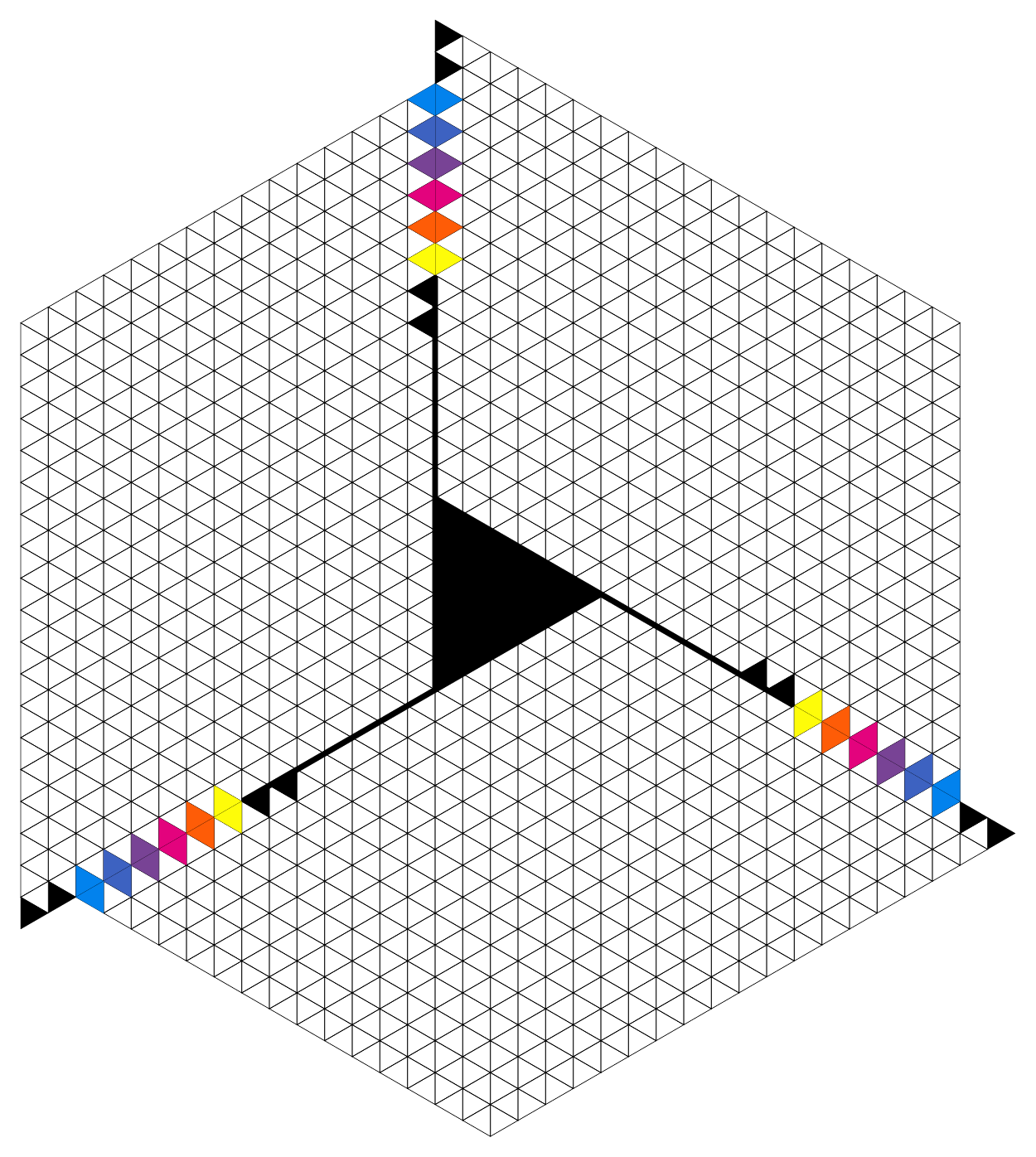}
\qquad
\includegraphics[width=0.45\textwidth]{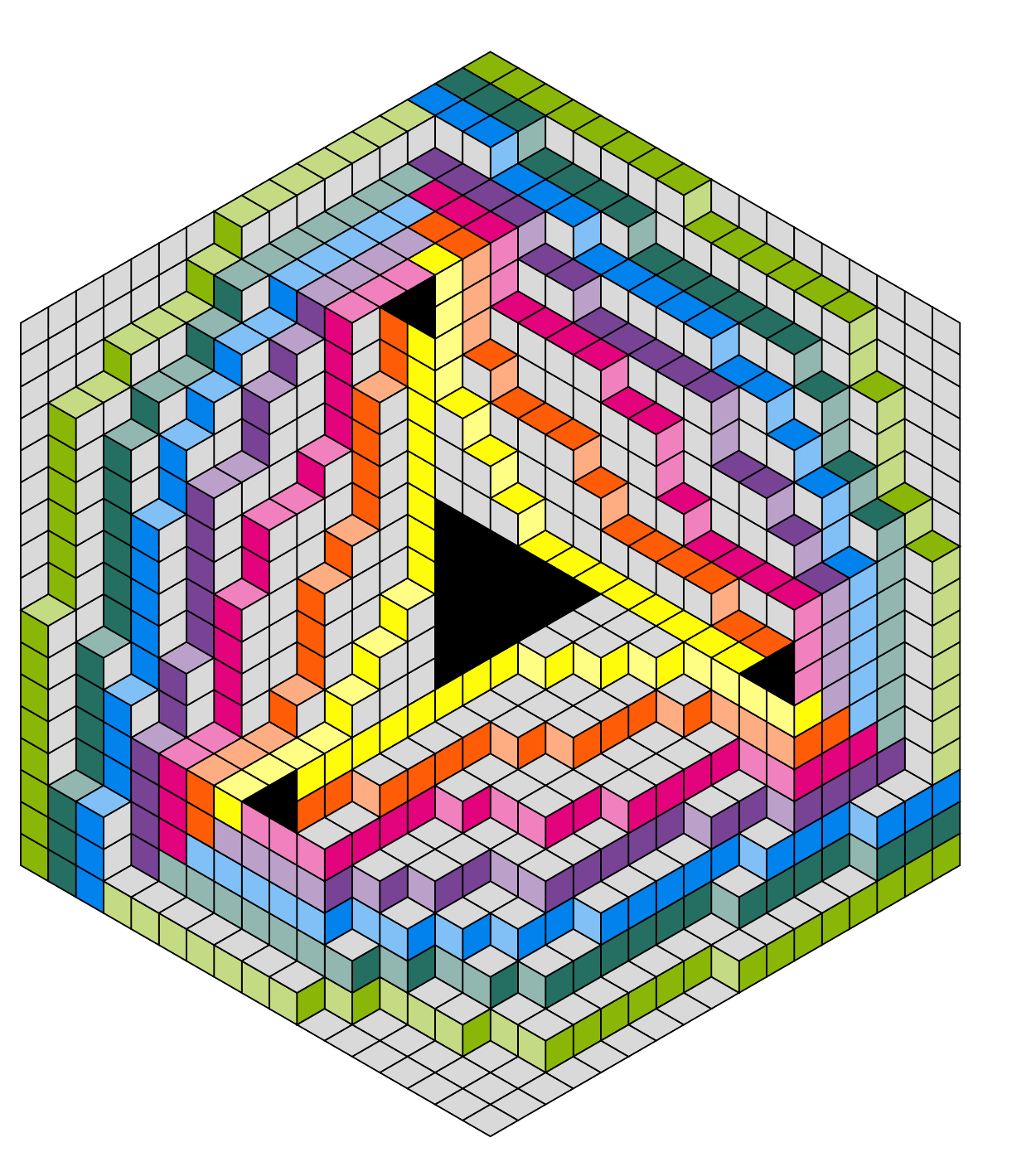}
\caption{The region on the left corresponds to the parameters
  $(s,t,n,\mu)=(5,7,8,8).$ On the right, we illustrate one cyclically
  symmetric tiling of this region. In this example, the optional starting
  points $1,2,3,4,6$ (corresponding to the colors yellow, orange, red, purple,
  and light blue) are ``activated'' while from number~$5$ no path is
  emerging.}
\lbl{fig:examplewithpaths}
\end{figure}

We summarize how we do this region construction more concretely in terms of our parameters under the given assumptions (with periodic commentary in brackets to indicate the analogous case), and provide some examples in \fig{regionexamples}. Set $\Delta=n-(s-t)$ to be the number of Kronecker deltas present in the matrix. We begin with a lozenge of size $n \times (\mu+s+n-2)$ with the longer edge on the bottom and shorter edge on the left (this is reversed in the analogous case, as in \fig{examplewithpaths} and the right side of \fig{hexfortriangle}). Then, we divide the bottom edge into five parts of the following lengths:
    \[
    \underbrace{\mu-2}_{\text{hole}}\quad+\underbrace{\vphantom{\mu}t}_{\text{border line}}+\underbrace{\vphantom{\mu}s-t}_{\text{no start points}}+\underbrace{\vphantom{\mu}\Delta}_{\text{optional start points}}+\underbrace{\vphantom{\mu}s-t}_{\text{mandatory start points}},
    \]
    and divide the left edge into three parts of the following lengths:
    \[
    \underbrace{t}_{\text{border line}}+\underbrace{s-t}_{\text{mandatory end points}}+\underbrace{\Delta}_{\text{optional end points}}.
    \]
Here, the start and end points refer to the start and end points of the paths we want to count and are represented by half-rhombi (i.e., triangles) rather than full tiles. These start points (ordered from left to right along the bottom edge) and end points (ordered from bottom to top along the side edge) are in one-to-one correspondence with the rows and columns of the original matrix, respectively. Their presence or absence is triggered by the set $I$ (see \fig{tilings}).

We proceed by copying/pasting this lozenge twice (so now there are three total), and then rotating the copied lozenges by $120^{\circ}$ and $240^{\circ}$, respectively. Next, we glue them together exactly at the positions corresponding to~$\Delta$. These paired points are indicated by the colored tiles. Furthermore, we apply a thickened border line of length $t$ on all edges starting from one vertex of the central triangular hole to the first triangle corresponding to the closest mandatory point. This is to prevent tiles from spilling over at these connections (there are no start/end points here). Lastly, we remove all forced tilings as described in \fig{forcedtiling} (the bottom row of \fig{regionexamples} illustrates the corresponding regions after such a removal). Our sum of minors formula can now be interpreted in three different ways:
    
\begin{figure}[ht]
\centering
\begin{tabular}{c|c|c|c}
\includegraphics[width=0.23\textwidth]{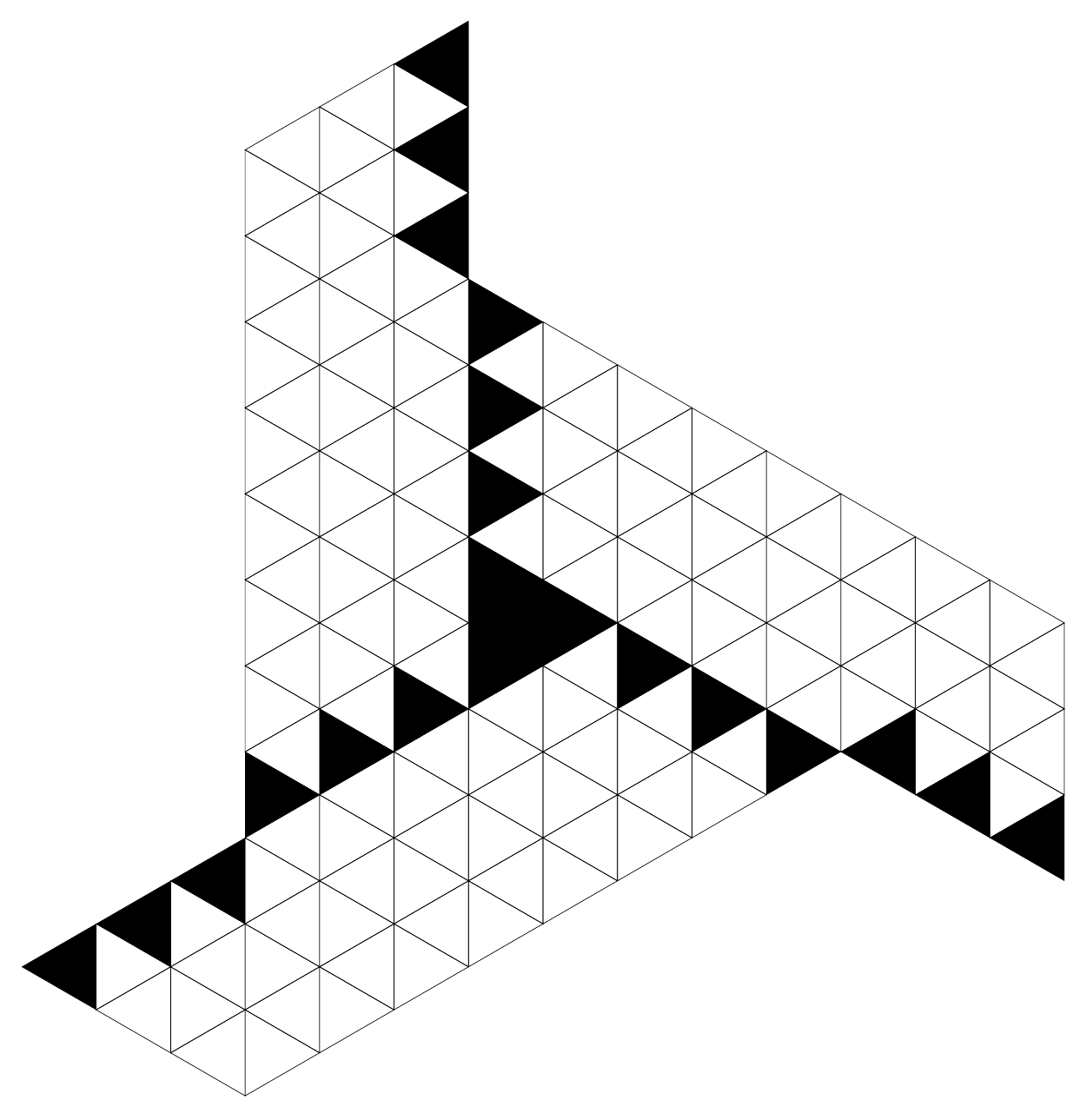}&
\includegraphics[width=0.23\textwidth]{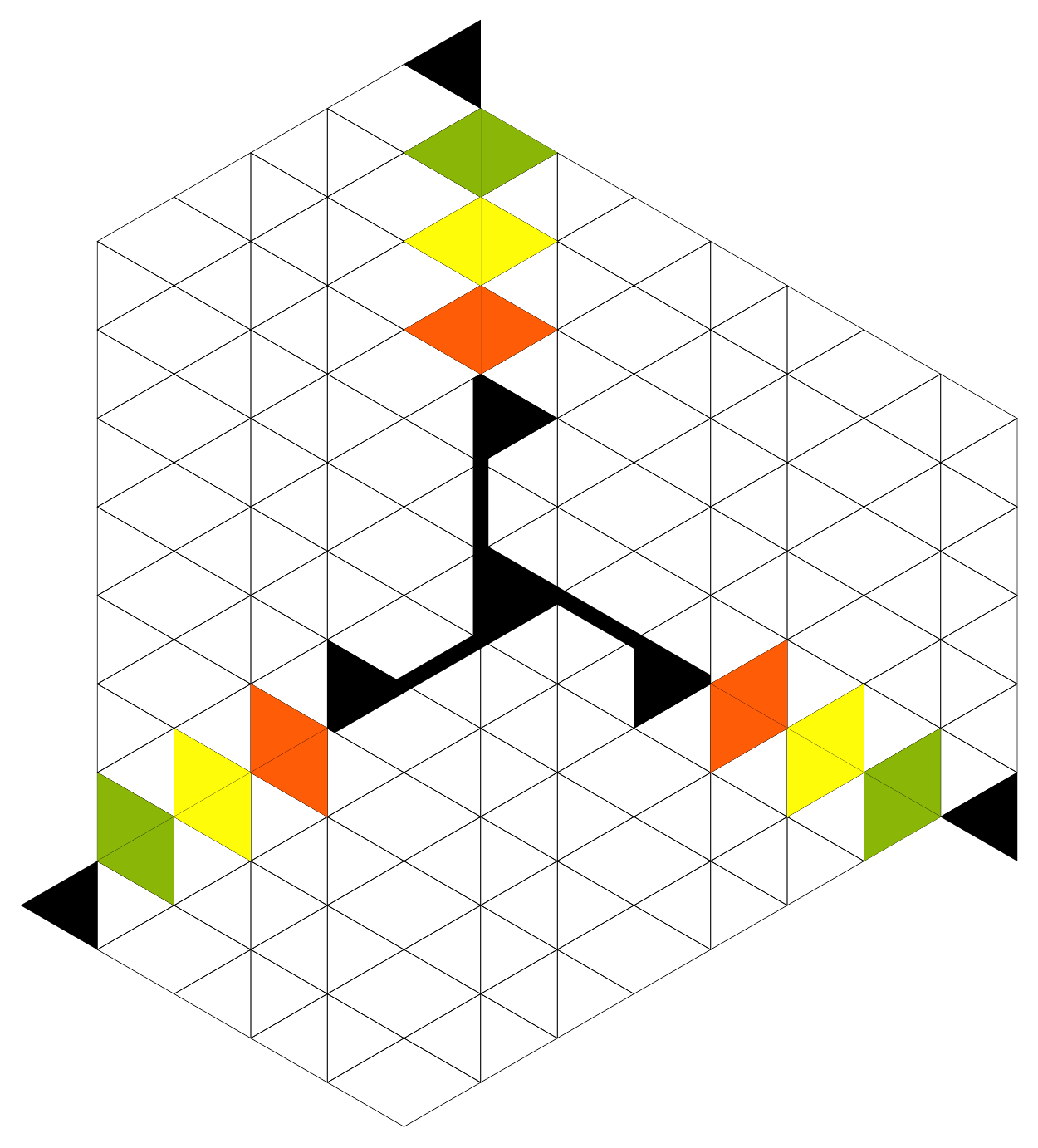}&
\includegraphics[width=0.23\textwidth]{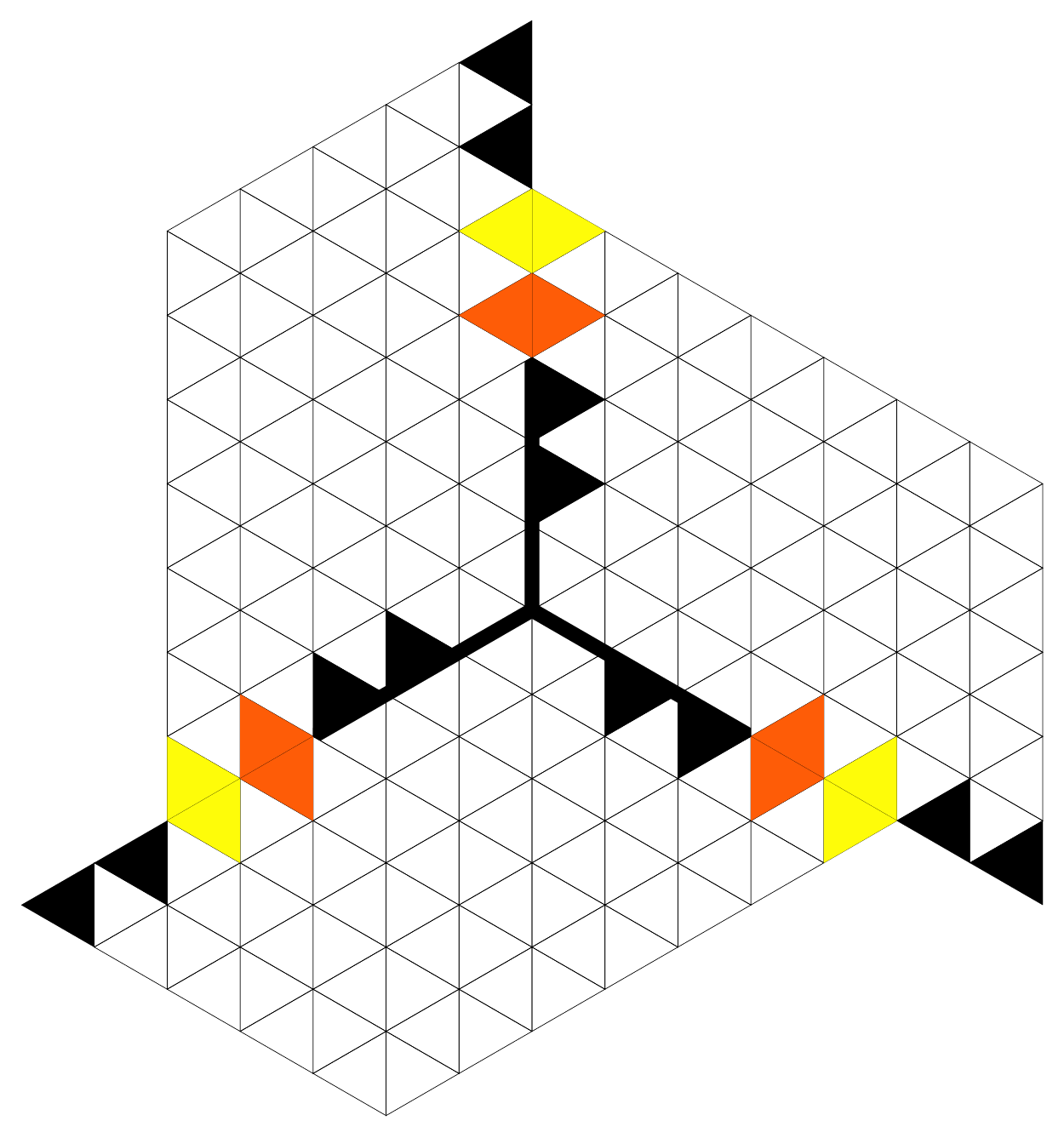}&
\includegraphics[width=0.23\textwidth]{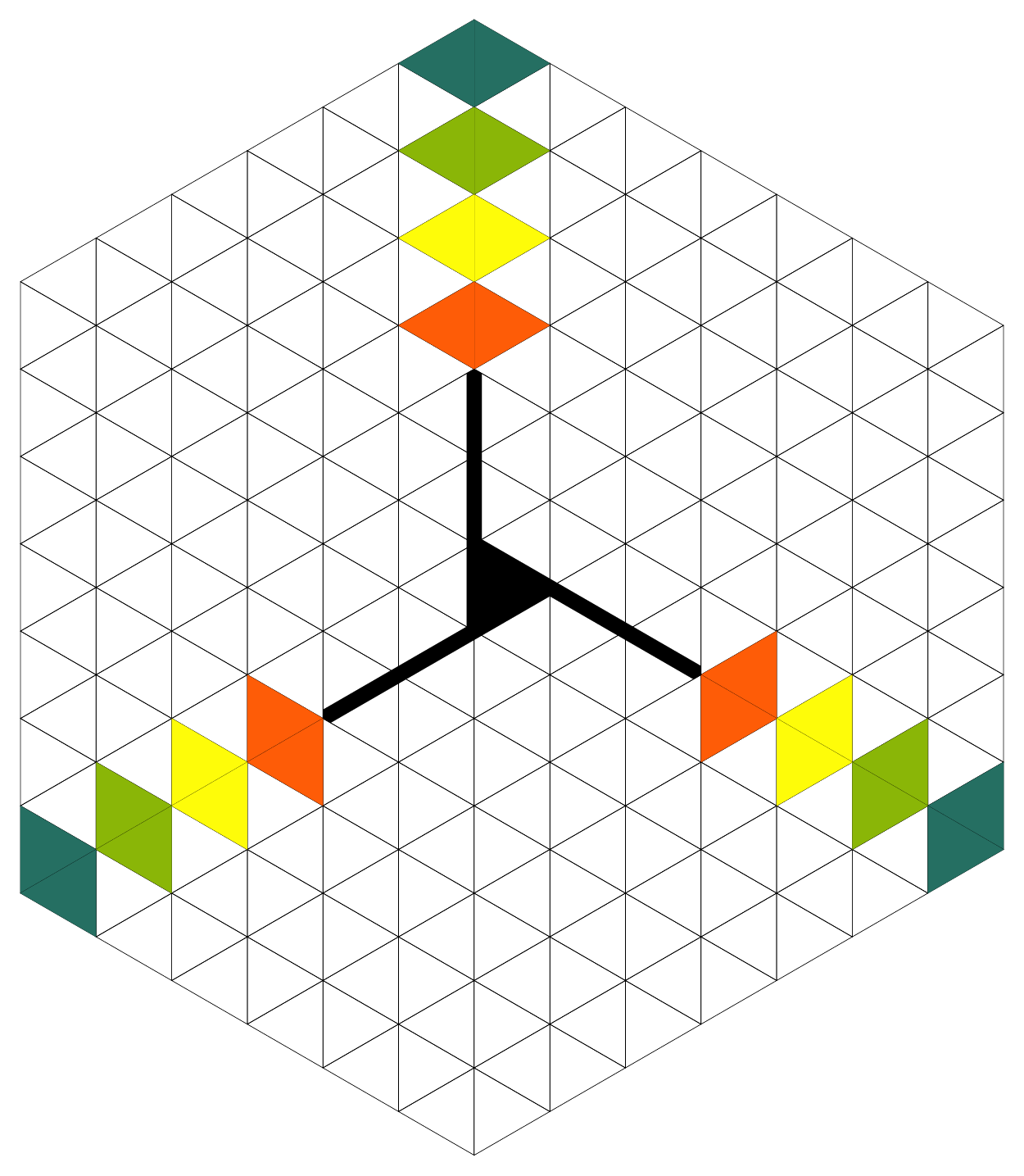}\\
$\downarrow$&$\downarrow$&$\downarrow$&$\downarrow$\\
\includegraphics[width=0.23\textwidth]{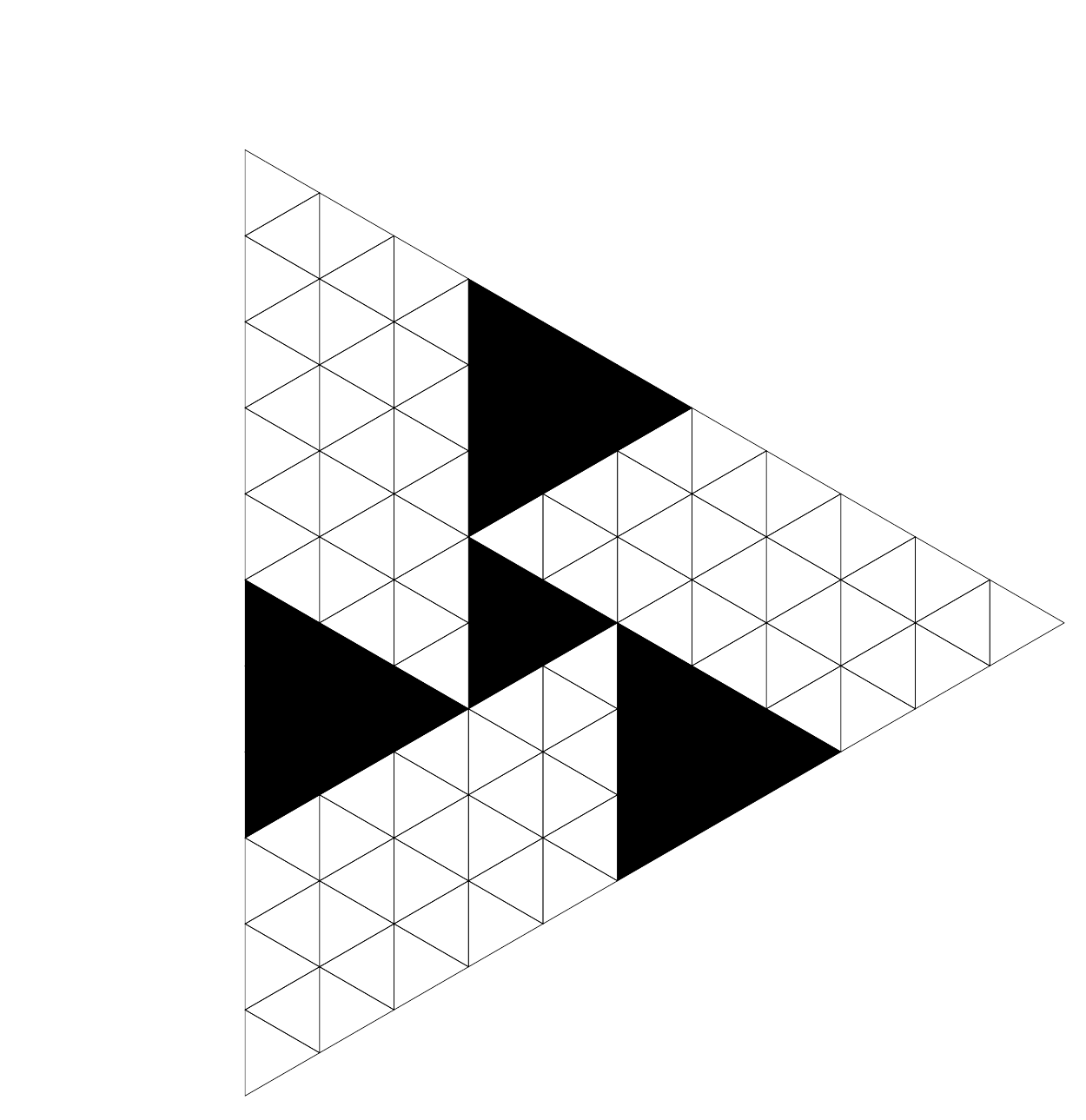}&
\includegraphics[width=0.23\textwidth]{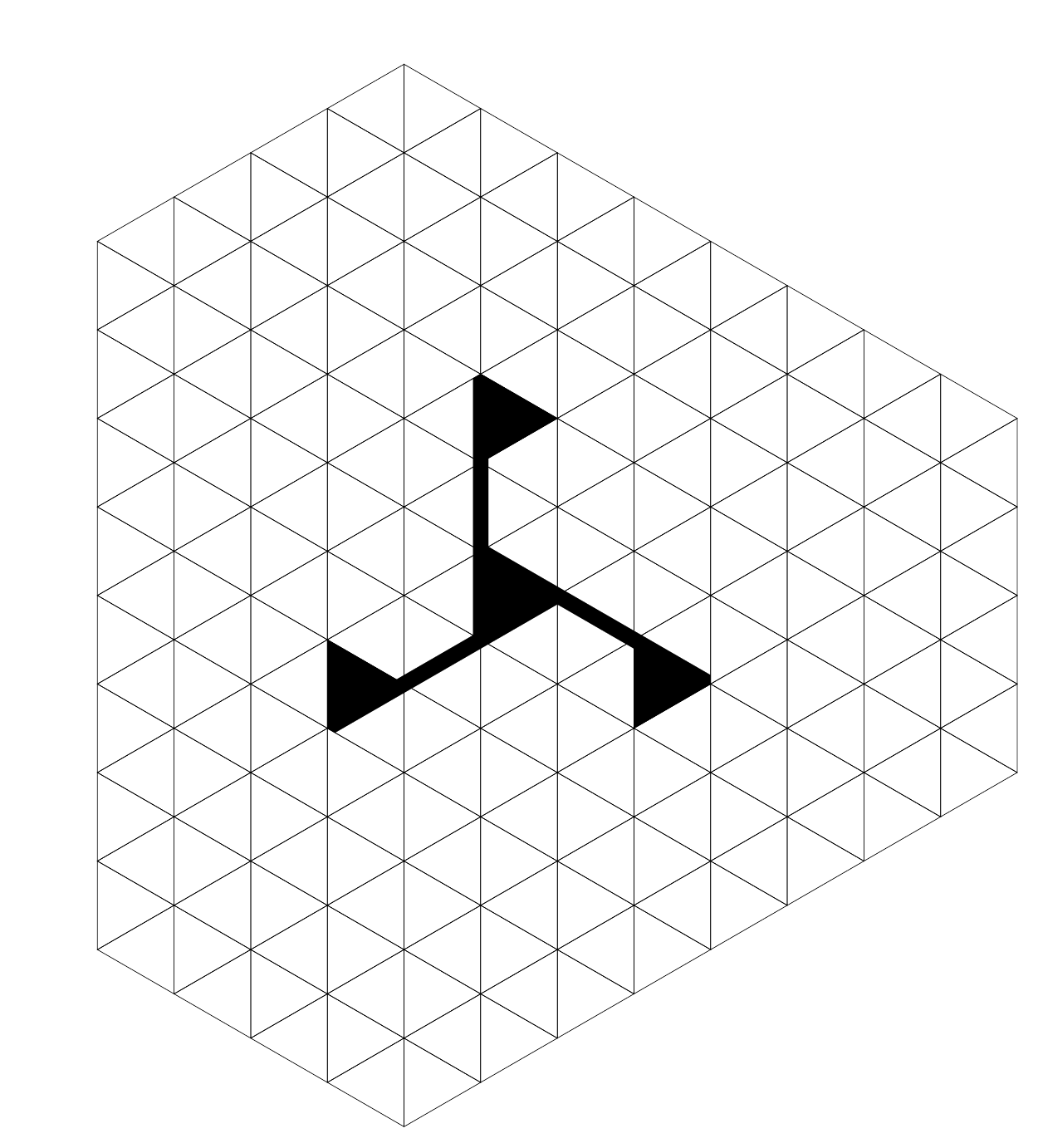}&
\includegraphics[width=0.23\textwidth]{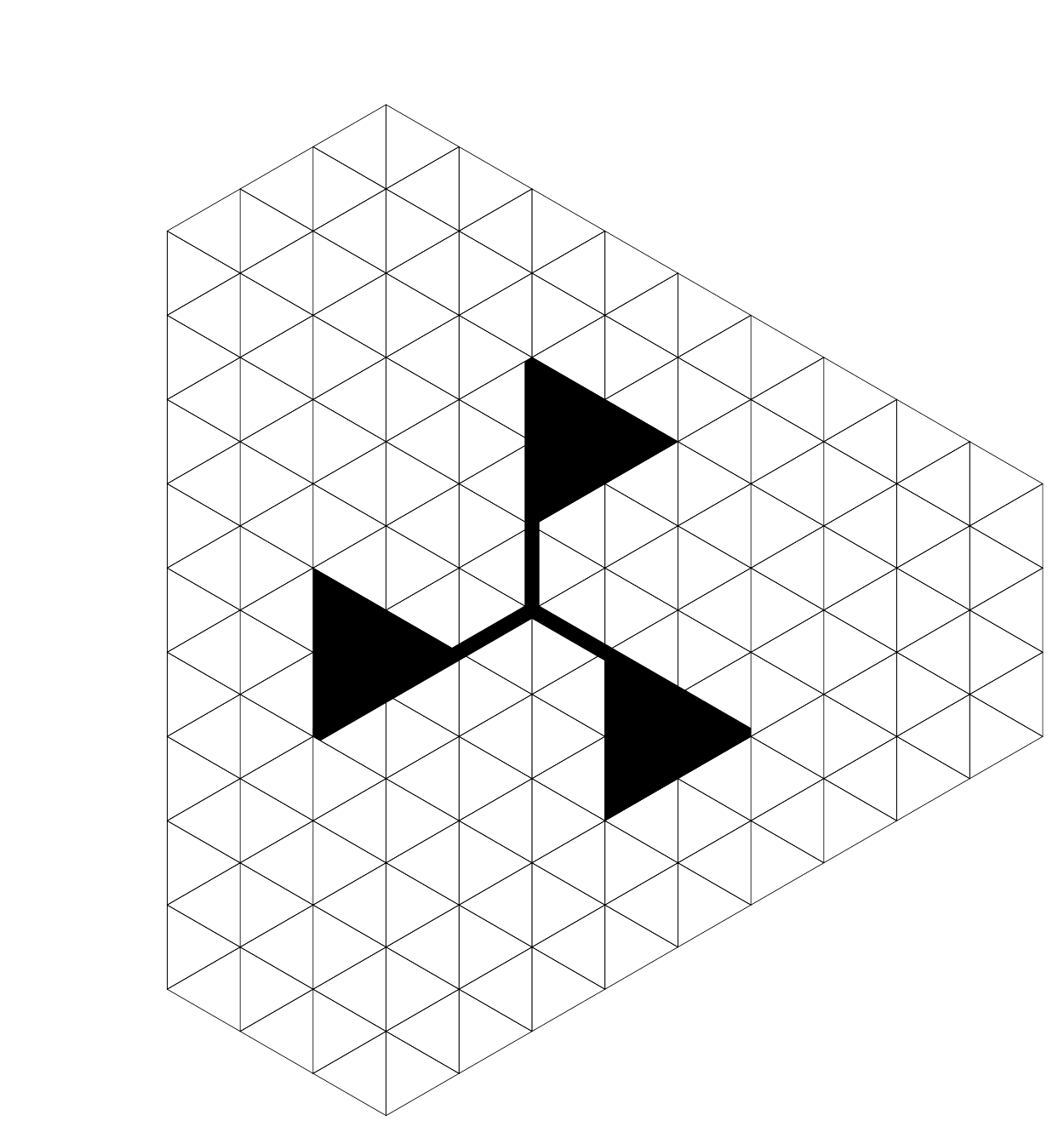}&
\includegraphics[width=0.23\textwidth]{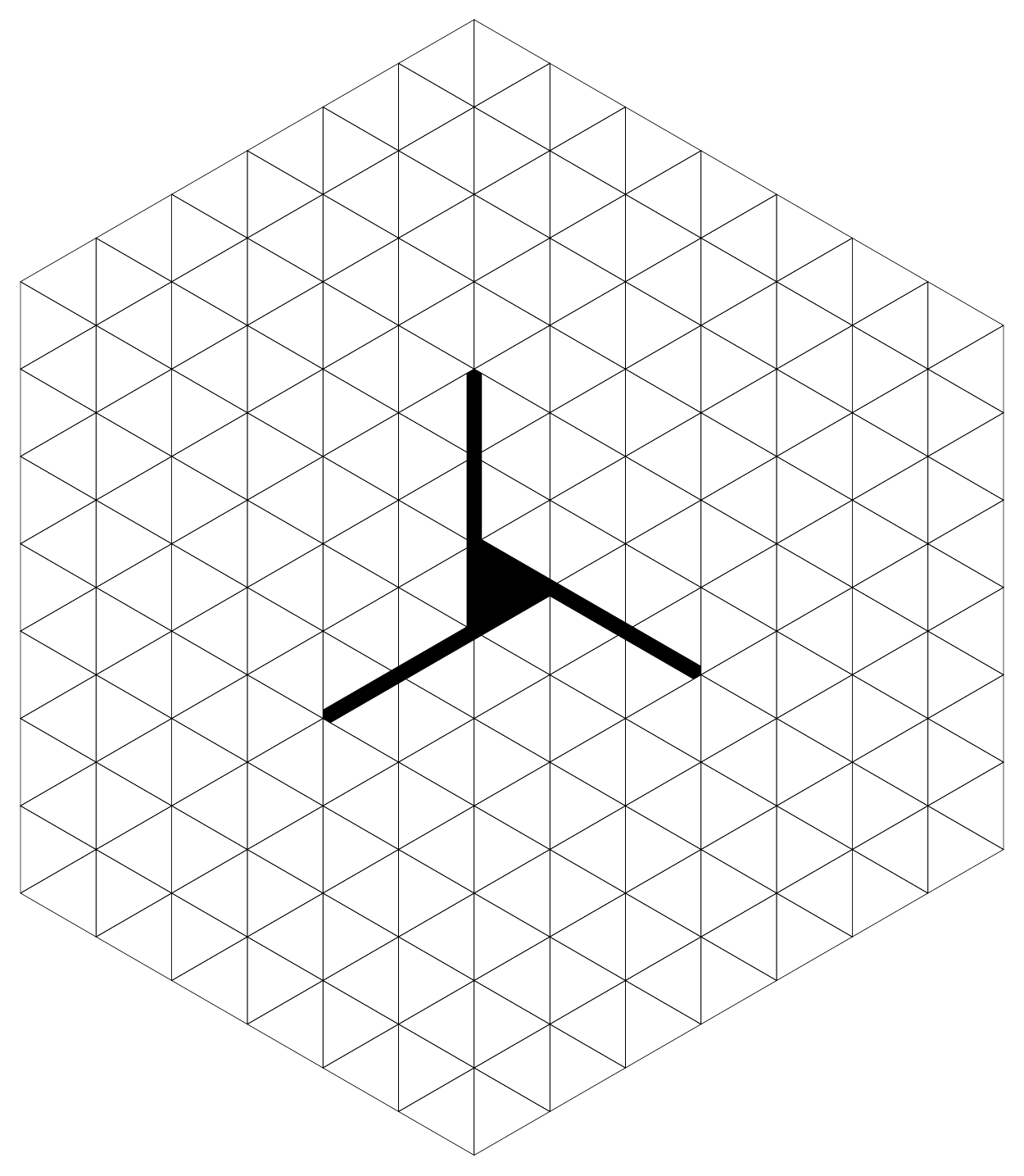}
\end{tabular}
\caption{Selected examples of the region related to the determinant
  $E_{s,t}^{\mu}(n)$. From left to right $(s,t,n,\mu): (3,0,3,4), (2,1,4,3),
  (3,1,4,2), (2,2,4,3).$ All regions are hexagonal after the forced tiling
  removal, except for the left-most one, which becomes triangular. The regions
  will typically have four triangular holes after the forced tiling removal.
  Exceptions are the right-most case, which has no mandatory points, and the
  one to its left, where the central hole has size~$0$.}
\lbl{fig:regionexamples}
\end{figure}

{\parskip=0pt
\begin{itemize}
\item $s=t$ : $\sum_{I\subseteq \lbrace 1,\ldots,n\rbrace}B_I^I$ counts all
  tuples of non-intersecting paths for all subsets of start points (and the
  same subsets of end points), all tilings of a lozenge-shaped region with the appearance of 
  optional points controlled by the set $I$, and equivalently the number of
  cyclically symmetric tilings of the corresponding hexagonal-shaped region
  (which may have a central triangular hole and some border lines as described
  in the above construction).
\item $s>t$ : $\sum_{I\subseteq \lbrace 1,\ldots,n-(s-t)\rbrace}B_{I+s-t}^I$
  counts all tuples of non-intersecting paths for all subsets that must
  contain the last $s-t$ start points and the first $s-t$ end points, all
  tilings of a lozenge-shaped region with the appearance of optional points controlled by the
  set $I$, and equivalently the number of cyclically symmetric tilings of the
  corresponding hexagonal or triangular shaped region (which may have up to
  four central triangular holes and some border lines as described in the
  above construction).
\item $s<t$ : Analogous to $s>t$.
\end{itemize}

Next, we think about the sign. If $s-t+1$ is even, $E^{\mu}_{s,t}(n)$ counts
exactly the number of cyclically symmetric rhombus tilings of the constructed
region as described above, regardless of the parity of $|I|$. This is because
an even sign implies that all of the possible paths/tilings that should be
counted are included in the summation.
}

In the case where $s-t+1$ is odd, we may want to consider which terms are
being cancelled in the sum. The sign $(-1)^{(s-t+1)|I|}$ indicates that we
should think more carefully about the set $I$. Recall that this set controls
the number of horizontal rhombi that crosses the border connections of the
lozenges, in other words, they control the number of optional points that are
present or absent. One way that we can take advantage of this fact is to use
the symmetry of the region to be tiled to deduce conditions on $n$ and $s$ for
which we can definitely see a cancellation.

We use the example of $s>t$ where $t=0$, and refer to \fig{canceltilings}
for a visual. We are now in the case where we do not have a border line. This
means that the four central triangular holes force the tiling of three
lozenge-shaped regions (which are indicated in gray in \fig{canceltilings}).
Their removal will not affect the final tiling count. Thus, we
can view our tiling region to be a hexagon with only one large central
triangular hole! We observe that there are three lines of symmetry of this new
triangular hole (in fact, they are the lines of symmetry for the whole region,
one of which is shown on the left in \fig{canceltilings}).

The collection of starting points that get removed by the set $I$ trigger
$n-s-|I|$ paths that start at the left edge of the green lozenge and exit it
at its bottom. They then enter the yellow lozenge and continue to meander
around the central hole, until they complete their cycle. The third (purple)
lozenge contains the red triangle. Note that the $n-s-|I|$ tiles crossing the
vertical side of the red triangle imply that exactly $|I|$ tiles will cross its
lower-left boundary (an example of such behavior is shown on
the right in \fig{canceltilings}). The symmetry of the figure (along the
dashed line) now implies that there is a one-to-one correspondence between
tilings with $n-s-|I|$ crossings and tilings with $|I|$ crossings. And so, if
$n-s$ is odd, we will get our cancellations in \eqref{eq:sumofminors}. This
allows us to deduce the identity
\[
E_{0,0}^{\mu}(2m-1)=0
\]
combinatorially (using $n=2m-1$ and $s=0$). Since the exact same reasoning \cite{KoutschanThanatipanonda19} was used to deduce that $D_{1,0}^\mu(2m)=0$, we can conjecture a more general identity to relate the $D$ and $E$ determinants and then use a combinatorial argument to resolve it.

\begin{figure}[ht]
  \centering
  \includegraphics[height=200pt]{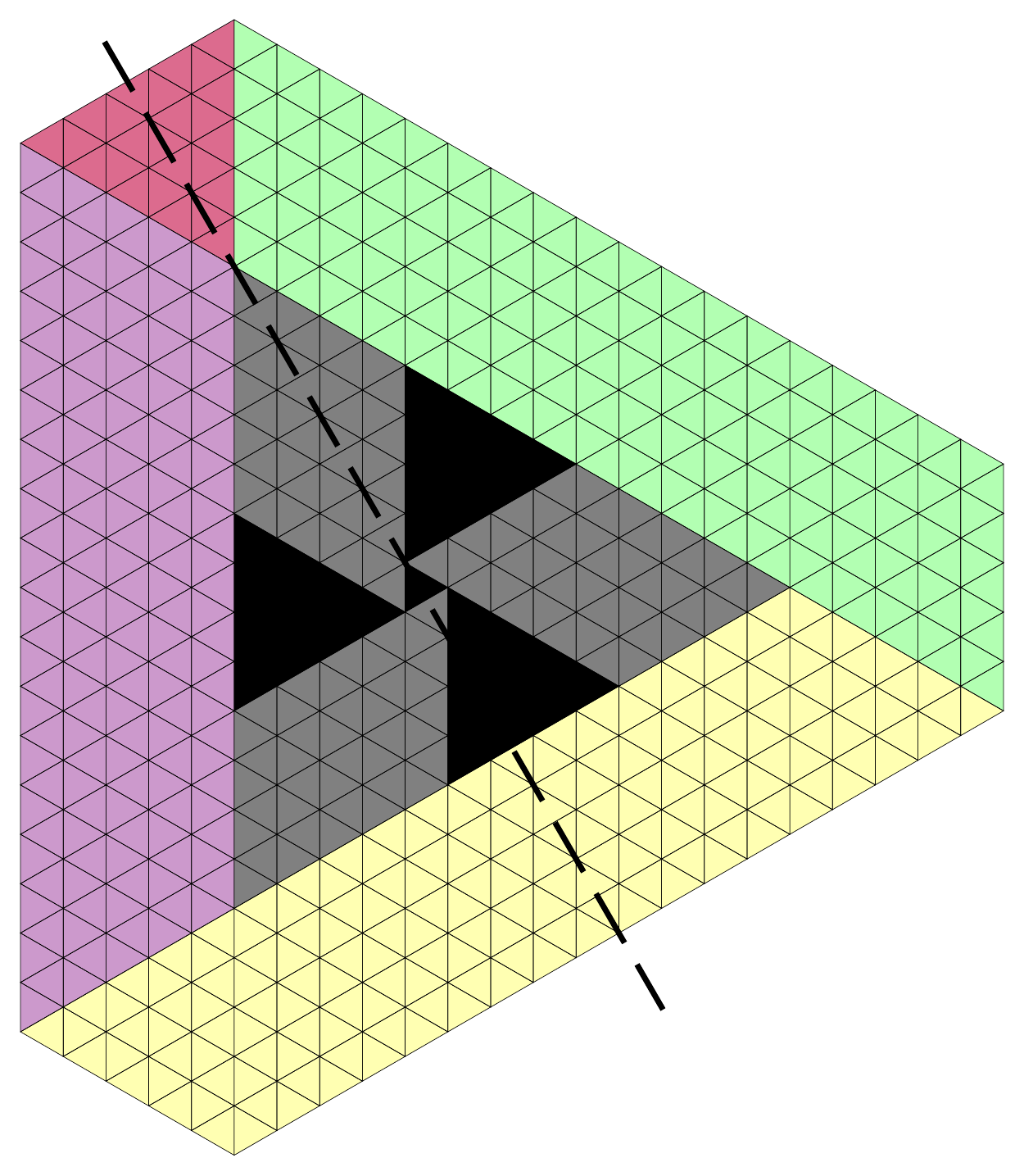}
  \hspace{60pt}
  \raisebox{45pt}{\includegraphics[height=120pt]{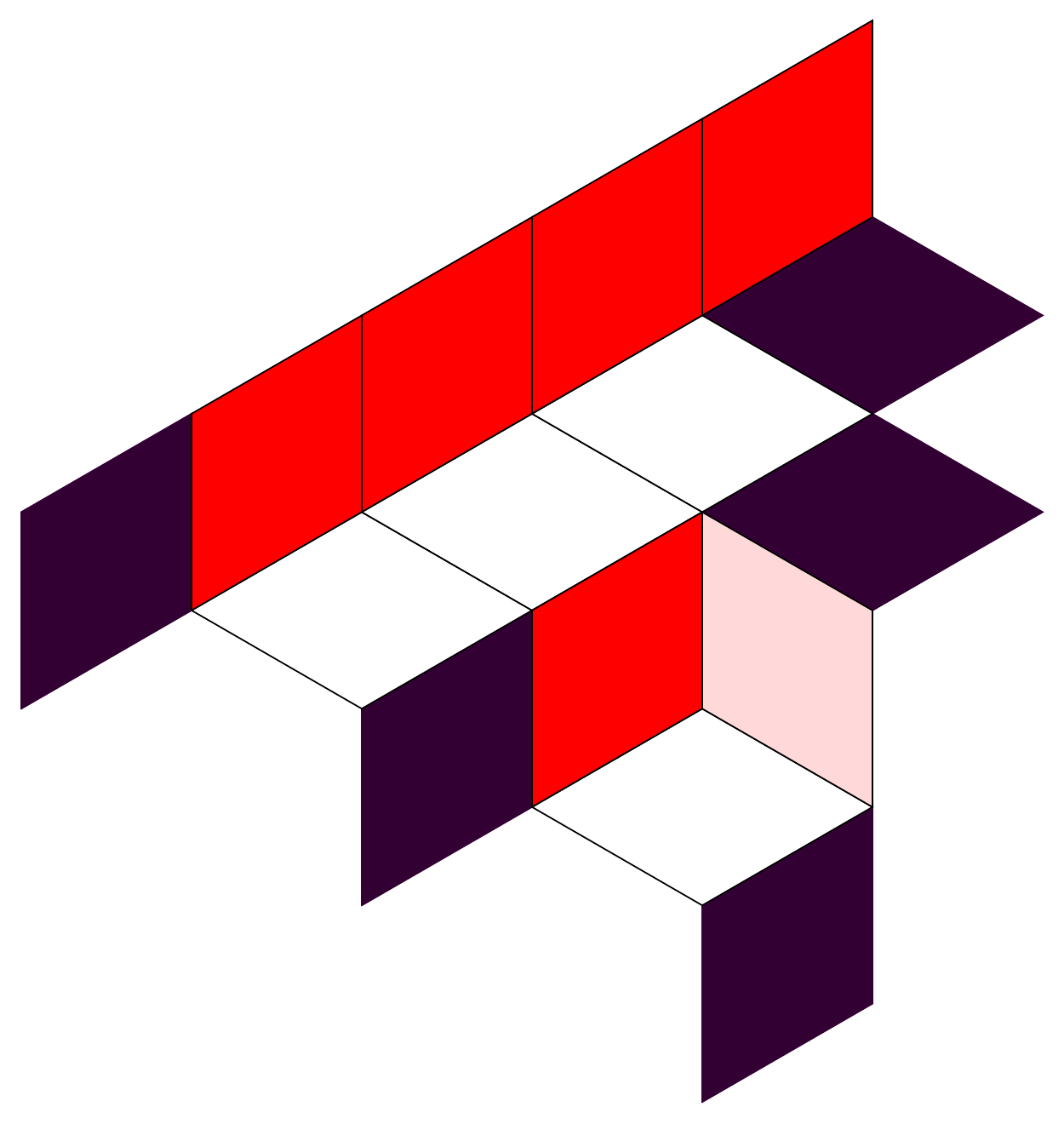}}
\caption{On the left, we present the hexagonal region associated to
  $(s,t,n,\mu)=(4,0,9,3)$. Observe that $s-t+1=3$ and $n-s=5$ are both
  odd. Forced tilings are shaded in gray, and the regions to be tiled are in
  different (lighter) colors. On the right is a zoomed in version of the red triangle,
  with an example tiling triggered by the set $I$ where $|I|=3$ tiles are
  removed. This gives us paths that are indicated by the $n-s-|I|=2$ dark
  tiles along its vertical edge that must end on the other side. Note that
  this forces $|I|=3$ dark tiles to appear on the other side.}
\lbl{fig:canceltilings}
\end{figure}

\begin{lemma}\lbl{lem:famA}
For an indeterminate $\mu$ and $n,s\in\Z$ such that $n\geq s \geq 1$ and $n >1$,
\begin{align*}
E_{s,0}^{\mu}(n)&=D_{s-1,0}^{\mu+3}(n-1),\\
D_{s,0}^{\mu}(n)&=E_{s-1,0}^{\mu+3}(n-1).
\end{align*}
\end{lemma}
\begin{proof}
If $n$ and $s$ are arbitrarily fixed integers such that $n\geq s \geq 1$ and
$n >1$, then all determinants in the statement of the lemma are polynomials in
$\mu$ (because the matrix entries themselves are polynomials in $\mu$). For
integral $\mu$ satisfying the condition $\mu+s\geq 2$, we can invoke the
combinatorial interpretation described above. For each identity, we will show
that the determinants (i.e., the polynomials) agree for infinitely many $\mu$,
and this will allow us to make our conclusion.

The first identity (the second identity is deduced analogously)
can be seen using the sum of minors formula, where we
simply need to observe that the number of Kronecker deltas and the sign
remain the same if $n$ and $s$ are both shifted by 1. Thus, the sign patterns
associated to the determinants in the sum are the same (which means we do not
need to separate odd/even cases). So it is enough to show that the smaller
determinants in the summands all count the same objects for all integral
$\mu\geq 2-s$. The construction as described above is applicable for this
purpose since the construction itself does not take into account these
weights, so we can use it for both $\Mat D{s-1}0{\mu+3}{n-1}$ and $\Mat
Es0{\mu}n$ and show that the resulting regions that need to be tiled
are the same.

The lozenge associated to $\Mat D{s-1}0{\mu+3}{n-1}$ actually has a longer
bottom edge (+1) and shorter left edge ($-1$). So as a single lozenge, it is
difficult to argue that we can get the same tilings. But viewed as a holey
hexagon, it is much easier! First, the removal of the triangular region (of
side length~$s$) associated to the mandatory points on the longer edge forces
the hexagon to be the same size as the one associated to $\Mat Es0{\mu}n$. So
now, we just need to argue about the four triangular holes in the middle,
which have different sizes! But we are in the case $t=0$, and the magic is
that there is only one way to tile the three lozenge-shaped regions that are forced by the four triangular holes. The removal of the additional forced tiling creates the larger triangular hole. For $\Mat
Es0{\mu}n$, this larger triangle has length $\mu-2+3s$ and for $\Mat
D{s-1}0{\mu+3}{n-1}$, this larger triangle has the same length: $(\mu+3)-2+3(s-1)=\mu-2+3s$
(see \fig{lem1pics}). If $n=s$, then there
is only one big triangle in both cases, so this is trivially equal.
We also remark that an algebraic proof of this lemma can be realized by applying certain row and column operations to the matrices. For this, we refer the reader to the discussion in the proof of \lem{biglemma1}.
\begin{figure}[ht]
\centering
\setlength{\tabcolsep}{20pt}
\begin{tabular}{c|c}
\includegraphics[scale=0.34]{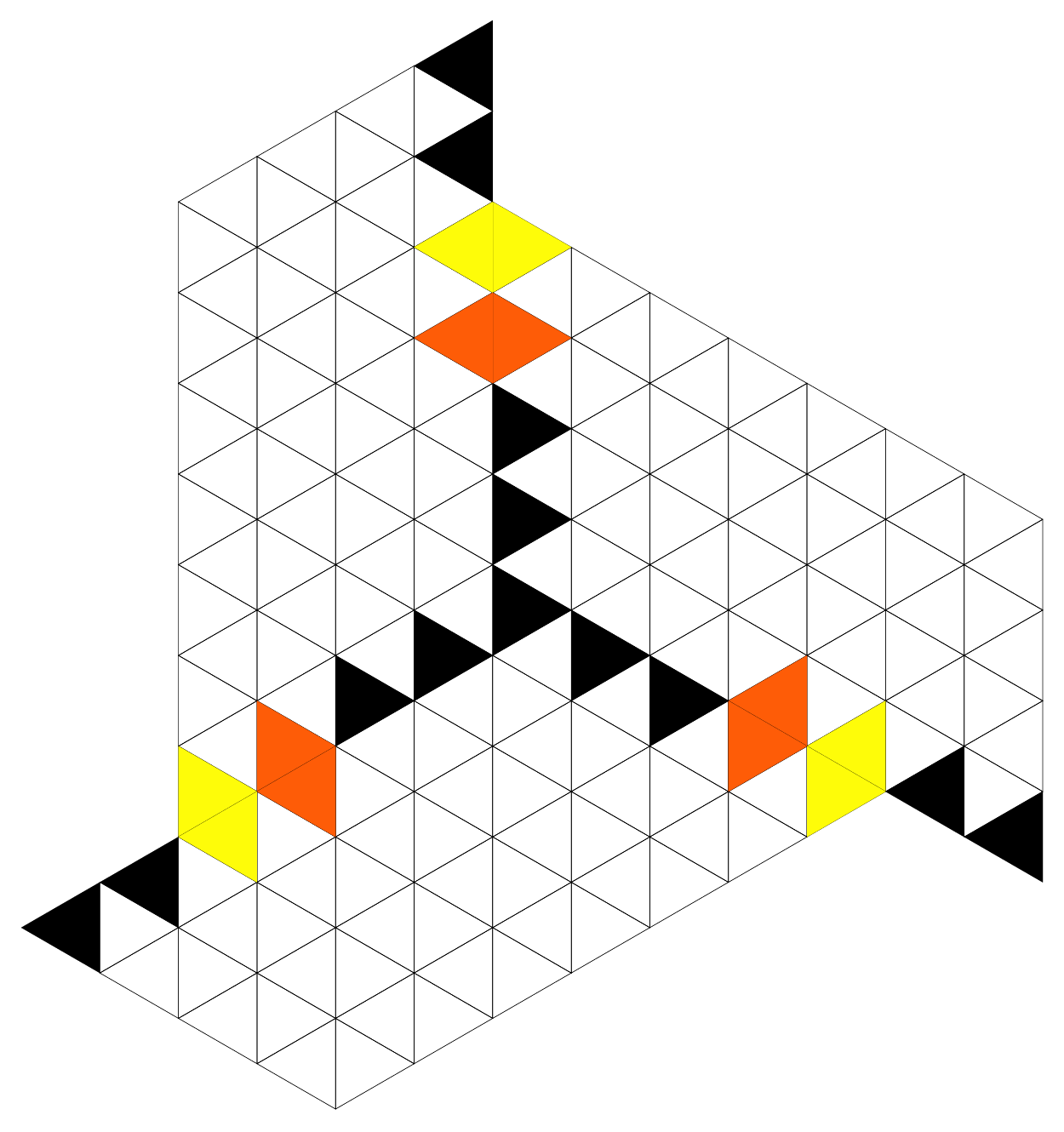}&
\includegraphics[scale=0.34]{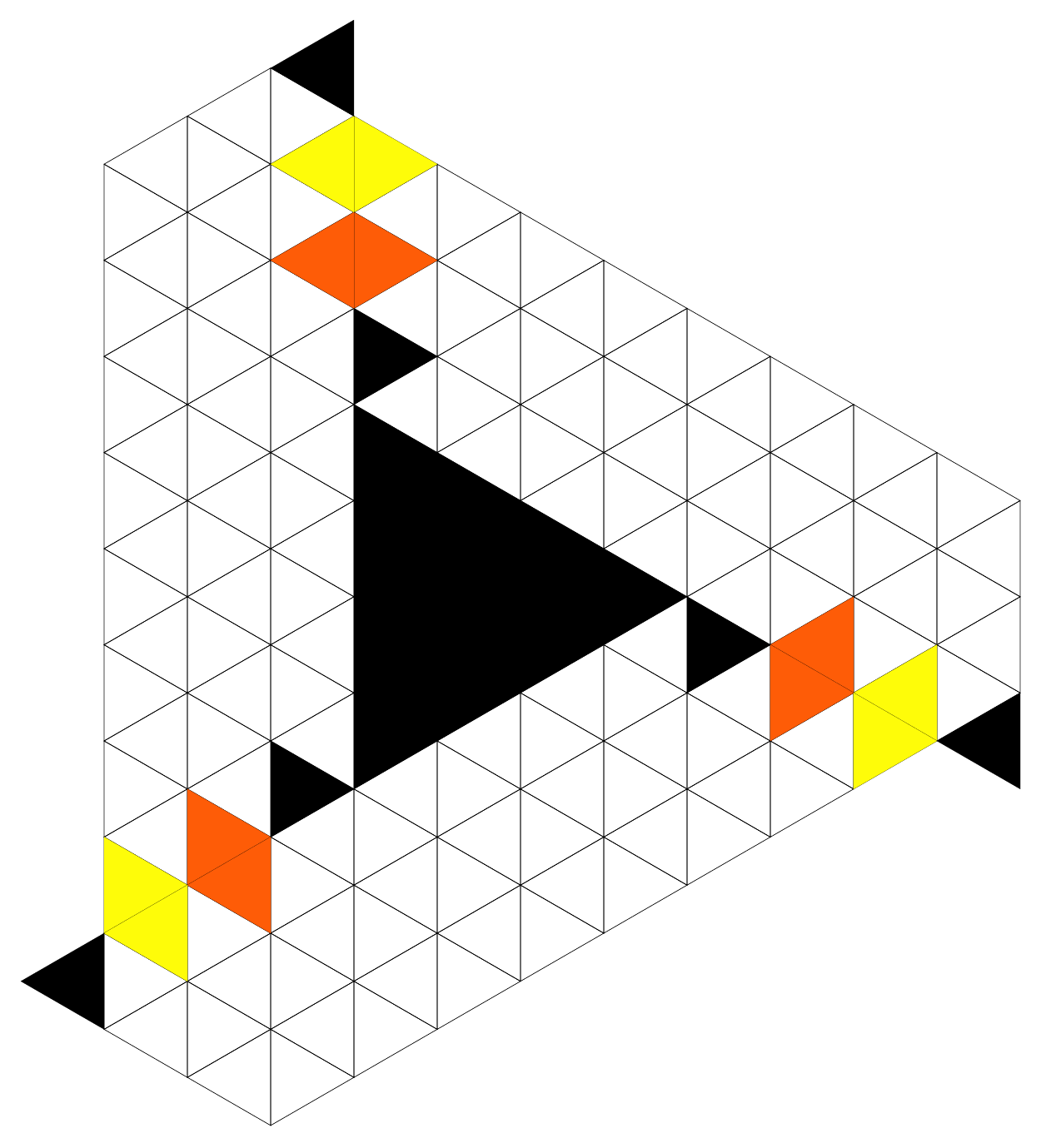}\\
$\downarrow$&$\downarrow$\\
\includegraphics[scale=0.3]{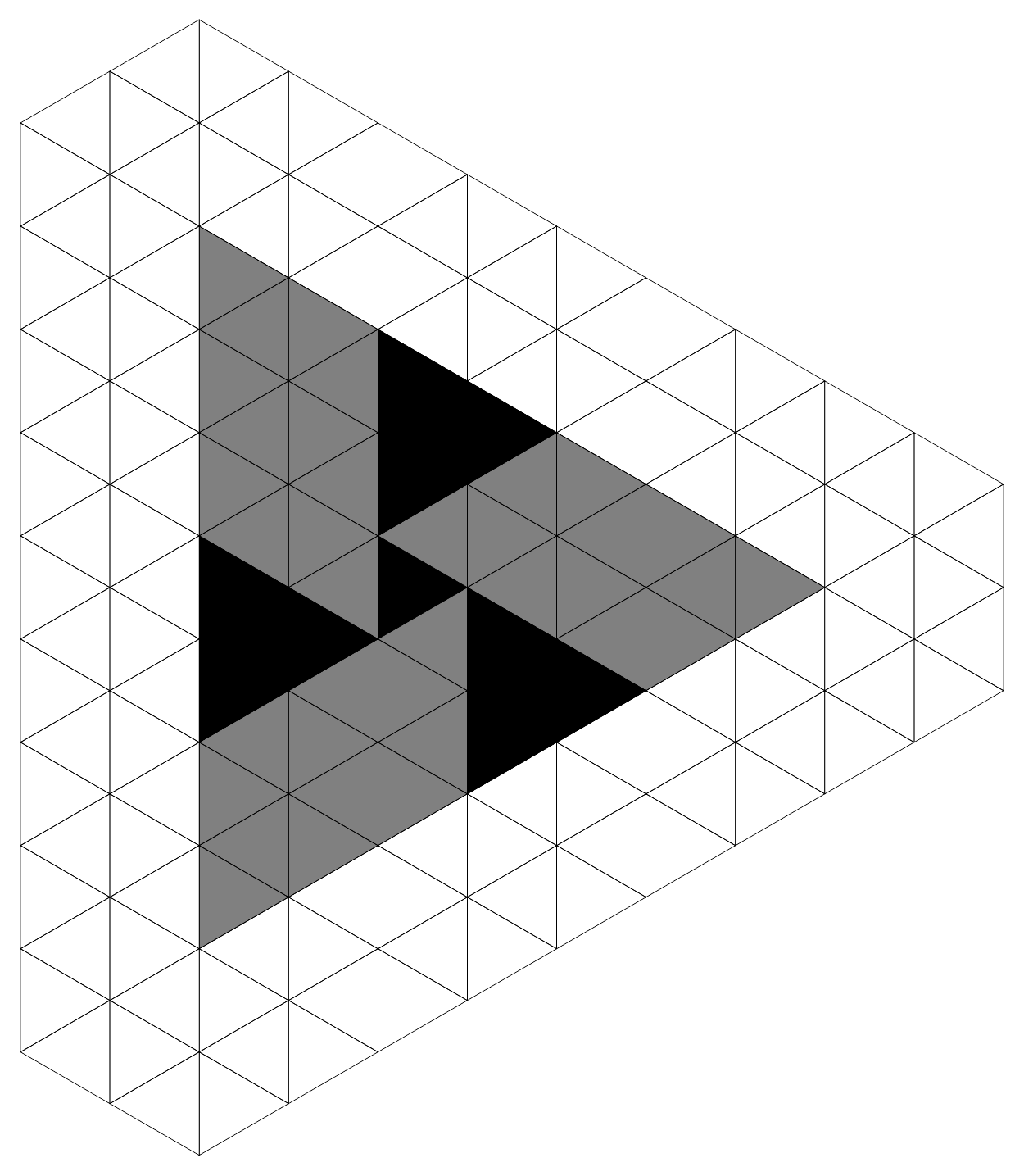}&
\includegraphics[scale=0.3]{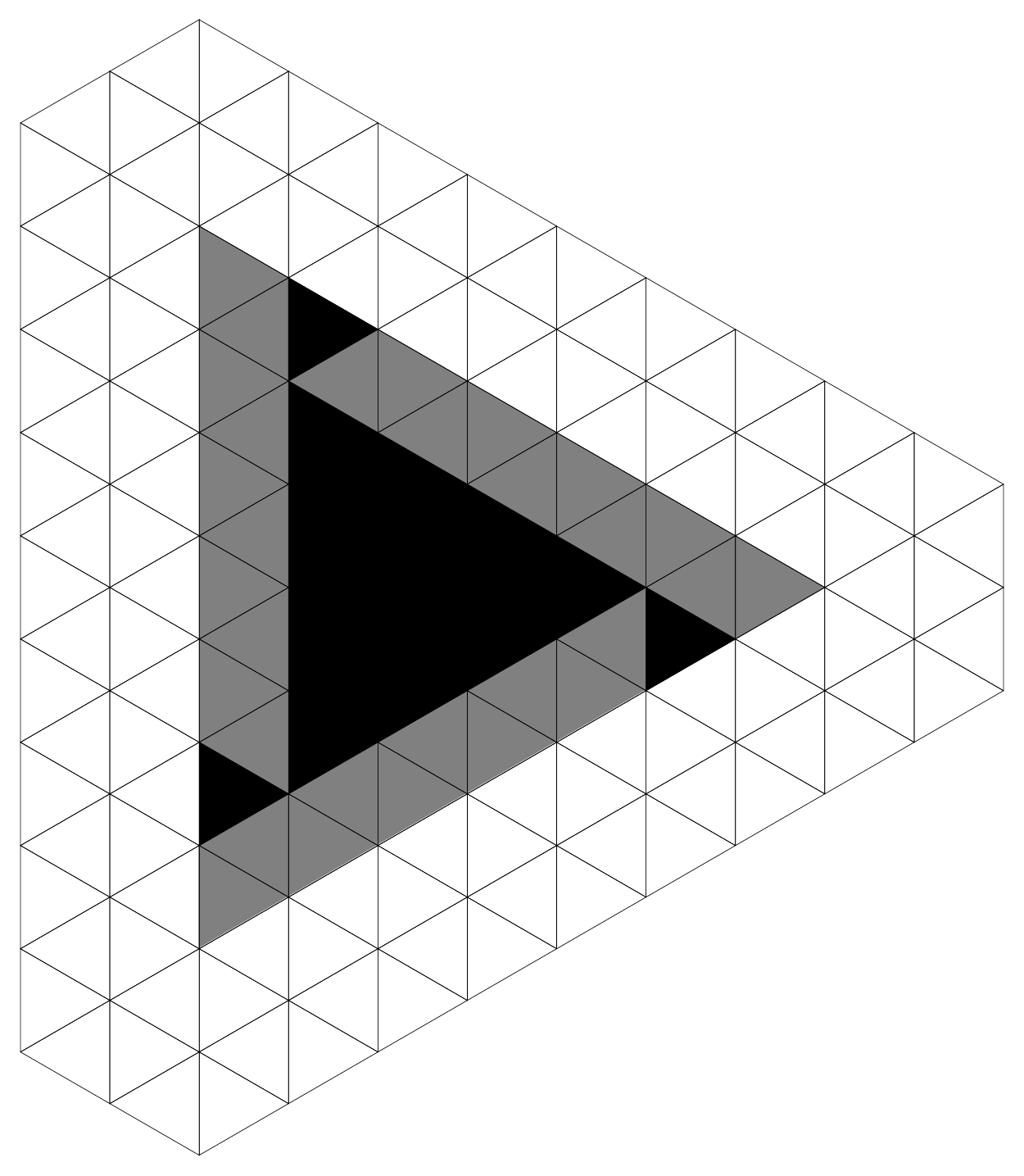}
\end{tabular}
\caption{Pinwheel and hexagon figures for $(s,t,n,\mu)=(2,0,4,3)$ in \lem{famA}.
  The left column is associated to $\Mat Es0{\mu}n$ and the right column is
  associated to $\Mat D{s-1}0{\mu+3}{n-1}$. Note that the region to be tiled
  is the same in both cases (after removing all regions of forced tilings).}
\lbl{fig:lem1pics}
\end{figure}
\end{proof}

\begin{corollary}\lbl{cor:Es0CF}
  For an indeterminate~$\mu$ and $m,r\in\Z$ such that $m>r\geq0$, denote
  \[
    P_{m,r} := \prod_{i=1}^{m-r-1} \frac{
      \bigl(\mu+2i+6r\bigr)_i^2 \, \bigl(\frac{\mu}{2}+2i+3r+1\bigr)_i^2}{
      \bigl(i+1\bigr)_i^2 \, \bigl(\frac{\mu}{2}+i+3r\bigr)_i^2}.
  \]
  Then $\E{2r}0\mu{2m-1}=0$ and
  \begin{align*}
    \E{2r}0\mu{2m} &= \frac{(-1)^{m-r} \, \bigl(\frac{\mu}{2}+3r-\frac12\bigr)_{m-r}}{
      \bigl(\frac12\bigr)_{m-r}} \cdot P_{m,r}, \\
    \E{2r+1}0\mu{2m-1} &= \frac{\bigl(m-r\bigr)_{m-r-1}}{
      \bigl(\frac{\mu}{2}+2m+r-1\bigr)_{m-r-1}} \cdot P_{m,r}, \\
    \E{2r+1}0\mu{2m} &= \frac{2 \, \bigl(\mu+2m+4r+1\bigr)_{m-r-1}}{
      \bigl(\frac{\mu}{2}+m+2r+1\bigr)_{m-r-1}} \cdot P_{m,r}.
  \end{align*}
\end{corollary}
\begin{proof}
  These formulas follow directly from \cite[Theorems~18 and~19]{KoutschanThanatipanonda19}
  by using \lem{famA}.
\end{proof}


\section{Closed Forms for \texorpdfstring{$\bm{\E{2r-1}1{\mu}{2m-1}}$}{E(2r-1,1)(2m-1)}
  and \texorpdfstring{$\bm{\D{2r}1{\mu}{2m}}$}{D(2r,1)(2m)}}
\lbl{sec:krat37ktconj20}

The main goal of this section is to derive closed forms for the
determinants $E_{2r-1,1}^{\mu}(2m-1)$ and $D_{2r,1}^{\mu}(2m)$. This allows us
to resolve two conjectures \cite[Conjecture~37]{Krattenthaler05} and
\cite[Conjecture~20]{KoutschanThanatipanonda19}. We note that this is the
first time that we are able to prove non-trivial results for whole families of
determinants (with $s$ or $t$ containing a parameter).  The roadmap for how we
do this can be seen in \fig{3Dfamilies} in the color \ckrat\ and summarized as follows:
{\parskip=0pt
\begin{itemize}
\item The key result is \lem{biglemma1}, where we establish the ratios between
  families $E_{2r-1,1}^{\mu}(2m-1)$ and $D_{2r,1}^{\mu}(2m)$.
\item This connection between the determinants along with the base case
  $E_{1,1}^{\mu}(2m-1)$, whose closed form was already derived in
  \cite[Theorem~2]{KoutschanThanatipanonda13} and presented in \prop{E11},
  allows us to realize the first main result, a closed form for
  $E_{2r-1,1}^{\mu}(2m-1)$ in \thm{Krat37nice}.
\item Applying \lem{switch} (``switching'') to this closed form and performing
  some algebraic manipulations, we demonstrate in \thm{Krat37ugly} that our
  result matches the conjectured formula for $\E1{2r-1}\mu{2m-1}$.
\item Finally, we also apply \lem{biglemma1} to $E_{2r-1,1}^{\mu}(2m-1)$ to
  deduce the second main result of this section, a closed form for
  $D_{2r,1}^{\mu}(2m)$ in \thm{KTConj20}.
\end{itemize}

In \lem{biglemma1}, we first process the matrix by multiplying with two
elementary matrices~$\matnot{L}_{n}$ and~$\matnot{R}_{n}$, which we define
in~\eqref{eq:LandR}. Then we apply a variant of the holonomic ansatz as
described in \sect{holonomicansatz}, to set up the problem so that the
computer can be used to prove our result with the machinery described in
\sect{computeralgebra}. The introduction of the new parameter~$r$ causes more
difficulties in the calculation than usual. We discuss these difficulties in
the proof below.
}
 
\begin{lemma}\lbl{lem:biglemma1}
Let $\mu$ be an indeterminate, and $m,r\in\Z$. If $m\geq r \geq 1$, then
\begin{align}
  \lbl{eq:biglemeq1}
  \frac{\D{2r}1{\mu}{2m}}{\E{2r-1}1{\mu+3}{2m-1}} &=
  \frac{(m+r-1)(\mu-1)(\mu+2m+1)(\mu+2r)}{2m(2r-1)(\mu+2)(\mu+2m+2r-1)}, \\
  \lbl{eq:biglemeq2}
  \frac{\E{2r+1}1{\mu}{2m+1}}{\D{2r}1{\mu+3}{2m}} &=
  \frac{(m+r)(\mu-1)(\mu+2m+2)(\mu+2r+1)}{2r(2m+1)(\mu+2)(\mu+2m+2r+1)}.
\end{align}
\end{lemma}
 
\begin{proof}
We first observe that the two identities can be presented in a uniform way:
\[
  \frac{\detnew As1{\mu}n}{\B{s-1}1{\mu+3}{n-1}} =
  \frac{(n+s-2) (\mu-1) (\mu+n+1) (\mu+s)}{2n (s-1) (\mu+2) (\mu+n+s-1)}
  =: R_{s,1}^{\mu}(n),
\]
where $(A,B,s,n)=(D,E,2r,2m)$ or $(A,B,s,n)=(E,D,2r+1,2m+1)$. Since we are
dealing with ratios of determinants, we first make sure that a division by
zero will not occur. To do this, we employ an inductive argument with respect
to~$n$, in order to show that all determinants that will be used in the proof
are nonzero: for the induction base, we note that $\E11{\mu}n\neq0$ by
\cite[Theorem~2]{KoutschanThanatipanonda13} (see also \prop{E11}), the induction
hypothesis is $\B{s-1}1{\mu+3}{n-1}\neq0$, and the induction step is
completed once the identities~\eqref{eq:biglemeq1} and~\eqref{eq:biglemeq2} are
established (note that both ratios on the right-hand sides are never identically
zero under the stated assumption $m\geq r\geq1$).  In each step of the
induction, the roles of $D$ and $E$ are interchanged, which is reflected by
the zigzag arrangement of the \ckrat\ connections in \fig{3Dfamilies}.

Next, we manipulate the matrix $\Mat As1{\mu}n$ so that its determinantal value
remains unaffected:
\[
  \matnot{L}_{n}\cdot \Mat As1{\mu}n\cdot \matnot{R}_{n}=:\Mattilde As1{\mu}n,
\]
where $\matnot{L}_n,\matnot{R}_n\in\R^{n\times n}$ are such that
\begin{equation}\lbl{eq:LandR}
  \matnot{L}_{n} := \begin{pmatrix}
    1&0&0&0&\cdots \\
    -1&1&0&0&\cdots \\
    0&-1&1&0&\cdots \\
    0&0&-1&1&\cdots \\
    \vdots & \vdots & \vdots & \vdots & \ddots
  \end{pmatrix}
  \qquad \mbox{and} \qquad
  \matnot{R}_{n} := \begin{pmatrix}
    1&1&1&1&\cdots \\
    0&1&1&1&\cdots \\
    0&0&1&1&\cdots \\
    0&0&0&1&\cdots \\
    \vdots & \vdots & \vdots & \vdots & \ddots
  \end{pmatrix}.
\end{equation}
The matrices $\matnot{L}_n$ and $\matnot{R}_n$ perform elementary row
resp.\ column operations that exploit the elementary
property~\eqref{eq:pascal} of the binomial coefficient. We also note that the
determinants of both matrices are~$1$.

Using \lem{pascal}, the resulting matrix is
\[
\Mattilde As1{\mu}n=
\left(
\begin{array}{c:c}
  \binom{\mu+s-1}{1} & \binom{\mu+j+s-1}{j}-1\pm\sum\limits_{k=1}^j\delta_{s,k} \\
  & \scriptstyle(2\;\leq\;j\;\leq\;n) \\
  \hdashline & \\[-12pt]
  1 & \binom{\mu+i+j+s-3}{j-1} \mp \delta_{s,j-i+2} \\
  \scriptstyle(2\;\leq\;i\;\leq\;n) & \scriptstyle(2\;\leq\;i,j\;\leq\;n)
\end{array}
\right),
\]
where $\pm$ is $+$ if $A=D$ and $-$ if $A=E$ (and $\mp$ is $-$ if $A=D$ and
$+$ if $A=E$). We observe that the bottom right $(n-1)\times(n-1)$ submatrix
is $\Mat B{s-1}1{\mu+3}{n-1}$. In other words, the ``other'' family (i.e., a matrix with the Kronecker delta of opposite sign modulo shifts in $\mu$ and $s$) appears. We can now adapt the holonomic
ansatz (see \sect{holonomicansatz}) to our problem. To compute the determinant
of $\Mattilde As1{\mu}n$, we choose to expand about the first row (rather than
the last row) to get
\[
\dettilde As1{\mu}n=\tilde{a}_{1,1}\cdot \Cof11n +\cdots+\tilde{a}_{1,n}\cdot \Cof1nn,
\]
where $\tilde{a}_{i,j}$ is the $(i,j)$-entry of $\Mattilde As1{\mu}n$, and
$\Cof1jn$ is the corresponding cofactor.  Then $\B{s-1}1{\mu+3}{n-1}=
\Cof11n$, which by the induction hypothesis is nonzero. Hence, we can define
\begin{equation}\lbl{eq:cnj}
  c_{n,j}:=\frac{\Cof1jn}{\Cof11n},
\end{equation}
and our formulas \eqref{eq:biglemeq1} and \eqref{eq:biglemeq2} will be
confirmed by showing that for all~$n\geq s$: 
\begin{equation}\lbl{eq:sys1eq3}
\sum\limits_{j=1}^{n}\tilde{a}_{1,j}\cdot c_{n,j} = R_{s,1}^{\mu}(n).
\end{equation}
Unfortunately, we cannot get a closed form for the $c_{n,j}$'s for symbolic
$n$ and~$j$, in order to prove~\eqref{eq:sys1eq3}. Instead, we will construct
an implicit description of the bivariate sequence $c_{n,j}$ in terms of
recurrences, and then employ the holonomic framework to
prove~\eqref{eq:sys1eq3}.

We can compute many values $c_{n,j}$ explicitly for fixed integers $n$ and
$j$, and then proceed to ``guess'' (i.e., interpolate) recurrences which
the $c_{n,j}$'s satisfy, using an appropriate guessing program (for our
purposes, we use \cite{Kauers09}). However, since these recurrences were
obtained from a finite amount of data, we need to substantiate their universal
validity, that is \emph{for all}~$n$ and~$j$. This is done by observing that
the following identities \textit{uniquely} characterize the $c_{n,j}$'s
\begin{equation}\lbl{eq:sys1}
\begin{cases}
  \; c_{n,1}=1, & n\geq 1, \\[1ex]
  \; \displaystyle\sum_{j=1}^n \tilde{a}_{i,j} \cdot c_{n,j}=0, & 2\leq i \leq n,
\end{cases}
\end{equation}
because by the induction hypothesis, the matrix
$\Mat B{s-1}1{\mu+3}{n-1}=(\tilde{a}_{i,j})_{2\leq i,j\leq n}$ has full
rank. Hence, if we confirm that a suitable solution of the guessed
recurrences satisfies~\eqref{eq:sys1}, then we can conclude that it completely
agrees with $c_{n,j}$. Of course, we will employ the holonomic framework
for this task. Less importantly, we remark that the data generation for the
guessing is achieved more efficiently using the system~\eqref{eq:sys1} rather
than using the definition~\eqref{eq:cnj} in terms of minors.

To prove \eqref{eq:biglemeq1}, we use $(A,B,s,n)=(D,E,2r,2m)$ and show that
the $c_{2m,j}$ satisfy the identities corresponding to \eqref{eq:sys1} and
\eqref{eq:sys1eq3} for all $m\geq r$:
\begin{align*}
  c_{2m,1} &= 1, \\
  \sum_{j=1}^{2m}\binom{\mu+i+j+2r-3}{j-1}\cdot c_{2m,j} - c_{2m,i+2r-2} &= 0,
  \qquad (2\leq i \leq 2m), \\
  \sum_{j=1}^{2m}\binom{\mu+j+2r-1}{j}\cdot c_{2m,j} - \sum_{j=1}^{2r-1}c_{2m,j}
  &= R_{2r,1}^{\mu}(2m).
\end{align*}
To prove \eqref{eq:biglemeq2}, we use $(A,B,s,n)=(E,D,2r+1,2m+1)$ and show
that the $c_{2m+1,j}$ satisfy the identities corresponding to \eqref{eq:sys1}
and \eqref{eq:sys1eq3} for all $m\geq r$:
\begin{align*}
  c_{2m+1,1} &= 1, \\
  \sum_{j=1}^{2m+1}\binom{\mu+i+j+2r-2}{j-1}\cdot c_{2m+1,j} + c_{2m+1,i+2r-1} &= 0,
  \qquad(2\leq i \leq 2m+1), \\
  \sum_{j=1}^{2m+1}\binom{\mu+j+2r}{j}\cdot c_{2m+1,j}
  -\sum_{j=1}^{2r} c_{2m+1,j} - \sum_{j=2r+1}^{2m+1}2\cdot c_{2m+1,j}
  &= R_{2r+1,1}^{\mu}(2m+1).
\end{align*}
At this point, the computer steps in to do some of the legwork for us, and we
briefly talk here about the computation part of the proof (see
\sect{computeralgebra} and \cite{EM2}). From the guessing step, we already
have the generators for a left annihilating
ideal of the $c$'s and we can see that all of the other constituents in these
identities are binomial coefficients or rational functions in the parameters,
which have the nice property of being holonomic.  We also note that the
summations in the identities have ``natural boundaries" in that the summands
evaluate to zero beyond the summation bounds. This means that when we apply
creative telescoping and closure properties for holonomic functions to our
objects (see \cite{Zeilberger90,Koutschan09}), we expect to be able to deduce
an annihilating ideal for the left-hand sides without further adjustments. We
can simplify things by moving terms that are not a summation to the right-hand
side and computing an annihilating ideal for them separately (this is to avoid
a possible slowdown from the need to apply additional closure properties). The
last step is to confirm that either the annihilating ideals on both sides are
equal, or one is a subideal of the other, along with comparing a
sufficient number of initial values.

In theory, the procedure described above is expected to be relatively
uncomplicated. Unfortunately, in practice it turned out to be a bit painful
and we take a couple of paragraphs to highlight two difficulties that were
encountered during the computation. All of the details can be found in the
online supplementary material \cite{EM2}.

{\parskip=0pt
\begin{enumerate}
\item Creative telescoping on the summation in the second identity
  of~\eqref{eq:sys1} did not finish (we left it running to see if it would, but
  in the third month a water leak in the building destroyed the node the
  computations were on). This meant that we needed a better way to speed up
  the process. This was achieved by interpolating/guessing telescoping
  relations for the generators of the annihilating ideals for the sum. We
  confirmed that our guesses are correct by showing that they lie in the
  annihilating ideal of the summands. We then extracted annihilators for the
  sum from these relations (i.e., the telescopers). A second trick to speed up
  this computation was not to construct the full Gr\"{o}bner basis this way,
  but only a few generators (concretely: two out of three), and then run
  Buchberger's algorithm to obtain the remaining ones.

  The timing to confirm the second identity of~\eqref{eq:sys1} was roughly
  8~hours in each of the two cases, and most of the time was taken to generate
  the data for interpolating the telescoping relations.

\item Applying creative telescoping on the summations in the third identity
  corresponding to~\eqref{eq:sys1eq3} resulted in the appearance of singularities in the
  certificates within the summation range. This meant that we were unable to
  certify that our telescopers were the correct annihilators of the
  sums. There is a way to fix this by hand (see \cite{KoutschanWong21} for
  examples and an easy-to-digest description) which involves removing the
  places where the singularity occurs and collecting inhomogeneous parts to
  compensate for the removal. Using this strategy, the final annihilator for
  each sum would consist of a ``left multiplication'' of the annihilator of
  these inhomogeneous parts to the original telescoper. When we applied this
  strategy, we encountered a problem in our computation because the
  annihilator of one of the inhomogeneous parts was unable to finish computing
  and this required another human interaction to complete the process. In
  particular, the difficulty occurred in a substitution step. So instead
  of applying the substitution command directly (which is an implementation of
  the corresponding closure property), we performed the substitution by hand
  on the coefficients of the computed annihilator and then searched for the
  final annihilator that had the support that we expected after substituting.

  The timing to achieve a ``grand'' recurrence for the left-hand sides of the
  identities \eqref{eq:biglemeq1} and \eqref{eq:biglemeq2} corresponding to
  \eqref{eq:sys1eq3} was roughly 30 hours each, with most of the time taken to
  deal with the inhomogeneous parts. In both cases, the recurrence is of order~$7$
  in~$m$ with an approximate byte count of 66,000,000. The degrees of the
  polynomial coefficients in the parameters $m,r,\mu$ are $47$, $37$, and~$38$,
  respectively. \qedhere
\end{enumerate}
}
\end{proof}

We now introduce a technical lemma that will enable us to convert the formula
for $\E11{\mu}{2m-1}$ given in \cite[Theorem~2]{KoutschanThanatipanonda13}
into a nicer form in \prop{E11}.

\begin{lemma}\lbl{lem:CancelPoch}
  Let $\mu$ be an indeterminate and $m\in\Z$ with $m\geq1$. Then
  \[
    2^{(m-1)(m-2)/2}\cdot \prod_{i=1}^{\left\lfloor m/2\right\rfloor}
    \bigl(\tfrac{\mu}{2}+3i-\tfrac{1}{2}\bigr)_{\!m-2i} \,
    \bigl(\tfrac{\mu}{2}+2m-i\bigr)_{\!m-2i+1}
    =
    \prod_{i=1}^{m-1} \frac{\bigl(\mu+2i+1\bigr)_{i-1} \, \bigl(\frac{\mu}{2}+2i+1\bigr)_{i}}
         {\bigl(\frac{\mu}{2}+i+1\bigr)_{i-1}}.
  \]
\end{lemma}
\begin{proof}
  The proof goes by induction with respect to~$m$. Let $L_m$ and $R_m$ denote
  the left-hand (resp. right-hand) side of the statement. For $m=1$, we get
  $L_1=1=R_1$. For all integers $m\geq1$, we have
  the relations
  \begin{align*}
    \frac{R_{m+1}}{R_{m}} &=
    \frac{\bigl(\mu+2m+1\bigr)_{m-1} \, \bigl(\frac{\mu}{2}+2m+1\bigr)_{m}}{
      \bigl(\frac{\mu}{2}+m+1\bigr)_{m-1}},
    \\[2ex]
    \frac{L_{m+1}}{L_{m}} &=
    \begin{cases}
       \displaystyle
       \frac{2^{m-1} \left(\frac{\mu}{2}+m+\frac{1}{2}\right)_{\frac{m}{2}}
         \left(\frac{\mu}{2}+\frac{3m}{2}+1\right)_{\frac{3m}{2}}
         }{
         \left(\frac{\mu}{2}+\frac{3m}{2}\right)_{\frac{m}{2}}
         \left(\frac{\mu}{2}+\frac{3m}{2}+1\right)_{\frac{m}{2}}
         }
       & \mbox{if $m$ is even,} \\[4ex]
       \displaystyle
       \frac{2^{m-1} \left(\frac{\mu}{2}+m+\frac{1}{2}\right)_{\frac{m-1}{2}}
         \left(\frac{\mu}{2}+\frac{3m}{2}+\frac{3}{2}\right)_{\frac{3m-1}{2}}}{
         \left(\frac{\mu}{2}+\frac{3m}{2}+\frac{1}{2}\right)_{\frac{m-1}{2}}
         \left(\frac{\mu}{2}+\frac{3m}{2}+\frac{3}{2}\right)_{\frac{m-1}{2}}
         }
       & \mbox{if $m$ is odd,}
     \end{cases}
  \end{align*}
  where $L_{m+1}/L_m$ comes from rearranging the Pochhammers so that
  \ref{itm:pochprodstep} can be applied.  Then, a strategic application of the
  Pochhammer properties \ref{itm:pochshift}, \ref{itm:pochconnect},
  \ref{itm:pochinterlace} to $L_{m+1}/L_m$ for both cases enables us to
  conclude that $L_{m+1}/L_m=R_{m+1}/R_m$. Thus, $L_m=R_m$ for all $m\geq1$.
\end{proof}

\prop{E11} presents a closed form for $E_{1,1}^{\mu}(2m-1)$, which will be used as a
base case for our main results (\thm{Krat37nice} and \thm{KTConj20}).

\begin{proposition}\lbl{prop:E11}
  Let $\mu$ be an indeterminate and $m\in\Z$ with $m\geq1$. Then
  \[
    \E11{\mu}{2m-1} = \frac{(-1)^{m-1} \, 2^{2m-1} \, \bigl(\frac{\mu-1}{2}\bigr)_m}{(m)_m}
    \cdot \prod_{i=1}^{m-1} \frac{
      \bigl(\mu+2i+1\bigr)_{i-1}^2 \, \bigl(\frac{\mu}{2}+2i+1\bigr)_i^2}{
      \bigl(i\bigr)_i^2 \, \bigl(\frac{\mu}{2}+i+1\bigr)_{i-1}^2}.
  \]
\end{proposition}
\begin{proof}
  We just need to rewrite the closed form for the determinant $E_{1,1}$, given
  in \cite[Theorem~2]{KoutschanThanatipanonda13}, into the above, more
  compact form. The formula in~\cite{KoutschanThanatipanonda13} reads
  (upon substituting $n\to2m-1$):
  \begin{equation}\lbl{eq:E11KT}
    (-1)^{m-1} \, 2^{m(m+1)} \,
    \bigl(\tfrac{\mu-1}{2}\bigr)_{m} \,
    \Biggl(\prod_{i=0}^{m-1} \frac{i!\,(i+1)!}{(2i)!\,(2i+2)!}\Biggr) \cdot
    \prod_{i=1}^{\lfloor \frac{m}{2}\rfloor}
      \biggl(\bigl(\tfrac{\mu}{2}+3i-\tfrac{1}{2}\bigr)_{m-2i}^2 \,
      \bigl(-\tfrac{\mu}{2}-3m+3i\bigr)_{m-2i+1}^2\biggr).
  \end{equation}
  Note that one advantage of the above formula is that there are no more
  cancellations in the second product. Nevertheless, we would like to bring
  this to a form that will be useful for us. We now reshape the first product as
  \begin{equation}\lbl{eq:iprod}
    \prod_{i=0}^{m-1} \frac{i!\,(i+1)!}{(2i)!\,(2i+2)!} =
    \frac12 \prod_{i=1}^{m-1} \frac{1}{(i+1)_i\,(i+2)_{i+1}} =
    \frac12 \prod_{i=1}^{m-1} \frac{1}{8\,(2i+1)\,(i)_i^2} =
    \frac{2^{1-2m}}{(m)_m}\cdot\prod_{i=1}^{m-1} \frac{1}{(i)_i^2},
  \end{equation}
  and after rewriting
  $\bigl(-\tfrac{\mu}{2}-3m+3i\bigr)_{m-2i+1}^2\stackrel{\ref{itm:negpoch}}{=}\bigl(\tfrac{\mu}{2}+2m-i\bigr)_{m-2i+1}^2$,
  we can apply \lem{CancelPoch} to the second product so that
  \eqref{eq:E11KT} turns into the asserted formula.
\end{proof}

And now, on to the main event: based on the results we have achieved so far,
we derive a closed form for $\E{2r-1}1{\mu}{2m-1}$, thereby resolving
Lascoux and Krattenthaler's conjecture \cite[Conjecture~37]{Krattenthaler05}.
Since the formula in that paper is quite different,
we make the effort in \thm{Krat37ugly} to show that the result here is indeed
equivalent to their formula.

\begin{theorem}\lbl{thm:Krat37nice}
Let $\mu$ be an indeterminate and $m,r\in\Z$. If $m\geq r\geq1$, then
\begin{multline*}
  \E{2r-1}1{\mu}{2m-1} = \\
  \frac{(-1)^{m-r} \, (\mu-1) \, \bigl(\mu+2r-1\bigr)_{2m-2}}
       {(2r-2)!\,\bigl(m+r-1\bigr)_{m-r+1}\bigl(\frac{\mu}{2}+r\bigr)_{m-r}} \cdot
  \prod_{i=1}^{m-r} \frac{
    \bigl(\mu+2i+6r-5\bigr)_{i-1}^2 \, \bigl(\frac{\mu}{2}+2i+3r-2\bigr)_{i}^2}{
    \bigl(i\bigr)_{i}^2 \, \bigl(\frac{\mu}{2}+i+3r-2\bigr)_{i-1}^2}.
\end{multline*}
\end{theorem}
\begin{proof}
  We apply \lem{biglemma1} $(2r-2)$ times:
  \begin{align*}
    \E{2r-1}1{\mu}{2m-1}
    &= R_{2r-1,1}^{\mu}(2m-1) \cdot \D{2r-2}1{\mu+3}{2m-2} \\
    &= R_{2r-1,1}^{\mu}(2m-1) \cdot R_{2r-2,1}^{\mu+3}(2m-2) \cdot \E{2r-3}1{\mu+6}{2m-3}
    = \ldots \\
    &= \left(\prod_{i=0}^{2r-3} R_{2r-1-i,1}^{\mu+3i}(2m-1-i)\right) \cdot
    \E11{\mu+6r-6}{2m-2r+1}.
  \end{align*}
  Next, we calculate the product:
  \begin{align}
    &\notag \prod_{i=0}^{2r-3} R_{2r-1-i,1}^{\mu+3i}(2m-1-i) = \\
    &\notag= \prod_{i=0}^{2r-3} \frac{(2m+2r-2i-4) (\mu+3i-1) (\mu+2m+2i) (\mu+2r+2i-1)}
      {2\,(2m-i-1) (2r-i-2) (\mu+3i+2) (\mu+2m+2r+i-3)} \\
    &\notag= \frac{\mu-1}{\mu+6r-7}\cdot \prod_{i=0}^{2r-3}
      \frac{8\,\bigl(m+r-(2r-3-i)-2\bigr)\bigl(\frac{\mu}{2}+m+i\bigr)
        \bigl(\frac{\mu}{2}+r+i-\frac{1}{2}\bigr)}
      {2\,(2m-(2r-3-i)-1)(2r-(2r-3-i)-2)(\mu+2m+2r+i-3)} \\
    &\lbl{eq:Rprod}= \frac{2^{4r-4} \, (\mu-1) \, \bigl(m-r+1\bigr)_{2r-2} \,
        \bigl(\frac{\mu}{2}+m\bigr)_{2r-2} \, \bigl(\frac{\mu}{2}+r-\frac{1}{2}\bigr)_{2r-2}}
      {(\mu+6r-7) \, (2r-2)! \, \bigl(2m-2r+2\bigr)_{2r-2} \, \bigl(\mu+2m+2r-3\bigr)_{2r-2}},
  \end{align}
  where the third line comes from reverting the order of multiplication for
  some factors and the last line comes from applying
  \ref{itm:pochprodstep}. For $\E11{\mu+6r-6}{2m-2r+1}$, \prop{E11} gives us:
  \[
    \frac{(-1)^{m-r} \, 2^{2m-2r+1} \, \bigl(\frac{\mu+6r-7}{2}\bigr)_{m-r+1}}{(m-r+1)_{m-r+1}} \cdot
    \prod_{i=1}^{m-r} \frac{
      \bigl(\mu+2i+6r-5\bigr)_{i-1}^2 \, \bigl(\frac{\mu}{2}+2i+3r-2\bigr)_{i}^2}{
      \bigl(i\bigr)_{i}^2 \, \bigl(\frac{\mu}{2}+i+3r-2\bigr)_{i-1}^2},
  \]
  and we realize that the product is exactly the same as in the statement of the
  theorem. Hence, it remains to simplify the product of the above prefactor
  times expression~\eqref{eq:Rprod}. After applying the following three rules:
  \begin{align*}
      &\frac{(m-r+1)_{2r-2}}{(m-r+1)_{m-r+1} \, (2m-2r+2)_{2r-2}}
    \stackrel{\ref{itm:pochconnect}}{=} \frac{(m-r+1)_{2r-2}}{(m-r+1)_{m+r-1}}
    \stackrel{\ref{itm:pochshift}}{=} \frac{1}{(m+r-1)_{m-r+1}},
    \\[1ex]
  &\bigl(\tfrac{\mu}{2}+r-\tfrac{1}{2}\bigr)_{2r-2}\cdot
    \frac{\bigl(\frac{\mu+6r-7}{2}\bigr)_{m-r+1}}{(\mu+6r-7)}
    \stackrel{\ref{itm:pochconnect}}{=}
    \tfrac12 \, \bigl(\tfrac{\mu}{2}+r-\tfrac{1}{2}\bigr)_{m+r-2}, \quad \text{and}
    \\[1ex]
    &\bigl(\tfrac{\mu}{2}+m\bigr)_{2r-2}
    \stackrel{\ref{itm:pochshift}}{=}
    \frac{\bigl(\frac{\mu}{2}+r\bigr)_{m+r-2}}{\bigl(\frac{\mu}{2}+r\bigr)_{m-r}},
  \end{align*}
  we obtain
  \[  
  \frac{(-1)^{m-r} \, (\mu-1)}{(2r-2)! \, (m+r-1)_{m-r+1} \,
      \bigl(\frac{\mu}{2}+r\bigr)_{m-r}}
      \cdot
    \frac{
      2^{2m+2r-4}\,
      \bigl(\frac{\mu}{2}+r\bigr)_{m+r-2} \, \bigl(\frac{\mu}{2}+r-\frac{1}{2}\bigr)_{m+r-2}
    }{
      \bigl(\mu+2m+2r-3\bigr)_{2r-2}
    }.
  \]
  We now apply \ref{itm:pochinterlace} followed by \ref{itm:pochshift} to the
  right quotient and arrive at the asserted expression in front of the product.
\end{proof}

We are now going to settle Conjecture~37 of~\cite{Krattenthaler05}, the last
open problem from that paper that concerns our families of matrices. Note that
the entries of the matrix were originally given in a slightly different form,
which can be easily adapted to our setting to see that it corresponds to
$E_{1,2r-1}^{\mu}(2m-1)$.  In \thm{Krat37nice} we already found a closed form
for $E_{2r-1,1}^{\mu}(2m-1)$. In \thm{Krat37ugly}, we show that the conjectured
determinant formula (stated below) is indeed equivalent to our formula,
modulo the switching of indices (by \lem{switch}).

\begin{theorem}\lbl{thm:Krat37ugly}
Let $\mu$ be an indeterminate and $m,r\in\Z$. If $m\geq r\geq 1$, then
\begin{equation}\lbl{eq:KratConj37}
  E_{1,2r-1}^{\mu}(2m-1)=2^{4m-3r}\cdot \ell_1 \cdot \ell_2 \cdot \ell_3
  \cdot \prod_{i=0}^{m-1}\frac{i!\,(i+1)!}{(2i)!\,(2i+2)!},
\end{equation}
where
\begin{align*}
  \ell_1 &:= \prod_{i=0}^{2r-3}i!\cdot
  \prod_{i=0}^{r-2}\frac{\bigl((2m-2i-3)!\bigr)^2}{\bigl((m-i-2)!\bigr)^2 \,
    (2m+2i-1)! \, (2m+2i+1)!}, \\
  \ell_2 &:= (\mu-1)\cdot\left(\tfrac{\mu}{2}+r-\tfrac{1}{2}\right)_{m-r}\cdot
  \prod_{i=1}^{2r-2}(\mu+i-1)_{2m+2r-2i-1},\\
  \ell_3 &:= (-1)^{m-r} \, 2^{(m-r)(m-r-1)}
  \prod_{i=0}^{\left\lfloor\frac{m-r-1}{2}\right\rfloor}
  \left(\tfrac{\mu}{2}+3i+3r-\tfrac{1}{2}\right)^2_{m-r-2i-1}
  \left(-\tfrac{\mu}{2}-3m+3i+3\right)^2_{m-r-2i}.
\end{align*}
Remark: This formula was obtained by first applying the transformation
$\mu\rightarrow \mu+r-1$ and then applying $n\rightarrow 2m-1$ and
$r\rightarrow 2r-1$ to~\cite[Conjecture~37]{Krattenthaler05}.
\end{theorem}
 
\begin{proof}
We would like to exploit the closed form found in \thm{Krat37nice}, but in order
to do so, we need to switch the indices by invoking \lem{switch}:
\[
  \E1{2r-1}{\mu}{2m-1} = \E{2r-1}1{\mu}{2m-1}\cdot
  \prod_{i=0}^{2r-3} \frac{(\mu+i)_{2m-1}}{(i+2)_{2m-1}}.
\]
We split the product in~\eqref{eq:KratConj37} at index $m-r$ and rewrite the
first part using the transformation $m\rightarrow m-r+1$ on the derivation
from~\eqref{eq:iprod}:
\begin{equation}\lbl{eq:splitprod}
  \prod_{i=0}^{m-1}\frac{i!\,(i+1)!}{(2i)!\,(2i+2)!} =
  \left(\prod_{i=1}^{m-r} \frac{1}{(i)_i^2}\right) \cdot
  \underbrace{%
  \frac{2^{2r-2m-1}}{(m-r+1)_{m-r+1}}\cdot
  \prod_{i=m-r+1}^{m-1}\frac{i!\,(i+1)!}{(2i)!\,(2i+2)!}
  }_{=:\,\ell_4}
\end{equation}
Instantiating \lem{CancelPoch} with $m\to m-r+1$ and $\mu\to\mu+6r-6$, we see
that $\ell_3$ combined with the parenthesized product in~\eqref{eq:splitprod}
yields exactly the product in the formula of \thm{Krat37nice}.
We would now like to show the equality of the remaining factors, that is,
\[
  2^{4m-3r} \cdot \ell_1 \cdot \ell_2 \cdot \ell_4
  = \underbrace{\frac{(\mu-1)\,\bigl(\mu+2r-1\bigr)_{2m-2}}
       {(2r-2)!\,\bigl(m+r-1\bigr)_{m-r+1}\bigl(\frac{\mu}{2}+r\bigr)_{m-r}}}%
    _{\text{prefactor from \thm{Krat37nice}}} \cdot
  \prod_{i=0}^{2r-3} \frac{(\mu+i)_{2m-1}}{(i+2)_{2m-1}}.
\]
We split this formula by separating the factors that contain $\mu$ and those
that do not (modulo some power of 2). The proof will therefore be complete
once we can prove the following two identities:
\begin{equation}\lbl{eq:Conj37withoutmu}
  2^{2m-r} \cdot \ell_1 \cdot \ell_4
  = \frac{1}{(2r-2)!\,\bigl(m+r-1\bigr)_{m-r+1}} \cdot\prod_{i=0}^{2r-3} \frac{1}{(i+2)_{2m-1}},
\end{equation}
\begin{equation} \lbl{eq:Conj37withmu}
  2^{2m-2r} \cdot \ell_2
  = \frac{(\mu-1)\,\bigl(\mu+2r-1\bigr)_{2m-2}}
       {\bigl(\frac{\mu}{2}+r\bigr)_{m-r}} \cdot\prod_{i=0}^{2r-3}(\mu+i)_{2m-1}.
\end{equation}
For identity~\eqref{eq:Conj37withoutmu}, we find that
\[
  \prod\limits_{i=0}^{2r-3}i! \cdot(2r-2)!\cdot \prod\limits_{i=0}^{2r-3}(i+2)_{2m-1} =
  \prod\limits_{i=1}^{2r-2} (i+2m-1)! =
  \prod\limits_{i=0}^{r-2}(2i+2m)!\cdot \prod\limits_{i=0}^{r-2}(2i+2m+1)!
\]
and
\begin{align*}
  & \left(\,\prod\limits_{i=0}^{r-2} \frac{\bigl((2m-2i-3)!\bigr)^2}{
    \bigl((m-i-2)!\bigr)^2 \, (2m+2i-1)! \, (2m+2i+1)!}\right) \cdot
  \left(\,\prod\limits_{i=m-r+1}^{m-1}\frac{i! \, (i+1)!}{(2i)! \, (2i+2)!}\right)
  \\
  &= \frac{2^{4-4r}}{\left(m-r+\frac{3}{2}\right)_{r-1}} \cdot
  \prod\limits_{i=0}^{r-2} \frac{1}{(2m+2i-1)! \, (2m+2i+1)!},
\end{align*}
with the latter obtained by rewriting the right product so that the limits in the products are the same and then taking out common factors (resulting in some
cancellations). Then using a repeated application of
\ref{itm:pochinterlace} and \ref{itm:pochconnect}, the quotient of both sides
of \eqref{eq:Conj37withoutmu} yields
\begin{align*}
  2^{-2r+2}\cdot\frac{(m)_{r-1}}{\left(m-r+\frac{3}{2}\right)_{r-1}} \cdot
  \frac{(m+r-1)_{m-r+1}}{(m-r+1)_{m-r+1}}=1.
\end{align*}
For identity \eqref{eq:Conj37withmu}, we proceed to simplify the ratio of
its left-hand side divided by its right-hand side:
\begin{equation*}
  \frac{{2^{2m-2r}} \bigl(\frac{\mu}{2}+r-\frac{1}{2}\bigr)_{m-r} \,
    \bigl(\frac{\mu}{2}+r\bigr)_{m-r}}{\bigl(\mu+2r-1\bigr)_{2m-2}} \cdot
  \frac{\prod_{i=1}^{2r-2} \bigl(\mu+i-1\bigr)_{2m+2r-2i-1}}{
    \prod_{i=0}^{2r-3}\bigl(\mu+i\bigr)_{2m-1}}=1,
\end{equation*}
where the equality was obtained by first combining the big products and
applying \ref{itm:pochinterlace} to simplify the big rational factor in front
of it, and then applying \ref{itm:pochshift} for another simplification that
enabled us to see that the factors can cancel.
\end{proof}

The following theorem gives a closed form for $\D{2r}1{\mu}{2m}$ and thereby
resolves \cite[Conjecture~20]{KoutschanThanatipanonda19}.  Note that the
result in~\cite{KoutschanThanatipanonda19} is stated in a slightly different,
but equivalent form, which can be verified by a routine calculation.

\begin{theorem}\lbl{thm:KTConj20}
Let $\mu$ be an indeterminate and $m,r\in\Z$. If $m\geq r\geq1$, then
\begin{multline*}
  \D{2r}1{\mu}{2m} = \\
  \frac{(-1)^{m-r} \, (\mu-1) \, \bigl(\mu+2r\bigr)_{2m-1}}{(2r-1)! \,
    \bigl(m+r\bigr)_{m-r+1} \, \bigl(\frac{\mu}{2}+r+\frac{1}{2}\bigr)_{m-r}} \cdot
  \prod_{i=1}^{m-r} \frac{
    \bigl(\mu+2i+6r-2\bigr)_{i-1}^2 \, \bigl(\frac{\mu}{2}+2i+3r-\frac{1}{2}\bigr)_{i}^2
  }{
    \bigl(i\bigr)_{i}^2 \, \bigl(\frac{\mu}{2}+i+3r-\frac{1}{2}\bigr)_{i-1}^2
  }.
\end{multline*}
\end{theorem}
\begin{proof}
  We employ the first equation of \lem{biglemma1} to connect this determinant to
  \thm{Krat37nice}:
  \[
    \D{2r}1{\mu}{2m} = R_{2r,1}^{\mu}(2m) \cdot \E{2r-1}1{\mu+3}{2m-1}.
  \]
  We observe that the product in \thm{Krat37nice} turns into the above product
  via the substitution $\mu\to\mu+3$, and the prefactor from \thm{Krat37nice}
  combines nicely with the rational function $R_{2r,1}^{\mu}(2m)$ to yield
  the prefactor in the claimed formula.
\end{proof}


\section{Closed Forms for \texorpdfstring{$\bm{\E{-1}{2r-1}{\mu}{2m-1}}$}{E(-1,2r-1)(2m-1)}
  and \texorpdfstring{$\bm{\D{-1}{2r}{\mu}{2m}}$}{D(-1,2r)(2m)}}
\lbl{sec:ktconj21}

In this section, we derive closed forms for the determinants
$E_{-1,2r-1}^{\mu}(2m-1)$ (\thm{Eneg1CF}) and $D_{-1,2r}^{\mu}(2m)$ (\thm{ktconj21}),
which allows us to
resolve \cite[Conjecture~21]{KoutschanThanatipanonda19} and give its
$E$-analog. The roadmap for how we do this can be seen in \fig{3Dfamilies} in
the color \cEDnegone. We tried to parallel \sect{krat37ktconj20} by introducing
a key lemma in order to establish a relationship between the two families that we
want closed forms for. However, we encountered serious problems with this strategy.
In the previous section, not only
did the families exhibit a simple ratio, but we were also able to make the
correct adjustments to the matrix (by multiplying by the elementary
matrices~$\matnot{L}_n$ and~$\matnot{R}_n$) to be able to apply the holonomic
ansatz. For the families in this section, the ratio was not just a
rational function with fixed numerator and denominator degrees (in~$\mu$), but
a quotient of Pochhammer symbols. As a consequence, we were unable to complete
the guessing step because the shape of the recurrences (namely, their
coefficient degrees) depended on the parameter~$r$, which prevented us from
finding recurrences with \emph{symbolic}~$r$.  In addition, we were unable to
make the necessary adjustments to modify our matrices to work with the
method. However, we were able to make it work \textit{after} switching the
parameters $s$ and~$t$, with the added bonus that the resulting ratio turned
out to be similarly simple as the one in \lem{biglemma1}!

The catch is that on its own, the switching of the parameters causes the
determinants to evaluate to zero, and this resulted in the ratio being of an
indeterminate form $\frac{0}{0}$.  So for this reason, we will introduce a new
parameter $\eps$ into the binomial coefficients to counteract the bad behavior
and then take the limit as $\eps\to 0$ to get the result. In particular, we
use \defn{newbc} to write $\binom{x+2\eps}{k+\eps}$ as a Taylor series in
$\eps$ around $\eps=0$ for integers $k<0$ to get
\begin{equation}\lbl{eq:bctaylor2eps}
\binom{x+2\eps}{k+\eps}=(-1)^{k+1}\cdot\frac{(-k-1)!}{(x+1)_{-k}}\cdot\eps+O(\eps^2),\\
\end{equation}
where the first (constant) term is zero and the coefficient of the $\eps$-term
is computed by exploiting the properties of the logarithmic derivative of
$\Gamma(z)$ \cite[5.2.2]{DLMF} to get the derivative:
\begin{align*}
  \frac{d}{d\eps}\binom{x+2\eps}{k+\eps} &=
  \frac{\Gamma(x+2\eps+1)}{\Gamma(k+\eps+1)\,\Gamma(x-k+\eps+1)} \\
  &\quad\times\bigl(2\, \psi(x+2\eps+1)-\psi(k+\eps+1)-\psi(x-k+\eps+1)\bigr).
\end{align*}
Taking $\eps\to 0$, the first and third terms vanish, leaving us with
\[
  \lim\limits_{\eps\to 0}\,\frac{d}{d\eps}\binom{x+2\eps}{k+\eps} =
  -\frac{\Gamma(x+1)}{\Gamma(x-k+1)}\cdot
  \lim\limits_{\eps\to 0}\,\frac{\psi(k+\eps+1)}{\Gamma(k+\eps+1)} =
  -\frac{1}{(x+1)_{-k}}\cdot(-1)^k\,(-k-1)!,
\]
where we use the fact that $\Gamma(z)$ and $\psi(z)$ are meromorphic functions
with simple poles of residue $(-1)^n/n!$ and $-1$ (respectively) at $z=-n$ for
$n\in\N_0$. For integers $k\geq 0$, the first (constant) term of the Taylor
expansion of $\binom{x+2\eps}{k+\eps}$ is the usual binomial coefficient~$\binom{x}{k}$.

We can now summarize the steps of this section.
{\parskip=0pt
\begin{itemize}
\item Analogous to \sect{krat37ktconj20}, we have a key result in
  \lem{biglemma2}, where we establish the ratios between the families
  $\D{2r+\eps}{-1+\eps}\mu{2m}$ and $\E{2r-1+\eps}{-1+\eps}\mu{2m-1}$. The
  introduction of $\eps$ causes more theoretical difficulties than usual,
  and we show in detail how to overcome them.
\item Once the connection is established, we apply an extra step to be able to
  connect the base case $\E{1+\eps}{-1+\eps}\mu{2m+1}$ to the known determinant
  $\D 10{\mu+3}{2m-1}$ in \lem{quoED1}.
\item These two lemmas, together with \lem{switch} (``switching''), enable us
  to realize the first main result of this section, a closed form for
  $E_{-1,2r-1}^{\mu}(2m-1)$ in \thm{Eneg1CF}.
\item In a similar fashion, we deduce the second main result of this section,
a closed form for $D_{-1,2r}^{\mu}(2m)$ in \thm{ktconj21}.
\end{itemize}
}\medskip

\begin{lemma}\lbl{lem:biglemma2}
Let $\mu$ be an indeterminate, $\eps\in \R$, $m,r\in\Z$. If $m>r\geq 1$, then
\begin{align}
  \lbl{eq:biglemeq3}
  \lim\limits_{\eps\to 0}
  \left(\frac{\D{2r+\eps}{-1+\eps}{\mu}{2m}}{\E{2r-1+\eps}{-1+\eps}{\mu+3}{2m-1}}\right)
  &=\frac{2r(2m-1)(\mu-3)(\mu+2m+2r-2)}{\mu(m+r)(\mu+2m-3)(\mu+2r-2)},\\
  \lbl{eq:biglemeq4}
  \lim\limits_{\eps\to 0}
  \left(\frac{\E{2r+1+\eps}{-1+\eps}{\mu}{2m+1}}{\D{2r+\eps}{-1+\eps}{\mu+3}{2m}}\right)
  &=\frac{2m(2r+1)(\mu-3)(\mu+2m+2r)}{\mu(m+r+1)(\mu+2m-2)(\mu+2r-1)}.
\end{align}
Remarks: Equations~\eqref{eq:biglemeq3} and~\eqref{eq:biglemeq4} are expressed
with the extra parameter $\varepsilon$ because the ratios of their left-hand sides are of the
indeterminate form $\tfrac{0}{0}$ otherwise. In~\eqref{eq:biglemeq4}, $r$
could be~$0$.
\end{lemma}

\begin{proof}
We can see that both identities can be presented in a uniform way:
\[
  \lim\limits_{\eps\to 0}
  \left(\frac{\detnew A{s+\eps}{-1+\eps}{\mu}{n}}{\B{s-1+\eps}{-1+\eps}{\mu+3}{n-1}}\right)
  = \frac{2s(n-1)(\mu-3)(\mu+n+s-2)}{\mu(n+s)(\mu+n-3)(\mu+s-2)}
  =: R_{s,-1}^{\mu}(n),
\]
where $(A,B,s,n)=(D,E,2r,2m)$ or $(A,B,s,n)=(E,D,2r+1,2m+1)$. Like in the
proof of \lem{biglemma1}, we use an inductive argument to ensure that
$\lim_{\eps\to0}\bigl(\frac1\eps\detnew A{s+\eps}{-1+\eps}{\mu}{n}\bigr)$
exists and is nonzero. As a base case, we use
$\E{1+\eps}{-1+\eps}{\mu}{2m-2r+1}$ (justified by the fact that once \lem{quoED1} is established, we can use the knowledge that the determinant $D_{1,0}^{\mu+3}(2m-1)$ is nonzero), and as induction
hypothesis we assume from now on that
$\lim_{\eps\to0}\bigl(\frac1\eps\B{s-1+\eps}{-1+\eps}{\mu+3}{n-1}\bigr)$
exists and is nonzero.

Like in the proof of \lem{biglemma1}, we manipulate
$\Mat{A}{s+\eps}{-1+\eps}{\mu}{n}$ by multiplying with the elementary matrices
$\matnot{L}_n,\matnot{R}_n$ defined in~\eqref{eq:LandR} such that its
determinantal value remains unaffected, and by employing \lem{pascal} to
simplify the entries:
\begin{equation}\lbl{eq:tildeA}
\matnot{L}_n\cdot \Mat{A}{s+\eps}{-1+\eps}{\mu}{n}\cdot \matnot{R}_n =
\left(
\begin{array}{c:c}
  \binom{\mu+s-3+2\eps}{-1+\eps} &
  \binom{\mu+j+s-3+2\eps}{j-2+\eps}-\binom{\mu+s-3+2\eps}{-2+\eps}
  \pm\smash{\sum\limits_{k=1}^j\delta_{s,k-2}}
  \\
  & \scriptstyle(2\;\leq\;j\;\leq\;n)
  \\
  \hdashline & \\[-12pt]
  \binom{\mu+i+s-5+2\eps}{-2+\eps}
  & \binom{\mu+i+j+s-5+2\eps}{j-3+\eps}-\binom{\mu+i+s-5+2\eps}{-3+\eps}\mp\delta_{s,j-i}
  \\
  \scriptstyle(2\;\leq\;i\;\leq\;n)
  & \scriptstyle(2\;\leq\;i,j\;\leq\;n)
\end{array}
\right),
\end{equation}
where $\pm$ is $+$ if $A=D$ and $-$ if $A=E$ (and $\mp$ is $-$ if $A=D$ and
$+$ if $A=E$). Next, we delete the second binomial coefficient from each entry
in the second column of this matrix, that is, we add the vector
$C=\bigl(\binom{\mu+s-3+2\eps}{-2+\eps},
\binom{\mu+i+s-5+2\eps}{-3+\eps}_{2\leq i\leq n}\bigr){}^\mathrm{T}$ to the second
column. The resulting matrix is displayed here in a form where
we express all of its entries in terms of their Taylor expansions with respect to
the variable~$\eps$ (around $\eps=0$), using the formula
in~\eqref{eq:bctaylor2eps}, and by omitting lower-order terms:
\[
\tilde{\matnot{A}} :=
\left(\begin{array}{c:c:c}
  \frac{1}{\mu+s-2}\cdot\eps
  & 1
  & \binom{\mu+j+s-3}{j-2} \pm \sum\limits_{k=1}^j\delta_{s,k-2}
  \\
  & & \scriptstyle(3\;\leq\;j\;\leq\;n)
  \\
  \hdashline & & \\[-12pt]
  \frac{-1}{(\mu+i+s-4)_2}\cdot\eps
  & \frac{1}{\mu+i+s-2}\cdot\eps
  & \binom{\mu+i+j+s-5}{j-3} \mp \delta_{s,j-i}
  \\
  \scriptstyle(2\;\leq\;i\;\leq\;n)
  & \scriptstyle(2\;\leq\;i\;\leq\;n)
  & \scriptstyle(2\;\leq\;i\;\leq\;n,\; 3\;\leq\;j\;\leq\;n)
\end{array}
\right) =:
\left(\begin{array}{c:c:c}
  \tilde{a}_{1,1}\cdot\eps & 1 & \tilde{a}_{1,j} \\[4pt]
  \hdashline & & \\[-14pt]
  \tilde{a}_{i,1}\cdot\eps & \tilde{a}_{i,2}\cdot\eps & \tilde{a}_{i,j}
\end{array}
\right).
\]
We let $\tilde{a}_{i,j}$ denote the first nonzero coefficient
in the Taylor expansion of the ($i,j)$-entry of~$\tilde{\matnot{A}}$.
Now imagine that the second column of~$\tilde{\matnot{A}}$ gets replaced by
the vector~$C$: the determinant of the resulting matrix is $O(\eps^2)$ because
all of the entries in its first two columns are $O(\eps)$. Hence,
\begin{equation}\lbl{eq:detAAtilde}
  \detnew A{s+\eps}{-1+\eps}{\mu}{n} =
  \det\bigl(\matnot{L}_n\cdot \Mat{A}{s+\eps}{-1+\eps}{\mu}{n}\cdot \matnot{R}_n\bigr) =
  \det\bigl(\tilde{\matnot{A}}\bigr) + O(\eps^2)
\end{equation}
by the linearity of the determinant in its columns.  We choose to compute the determinant
of~$\tilde{\matnot{A}}$ by expanding along the first column:
\begin{equation}\lbl{eq:LaplaceEps}
  \det\bigl(\tilde{\matnot{A}}\bigr) =
  \sum_{i=1}^n \tilde{a}_{i,1}\cdot\eps\cdot\Cof i1n,
\end{equation}
where $\Cof i1n$ are the corresponding cofactors
from~$\tilde{\matnot{A}}$. By noting that the lower-right
$(n-1)\times(n-1)$-submatrix of~$\tilde{\matnot{A}}$ is equal to the matrix
$\Mat B{s-1+\eps}{-1+\eps}{\mu+3}{n-1}$ (after the omission of lower-order
terms), we see that
\begin{equation}\lbl{eq:limBCof1}
  \lim\limits_{\eps\to 0}\frac{\Cof 11n}{\B{s-1+\eps}{-1+\eps}{\mu+3}{n-1}}=1.
\end{equation}
Our induction hypothesis tells us that
$\lim_{\eps\to0}\bigl(\frac{1}{\eps}\Cof 11n\bigr)$ exists and is nonzero.
By defining
\begin{equation}\lbl{eq:defcni}
  c_{n,i} := \lim_{\eps\to0}\frac{\eps\cdot\Cof i1n}{\Cof 11n}
\end{equation}
and by using \eqref{eq:detAAtilde}, \eqref{eq:limBCof1},
and~\eqref{eq:LaplaceEps}, we can express our desired quotient of determinants
in terms of these quantities:
\begin{align*}
  \frac{\detnew A{s+\eps}{-1+\eps}{\mu}{n}}{\B{s-1+\eps}{-1+\eps}{\mu+3}{n-1}}
  = \frac{\det\bigl(\tilde{\matnot{A}}\bigr)}{\Cof11n}+O(\eps)
  &= \tilde{a}_{1,1}\cdot\eps +
    \sum_{i=2}^n \tilde{a}_{i,1}\cdot\bigl(c_{n,i}+O(\eps)\bigr) + O(\eps) \\
  &= \sum_{i=2}^n \tilde{a}_{i,1}\cdot c_{n,i} + O(\eps).
\end{align*}
Like in the proof of \lem{biglemma1}, we aim at characterizing the $c_{n,i}$
as the unique solution of a certain linear system.  If we replace the first
column of~$\tilde{\matnot{A}}$ by its second column, then the corresponding
determinant is~$0$. By modifying \eqref{eq:LaplaceEps} accordingly, we obtain
\[
  0 = 1\cdot\Cof11n + \sum_{i=2}^n\tilde{a}_{i,2}\cdot\Cof i1n,
\]
which, after dividing by $\Cof11n$, turns into
\begin{equation}\lbl{eq:syseps1}
  0 = 1 + \sum_{i=2}^n \tilde{a}_{i,2}\cdot\bigl(c_{n,i}+O(\eps)\bigr)
  = 1 + \sum_{i=2}^n \tilde{a}_{i,2}\cdot c_{n,i} + O(\eps).
\end{equation}
Similarly, for $3\leq j\leq n$, we replace the first column
of~$\tilde{\matnot{A}}$ by its $j$-th column to get
\[
  0 = \tilde{a}_{1,j}\cdot\Cof11n + \sum_{i=2}^n\tilde{a}_{i,j}\cdot\Cof i1n,
\]
which, after dividing both sides by $\frac{1}{\eps}\Cof11n$, turns into
\begin{equation}\lbl{eq:syseps2}
  0 = \eps\cdot\tilde{a}_{1,j} +
  \sum_{i=2}^n\tilde{a}_{i,j}\cdot\bigl(c_{n,i}+O(\eps)\bigr)
  = \sum_{i=2}^n \tilde{a}_{i,j}\cdot c_{n,i} + O(\eps).
\end{equation}
By our induction hypothesis, the matrix $\bigl(\tilde{a}_{i,j}\bigr){}_{2\leq i,j\leq n}$
has full rank. Hence the system \eqref{eq:syseps1}, \eqref{eq:syseps2}, after
removing the unnecessary $O(\eps)$ terms, has a unique solution
$(c_{n,2},\dots,c_{n,n})$ for all $n>s$.

We proceed now in the usual way, as described in \sect{holonomicansatz}:
compute the $c_{n,i}$ explicitly for several fixed $n$ and~$i$, guess/interpolate
recurrences that are satisfied by this data, view these recurrences as an
implicit definition of some bivariate sequence, show that this sequence
satisfies the characterizing linear system and therefore agrees with the
$c_{n,i}$ as defined in~\eqref{eq:defcni}, and finally use it to obtain the
desired quotient of determinants. More explicitly, we employ the holonomic
framework to prove the following three identities:
\begin{alignat}{2}
  &\sum_{i=2}^{n}\frac{1}{\mu+i+s-2}\cdot c_{n,i} &&= -1, \notag \\
  &\sum_{i=2}^{n}\binom{\mu+i+j+s-5}{j-3}\cdot c_{n,i} &&= \pm c_{n,j-s},
  \qquad (3\leq j \leq n), \notag \\
  \lbl{eq:bl2eq3}
  &\sum_{i=2}^{n}\frac{-1}{(\mu+i+s-4)_2}\cdot c_{n,i} &&= R_{s,-1}^{\mu}(n),
\end{alignat}
where $c_{n,j-s}=0$ for $j\leq s$.
The computations for these identities turned out to be very similar to the
computations for the identities in \lem{biglemma1} so we will not repeat the
exposition. All of the computational details can be found in the accompanying
electronic material \cite{EM2}. However, we remark that the third identities
were much easier as there are no singularities in the certificates. Thus, the
annihilating ideal for the summation could be directly read off and certified
from the computation without further adjustments.

We can hence conclude that \eqref{eq:biglemeq3} and \eqref{eq:biglemeq4} hold,
which also completes our induction step.
\end{proof}

We now connect our base case to another determinant that will enable us to
prove the main theorems.

\begin{lemma}\lbl{lem:quoED1}
Let $\mu$ be an indeterminate, $\eps\in\R\setminus \lbrace 0\rbrace$,
and $m\in\N$. Then
\[
  \lim\limits_{\eps\to 0}\left(
  \frac{\E{1+\eps}{-1+\eps}{\mu}{2m+1}}{\eps\cdot \D10{\mu+3}{2m-1}}\right)
  =-\frac{(4m-2)(\mu-3)(\mu+2m+1)}{(m+1)(\mu-1)(\mu+1)(\mu+3)(\mu+2m-2)}.
\]
\end{lemma}
 
\begin{proof}
We would like to adapt the holonomic ansatz (see \sect{holonomicansatz}) to
confirm the claimed identity. First we do some basic row and column operations
for $\Mat{E}{1+\eps}{-1+\eps}{\mu}{2m+1}$ by multiplying with the elementary
matrices $\matnot{L}_{2m+1}$ from~\eqref{eq:LandR} and 
\[
\tilde{\matnot{R}}_{2m+1}:=\begin{pmatrix} 
0&-1&0&0&0&\cdots \\ 
1&0&1&1&1&\cdots \\ 
0&0&1&1&1&\cdots \\ 
0&0&0&1&1&\cdots \\ 
\vdots & \vdots & \vdots &  & \ddots & \ddots \end{pmatrix},\
\]
and then apply \lem{pascal} and the Taylor expansion~\eqref{eq:bctaylor2eps}
such that the transformed matrix
$\matnot{L}_{2m+1} \cdot \Mat{E}{1+\eps}{-1+\eps}{\mu}{2m+1}\cdot \tilde{\matnot{R}}_{2m+1}$
becomes
\begin{equation}\lbl{eq:matepsED1}
  \left(
  \begin{array}{c:c:c}
  1+O(\eps) & \frac{\eps}{1-\mu}+O(\eps^2) & \binom{\mu+j-2}{j-2}-1 + O(\eps) \\
  & & \scriptstyle (3\;\leq\;j\;\leq\;2m+1) \\
  \hdashline & & \\[-12pt]
  \frac{\eps}{\mu}+O(\eps^2) & \frac{\eps}{(\mu-1)_2}+O(\eps^2) &
  \binom{\mu+j-2}{j-3}+\delta_{1,j-2}+ O(\eps) \\
  & & \scriptstyle(3\;\leq\;j\;\leq\;2m+1) \\
  \hdashline & & \\[-12pt]
  \begin{array}{c} \frac{\eps}{\mu+i-2}+O(\eps^2) \\ \null \end{array} &
  \begin{array}{c}
    \frac{\eps}{(\mu+i-3)_2}+O(\eps^2) \\ \scriptstyle(3\;\leq\;i\;\leq\;2m+1)
  \end{array} &
  \Mat{D}{1}{0}{\mu+3}{2m-1}+ O(\eps)
\end{array}
\right).
\end{equation}
Note that the determinantal value remained unaffected under this
transformation, and that now the matrix $\Mat{D}10{\mu+3}{2m-1}$ appears as
a submatrix (the $O(\eps)$ added to this matrix means that it is added to
every entry). Since the determinant behaves like a linear function in the
columns of the matrix, the determinant of~\eqref{eq:matepsED1}
is equal to $\eps\cdot\tilde{E}+O(\eps^2)$, where
\[
  \tilde{E} := \det
\left(
\begin{array}{c:c:c}
  \mute{1} & \mute{\frac{1}{1-\mu}} &
  \mute{\binom{\mu+j-2}{j-2}-1} \\
  & & \scriptstyle\mute{(3\;\leq\;j\;\leq\;2m+1)} \\
  \hdashline & & \\[-12pt]
  \mute{0} & \frac{1}{(\mu-1)_2} & \binom{\mu+j-2}{j-3}+\delta_{1,j-2} \\
  & & \scriptstyle(3\;\leq\;j\;\leq\;2m+1) \\
  \hdashline & & \\[-12pt]
  \begin{array}{c} \mute{0} \\ \null \end{array} &
  \begin{array}{c}
    \frac{1}{(\mu+i-3)_2} \\ \scriptstyle(3\;\leq\;i\;\leq\;2m+1)
  \end{array} &
  \Mat{D}{1}{0}{\mu+3}{2m-1}
\end{array}
\right).
\] 
Then denote by $\tilde{\matnot{E}}$ the bottom right $2m \times 2m$ submatrix
of the above matrix, whose determinant also
equals~$\tilde{E}$. Since $\D10{\mu+3}{2m-1}$ is nonzero
by~\cite[Proposition~9]{KoutschanThanatipanonda19}, we have that
\[
  \lim_{\eps\to 0}
  \left(\frac{\E{1+\eps}{-1+\eps}{\mu}{2m+1}}{\eps\cdot \D10{\mu+3}{2m-1}}\right)
  =\lim_{\eps\to 0}\left(\frac{\eps\cdot\tilde{E}+O(\eps^2)}{\eps\cdot\D10{\mu+3}{2m-1}}\right)
  =\frac{\tilde{E}}{\D10{\mu+3}{2m-1}}.
\]
It is now sufficient to apply the holonomic ansatz argument to
$\tilde{\matnot{E}}$ (for which $\Mat{D}10{\mu+3}{2m-1}$ is its
bottom right submatrix) and to expand along its first column:
\[
  \tilde{E}=\frac{1}{(\mu-1)_2}\cdot \Cof11{2m}+\ldots+
  \frac{1}{(\mu+2m-2)_2}\cdot \Cof{2m}1{2m},
\]
where $\Cof i1{2m}$ denotes the corresponding cofactor. Define
\begin{equation}\lbl{eq:cquoED1}
  c_{2m,i}:=\frac{\Cof i1{2m}}{\Cof11{2m}},
\end{equation}
and note that $\Cof11{2m}=\D10{\mu+3}{2m-1}$, then the assertion will be
confirmed provided that we can show that for all $m\geq 1$:
\begin{equation}\lbl{eq:sys3eq3}
  \sum\limits_{i=1}^{2m} \frac{c_{2m,i}}{(\mu+i-2)_2} = 
  -\frac{(4m-2)(\mu-3)(\mu+2m+1)}{(m+1)(\mu-1)(\mu+1)(\mu+3)(\mu+2m-2)}.
\end{equation}
Like before, we note that for each fixed $m$, $(c_{2m,1},\ldots,c_{2m,2m})$
satisfy the system of equations
\begin{equation}\lbl{eq:sys3}
\begin{cases}
  \; c_{2m,1}=1, & m\geq 1,\\[1ex]
  \; \displaystyle\sum_{i=1}^{2m} c_{2m,i}\cdot
  \left( \binom{\mu+i+j-2}{j-2}+\delta_{i,j-1}\right) =0, & 2\leq j \leq 2m,
\end{cases}
\end{equation}
and that this solution is unique since $\Mat{D}10{\mu+3}{2m-1}$ has full rank.
We use \eqref{eq:sys3} to generate data, and then guess recurrences in order
to re-define $c_{2m,i}$ as their solution with suitable initial values.  In
the end, proving the lemma reduces to confirming the
three identities corresponding to \eqref{eq:sys3} and \eqref{eq:sys3eq3} for all
$m\geq 1$ by the same method as in \lem{biglemma1}.
The advantage here of course is that we have one parameter less 
(no~$r$). However, we did encounter a singularity in the certificates at
$i=2m+1$ for both summations, which needed to be treated with some additional
adjustments as in \lem{biglemma1}.
\end{proof}

\begin{remark}
The proof of \lem{quoED1} shows how a slight modification in column
operations can produce a new setting in which the holonomic ansatz applies. In
particular, using the matrix $\tilde{\matnot{R}}$ leads to the appearance of
$\Mat D10{\mu+3}{2m-1}$, whose determinant appears in the statement of the
lemma. We note that this matrix $\tilde{\matnot{R}}$ can also be used to give
an alternative proof to \lem{biglemma2}, in which the multiplication
\eqref{eq:tildeA} (with $\tilde{\matnot{R}}$ replacing $\matnot{R}$)
introduces a new submatrix $\Mat B{s-1}0{\mu+3}{n-1}$ that is not exactly the
denominator in the statement of the lemma, but related to it in such a way that the
holonomic ansatz can be adapted without needing to deal with $\eps$ explicitly
(similarly to the second part of the proof of \lem{quoED1}). 
However, this method introduces a new summation on the right-hand side
of~\eqref{eq:bl2eq3}, and does not improve the overall computational time. The
reader can find this alternative proof of \lem{biglemma2} in \ref{ap:altproof}.
\end{remark}

\begin{theorem}
\lbl{thm:Eneg1CF}
Let $\mu$ be an indeterminate and $m,r\in\Z$. If $m\geq r\geq1$, then
\begin{align*}
  E_{-1,2r-1}^\mu(2m+1) &=
  \frac{(-1)^{m-r} \, (3-\mu) \, (m+r+1)_{m-r}}{
    2^{2m-2r+1} \, \bigl(\frac{\mu}{2}+r-\frac32\bigr)_{m-r+1}} \cdot
  \prod_{i=1}^{2m} \frac{(\mu+i-3)_{2r}}{(i)_{2r}} \\
  &\quad\times \prod_{i=1}^{m-r} \frac{\bigl(\mu+2i+6r-3\bigr)_i^2 \,
    \bigl(\frac{\mu}{2}+2i+3r-1\bigr)_{i-1}^2}{\bigl(i\bigr)_i^2 \,
    \bigl(\frac{\mu}{2}+i+3r-1\bigr)_{i-1}^2}.
\end{align*}
\end{theorem}

\begin{proof}
By applying \lem{biglemma2} $(2r-2)$ times we get the relation
\begin{align*}
  & \lim\limits_{\eps\to 0}\biggl(\frac{\E{2r-1+\eps}{-1+\eps}{\mu}{2m-1}}{
    \E{1+\eps}{-1+\eps}{\mu+6r-6}{2m-2r+1}}\biggr)
  \\
  &=\lim\limits_{\eps\to 0}\biggl(
  \frac{\E{2r-1+\eps}{-1+\eps}{\mu}{2m-1}}{\D{2r-2+\eps}{-1+\eps}{\mu+3}{2m-2}}\cdot
  \frac{\D{2r-2+\eps}{-1+\eps}{\mu+3}{2m-2}}{\E{2r-3+\eps}{-1+\eps}{\mu+6}{2m-3}}
  \cdots\frac{\D{2+\eps}{-1+\eps}{\mu+6r-9}{2m-2r+2}}{\E{1+\eps}{-1+\eps}{\mu+6r-6}{2m-2r+1}}
  \biggr)
  \\[1ex]
  &= R_{2r-1,-1}^{\mu}(2m-1) \cdot R_{2r-2,-1}^{\mu+3}(2m-2) \cdots R_{2,-1}^{\mu+6r-9}(2m-2r+2) \\
  &= \prod_{i=0}^{2r-3} R_{2r-1-i,-1}^{\mu+3i}(2m-1-i)=:P_r(m).
\end{align*}
Then
\begin{equation}\lbl{eq:eneg1ratio}
  \lim_{\eps\to 0}\biggl(\frac{E_{2r-1+\eps,-1+\eps}^\mu(2m+1)}{E_{2r-1+\eps,-1+\eps}^\mu(2m-1)}\biggr)
  = \frac{P_r(m+1)}{P_r(m)}\cdot
  \lim_{\eps\to 0}\biggl(\frac{\E{1+\eps}{-1+\eps}{\mu+6r-6}{2m-2r+3}}{
    \E{1+\eps}{-1+\eps}{\mu+6r-6}{2m-2r+1}}\biggr).
\end{equation}
By applying \lem{quoED1} and some simplifications of the ratio of products,
we note that the right-hand side of~\eqref{eq:eneg1ratio} is equal to
\pagebreak[1]
\[
  \frac{m(2m-1)(\mu+2m-4)(\mu+2m+4r-3)}
  {(m+r)(2m-2r-1)(\mu+2m+2r-3)(\mu+2m+2r-4)} \cdot
  \frac{D_{1,0}^{\mu+6r-3}\big(2m-2r+1\big)}{D_{1,0}^{\mu+6r-3}\big(2m-2r-1\big)}.
\]
According to \cite[Proposition~9]{KoutschanThanatipanonda19}, the ratio of
$D$'s can be reduced so that the previous expression simplifies to
\[
  \frac{-m (2m-1) (\mu+2m-4) \, \bigl(\mu+2m+4r-3\bigr)_{m-r}^2 \,
    \bigl(\frac{\mu}{2}+2m+r-1\bigr)_{m-r-1}^2
  }{
    (m+r) (\mu+2m+2r-3) (\mu+2m+2r-4) \, \bigl(m-r\bigr)_{m-r}^2 \,
    \bigl(\frac{\mu}{2}+m+2r-1\bigr)_{m-r-1}^2}.
\]
Now we apply \lem{switch} to get the form that we want, in other words:
\begin{align*}
  \frac{E_{-1,2r-1}^\mu(2m+1)}{E_{-1,2r-1}^\mu(2m-1)}
  &= \lim_{\eps\to 0}\biggl(\frac{E_{-1+\eps,2r-1+\eps}^\mu(2m+1)}{E_{-1+\eps,2r-1+\eps}^\mu(2m-1)}\biggr) \\
  &= \lim_{\eps\to 0}\left(\frac{E_{2r-1+\eps,-1+\eps}^\mu(2m+1)}{E_{2r-1+\eps,-1+\eps}^\mu(2m-1)} \cdot
    \prod_{i=0}^{2r-1}\frac{\bigl(\mu+i+\eps-2\bigr)_{2m+1} \, \bigl(i+\eps\bigr)_{2m-1}}{
      \bigl(i+\eps\bigr)_{2m+1} \, \bigl(\mu+i+\eps-2\bigr)_{2m-1}} \right).
\end{align*}
For the limit of the product we obtain, after applying \ref{itm:pochshift},
\[
  \lim_{\eps\to 0}
  \prod_{i=0}^{2r-1} \frac{(\mu+i+\eps+2m-3)(\mu+i+\eps+2m-3)}{
    (i+\eps+2m-1)(i+\eps+2m)} =
  \frac{(\mu+2m-3)_{2r} \, (\mu+2m-2)_{2r}}{(2m-1)_{2r} \, (2m)_{2r}}.
\]
Furthermore, we also have that
\begin{align*}
  &\frac{m (2m-1) (\mu+2m-4)}{(m+r) (\mu+2m+2r-3) (\mu+2m+2r-4)} \cdot
  \frac{(\mu+2m-3)_{2r} \,  (\mu+2m-2)_{2r}}{(2m-1)_{2r} \, (2m)_{2r}} \\
  &=\frac{(\mu+2m-4)_{2r} \, (\mu+2m-2)_{2r-1}}{(2m)_{2r-1}\, (2m+1)_{2r}}.
\end{align*}
Therefore, we can express $\E{-1}{2r-1}\mu{2m+1}/\E{-1}{2r-1}\mu{2r+1}$ in the form:
\[
  \prod_{i=r+1}^{m} \left(
  -\frac{\bigl(\mu+2i-4\bigr)_{2r} \, \bigr(\mu+2i-2\bigr)_{2r-1} \,
    \bigl(\mu+2i+4r-3\bigr)_{i-r}^2 \, \bigl(\frac{\mu}{2}+2i+r-1\bigr)_{i-r-1}^2
  }{
    \bigl(2i\bigr)_{2r-1} \, \bigl(2i+1\bigr)_{2r} \, \bigl(i-r\bigr)^2_{i-r} \,
    \bigl(\frac{\mu}{2}+i+2r-1\bigr)_{i-r-1}^2}
  \right).
\]
By the sums-of-minors formula \eqref{eq:sumofminors} and using \prop{detwithnoKD}
we get
\begin{align*}
  \E{-1}{2r-1}\mu{2r+1} &=
  \det_{\genfrac{}{}{0pt}{}{1\leq i\leq 2r+1}{1\leq j\leq 2r+1}} \binom{\mu+i+j+2r-6}{j+2r-2}
  - \det_{\genfrac{}{}{0pt}{}{1\leq i\leq 2r}{2\leq j\leq 2r+1}} \binom{\mu+i+j+2r-6}{j+2r-2} \\
  &= \prod_{i=1}^{2r-1} \frac{(\mu+i-3)_{2r+1}}{(i)_{2r+1}}
  - \prod_{i=1}^{2r} \frac{(\mu+i-3)_{2r}}{(i)_{2r}}
  = \frac{3-\mu}{2r} \cdot \prod_{i=1}^{2r-1} \frac{(\mu+i-3)_{2r+1}}{(i)_{2r+1}},
\end{align*}
because the Kronecker delta affects only the $(2r+1,1)$-entry of the matrix.
We rewrite the products
\begin{align*}
  &\qquad \frac{3-\mu}{2r} \cdot
  \left(\prod_{i=1}^{2r-1} \frac{(\mu+i-3)_{2r+1}}{(i)_{2r+1}}\right) \cdot
  \prod_{i=r+1}^{m} \frac{\bigl(\mu+2i-4\bigr)_{2r} \, \bigr(\mu+2i-2\bigr)_{2r-1}}{
    \bigl(2i\bigr)_{2r-1} \, \bigl(2i+1\bigr)_{2r}} \\
  &= \frac{3-\mu}{\mu+2r-3} \cdot
  \left(\prod_{i=1}^{2r} \frac{(\mu+i-3)_{2r}}{(i)_{2r}}\right) \cdot
  \prod_{i=r+1}^{m} \frac{i \, (2i-1) \, \bigl(\mu+2i-4\bigr)_{2r} \, \bigr(\mu+2i-3\bigr)_{2r}}{
    (i+r) \, (\mu+2i-3) \, \bigl(2i-1\bigr)_{2r} \, \bigl(2i\bigr)_{2r}} \\
  &= \frac{3-\mu}{\mu+2r-3} \cdot
  \left(\prod_{i=1}^{2m} \frac{(\mu+i-3)_{2r}}{(i)_{2r}}\right) \cdot
  \prod_{i=r+1}^{m} \frac{i \, (2i-1)}{(i+r) \, (\mu+2i-3)},
\end{align*}
so that we arrive at the claimed formula, after putting everything together
and performing some final Pochhammer simplifications.
\end{proof}

\begin{theorem}\lbl{thm:ktconj21}
Let $\mu$ be an indeterminate and $m,r\in\Z$. If $m>r\geq0$, then
\begin{align*}
  \D{-1}{2r}\mu{2m} &= 
  \frac{(-1)^{m-r} \, (\mu-3) \, \bigl(\frac{\mu }{2}+r-\frac{1}{2}\bigr)_{m-r-1}}{(2r+1)_{m-r}}
  \cdot\prod_{i=1}^{2m} \frac{(\mu+i-3)_{2r}}{(i)_{2r}} \\
  &\quad\times\prod_{i=1}^{m-r-1}
   \frac{\bigl(\mu+2i+6r\bigr)_i^2 \, \bigl(\frac{\mu}{2}+2i+3r+\frac{1}{2}\bigr)_{i-1}^2}{
     \bigl(i\bigr)_i^2 \, \bigl(\frac{\mu}{2}+i+3r+\frac{1}{2}\bigr)_{i-1}^2}.
\end{align*}
Remark: This proves~\cite[Conjecture 21]{KoutschanThanatipanonda19}.
\end{theorem}
\begin{proof}
We use~\eqref{eq:biglemeq4} to get
\begin{align*}
  \lim_{\eps\to 0}
  \biggl(\frac{D_{2r+\eps,-1+\eps}^\mu(2m+2)}{D_{2r+\eps,-1+\eps}^\mu(2m)}\biggr) &=
  \lim_{\eps\to 0}
  \biggl(\frac{E_{2r+1+\eps,-1+\eps}^{\mu-3}(2m+3)}{E_{2r+1+\eps,-1+\eps}^{\mu-3}(2m+1)}\biggr) \\
  &\quad\times\frac{m(m+r+2)(\mu+2m-3)(\mu+2m+2r-3)}{(m+1)(m+r+1)(\mu+2m-5)(\mu+2m+2r-1)}.
\end{align*}
Furthermore, using \lem{biglemma2} as in~\eqref{eq:eneg1ratio}, and then~\lem{quoED1},
the right-hand side of the previous equation can be simplified to
\[
  -\frac{m(2m+1)(\mu+2m-3) \, \bigl(\mu+2m+4r\bigr)_{m-r}^2 \,
    \bigl(\frac{\mu}{2}+2m+r+\frac{1}{2}\bigr)_{m-r-1}^2}{
    (m+r+1)(\mu+2m+2r-2)(\mu+2m+2r-1) \, \bigl(m-r\bigr)^2_{m-r} \,
    \bigl(\frac{\mu}{2}+m+2r+\frac{1}{2}\bigr)_{m-r-1}^2}.
\]
After applying \lem{switch}, we obtain for the ratio
\[
  \frac{\D{-1}{2r}\mu{2m+2}}{\D{-1}{2r}\mu{2m}} =
  -\frac{\bigl(\mu+2m-3\bigr)_{2r+1} \, \bigl(\mu+2m-1\bigr)_{2r} \,
    \bigl(\mu+2m+4r)_{m-r}^2 \, \bigl(\frac{\mu}{2}+2m+r+\frac{1}{2}\bigr)_{m-r-1}^2
  }{
    \bigl(2m+1\bigr)_{2r} \, \bigl(2m+2\bigr)_{2r+1} \, \bigl(m-r\bigr)_{m-r}^2 \,
    \bigl(\frac{\mu}{2}+m+2r+\frac{1}{2}\bigr)_{m-r-1}^2},
\]
which is exactly the form stated in \cite[Conjecture 21]{KoutschanThanatipanonda19}.
Continuing in an analogous way as in the proof of \thm{Eneg1CF}, we end up with
the claimed formula.
\end{proof}


\section{New Relationships Between the Families}
\lbl{sec:misc}

In the previous two sections, we had a very clear goal in mind: we wanted to
obtain closed forms for certain determinants. This was accomplished by
\textit{first} recognizing that the determinants we wanted closed forms for
already had a relationship (that is, their ratios or the limit of their ratios
are equal to nice rational functions), and \textit{then} exploiting those
relationships (e.g., \lem{biglemma1} and \lem{biglemma2}) to get what we want.
In this section, we explore the opposite direction. Can we find nice
relationships between determinants whose closed forms we already know?
Moreover, could these relationships be used to understand better what is
happening from a combinatorics perspective? 

As a first example, we can show a connection between two determinants of different types: one
which exactly counts cyclically symmetric rhombus tilings, and one which
performs a weighted count. However, note that whatever values are substituted
for the parameters, at least one of the two determinants does not allow for a
combinatorial interpretation, because the central hole is larger than the
whole hexagon. This combinatorial reciprocity was first observed in
\cite[Conjecture~24]{KoutschanThanatipanonda19} in the case of
$D_{s,t}^{\mu}(n)$.  The following theorem resolves this conjecture, using
relations between the $D$- and $E$-determinants that were derived above, and
in addition states an analogous formula for $E_{s,t}^{\mu}(n)$. To visualize the
idea of this proof, the roadmap is highlighted in \fig{3Dfamilies} with the
color \cEDzeroref.

\begin{theorem}\lbl{thm:KTConj24}
  Let $\mu$ be an indeterminate and $m,r\in\Z$ such that $m\geq r \geq 1$.
  Then
  \begin{align*}
    D_{2r-1,0}^{\mu}(2m+1) &= D_{0,0}^{1-\mu-6m}(2m-2r+2), \\
    E_{2r-1,0}^{\mu}(2m+1) &= E_{0,0}^{1-\mu-6m}(2m-2r+2).
  \end{align*}
\end{theorem}

\begin{proof}
Applying \lem{famA} to the left-hand sides of both identities $2r-1$ times, we
can reduce the problem to confirming
\begin{align}
  \lbl{E00equalsD00}
  E_{0,0}^{\mu+6r-3}(2m-2r+2) &= D_{0,0}^{1-\mu-6m}(2m-2r+2), \\
  \notag
  D_{0,0}^{\mu+6r-3}(2m-2r+2) &= E_{0,0}^{1-\mu-6m}(2m-2r+2).
\end{align}
Note that the above two identities are equivalent via the
substitution $\mu\to4-\mu-6r-6m$. So it is enough to prove the first identity.
Instantiating \cor{Es0CF} and rewriting, we obtain
\[
  E_{0,0}^{\mu}(2m) =
  4\cdot\prod_{i=0}^{m-1}
  \frac{-(\mu+2i-1)\,\bigl(\mu+2i\bigr)_i^2\,\bigl(\frac{\mu}{2}+2i\bigr)_{i+1}^2}
  {4(2i+1) \, \bigl(i\bigr)_i^2 \, \bigl(\frac{\mu}{2}+i\bigr)_{i+1}^2},
\]
whose right-hand side, after substituting  $\mu\rightarrow \mu+6r-3, m\rightarrow m-r+1$, turns into
\begin{equation*}
  -(\mu+6r-4)\cdot\prod_{i=1}^{m-r}\frac{-(\mu+2i+6r-4) \, \bigl(\mu+2i+6r-3\bigr)_i^2 \,
    \bigl(\frac{\mu}{2}+2i+3r-\frac{3}{2}\bigr)_{i+1}^2}{4(2i+1) \, \bigl(i\bigr)_i^2 \,
    \bigl(\frac{\mu}{2}+i+3r-\frac{3}{2}\bigr)_{i+1}^2},
\end{equation*}
and simplifying the formula $2\cdot\prod_{i=1}^{2m-2r+1}R_{0,0}^{1-\mu-6m}(i)$ from
\cite[Proposition~8]{KoutschanThanatipanonda19} for the right-hand side of
\eqref{E00equalsD00} (taking into account the even/odd behavior) gives
\begin{align*}
  & \frac{2 \, \bigl(-\mu-4m-2r+1\bigr)_{m-r} \, \bigl(-\frac{\mu}{2}-m-2r+2\bigr)_{m-r+1}}
  {\bigl(m-r+1\bigr)_{m-r+1} \, \bigl(-\frac{\mu}{2}-2m-r+1\bigr)_{m-r}} \\
  & \times\prod\limits_{i=1}^{m-r}\frac{\bigl(-\mu+2i-6m+1\bigl)_i \,
    \bigl(-\frac{\mu}{2}+2i-3m+1\bigr)_{i-1} \, \bigl(-\mu+2i-6m-1\bigr)_{i-1} \,
    \bigl(-\frac{\mu}{2}+2i-3m\bigr)_i}{\bigl(i\bigr)_i^2 \,
    \bigl(-\frac{\mu}{2}+i-3m+1\bigr)_{i-1} \, \bigl(-\frac{\mu}{2}+i-3m\bigr)_{i-1}}.
\end{align*}
Then the ratio of the above two formulas is
\begin{align}
  \lbl{eq:ratio1}
  &\frac{\bigl(m-r+1\bigr)_{m-r+1} \, \bigl(-\frac{\mu}{2}-2m-r+1\bigr)_{m-r} \,
    \prod_{i=0}^{m-r}-(\mu+2i+6r-4)}{2 \, \bigl(-\mu-4m-2r+1\bigr)_{m-r} \,
    \bigl(-\frac{\mu}{2}-m-2r+2\bigr)_{m-r+1} \, \prod_{i=1}^{m-r} 4(2i+1)} \\
  \lbl{eq:ratio2}
  &\times\prod\limits_{i=1}^{m-r}\frac{\bigl(-\frac{\mu}{2}+2i-3m-1\bigr)_2 \,
    \bigl(-\mu+2i-6m-1\bigr)_2}{\bigl(-\frac{\mu}{2}+i-3m\bigr) \,
    \bigl(-\mu+3i-6m-2\bigr)_3}\\
  \lbl{eq:ratio3}
  &\times\prod\limits_{i=1}^{m-r}\frac{\bigl(\frac{\mu}{2}+2i+3r-\frac{3}{2}\bigr)_{i+1}^2 \,
    \bigl(\mu+2i+6r-3\bigr)_i^2 \, \bigl(-\frac{\mu}{2}+i-3m\bigr)_{i-1}^2}{%
    \bigl(\frac{\mu}{2}+i+3r-\frac{3}{2}\bigr)_{i+1}^2 \, \bigl(-\mu+2i-6m-1\bigr)_{i-1}^2
    \, \bigl(-\frac{\mu}{2}+2i-3m\bigr)_i^2}.
\end{align}
Applying \ref{itm:pochshift}, we simplify the first line and get that
\eqref{eq:ratio1}${}=\frac{(\frac{\mu}{2}+m+2r)_{m-r}}{(\mu+3m+3r)_{m-r}}$.
Using \ref{itm:pochprodstep}, \ref{itm:pochconnect} and
\ref{itm:pochshift}, we get
\eqref{eq:ratio2}${}=\frac{(\frac{\mu}{2}+m+2r)_{m-r}}{(\mu+3m+3r)_{m-r}}$,
too. Performing induction on $m-r$ yields
\eqref{eq:ratio3}${}=\frac{(\mu+3m+3r)_{m-r}^2}{(\frac{\mu}{2}+m+2r)_{m-r}^2}$.
Therefore the ratio is equal to~$1$ and the theorem holds.
\end{proof}

The following two corollaries highlight two more relationships between the $E$
and $D$ determinants not found elsewhere in this paper: they look like
special cases of \lem{biglemma1}, but in \cor{EDCorollary1} the parity of~$n$
is reversed, while in \cor{EDCorollary2} the lower indices are shifted.
In \fig{3Dfamilies}, they are depicted with the color \cEDcor. The proofs
involve simplifications of Pochhammers and other algebraic manipulations of
known formulas. The obtained ratios are remarkable because they have fixed
degrees in~$\mu$, independent of the sizes of the matrices. However, it
may not be easy to find a combinatorial explanation for these astonishingly
simple quotients, partly because these determinants perform weighted counts
(and the formulas tell us that exactly one of each two determinants is negative).

\begin{corollary}\lbl{cor:EDCorollary1}
Suppose $\mu$ is an indeterminate and $m$ is a positive integer. Then
\[
  \frac{\E11{\mu}{2m}}{\D01{\mu+3}{2m-1}} =
  \frac{-2\mu(2m-1)(\mu+2m+1)}{m(\mu+3)(\mu+2m)}.
\]
\end{corollary}
\begin{proof}
We use the formulas derived in \cite[Theorem~2]{KoutschanThanatipanonda13} for
the numerator and \cite[Proposition~10]{KoutschanThanatipanonda19} for the
denominator and we need to show that
\begin{align*}
  &\underbrace{\frac{2^{m-1}\prod\limits_{i=0}^{m-1} (i!)^2}
    {m! \prod\limits_{i=0}^{m-1}\bigl((2i)!\bigr)^2}}_{\ell_1} \cdot
  \underbrace{\vphantom{\frac{\prod\limits_{i=0}^{m-1}}{\prod\limits_{i=0}^{m-1}}}
    2^{m^2}  \bigl(\tfrac{\mu}{2} \bigr)_m \prod_{i=1}^{\lfloor \frac{m}{2}\rfloor}
    \bigl( \tfrac{\mu}{2} + 3i -\tfrac{1}{2}\bigr)_{m-2i+1}^2 \,
    \bigl( -\tfrac{\mu}{2} -3m + 3i \bigr)_{m-2i}^2}_{\ell_2} \\
  =&\underbrace{\frac{2m-1}{m\prod\limits_{i=1}^{m-1}(i)_{i+2} \, (i)_{i-1}}}_{r_1} \cdot
  \underbrace{\vphantom{\frac{2}{\prod\limits_{i=1}^{m-1}}}
    \frac{\mu(\mu+2)(\mu+2m+1)}{(\mu+3)(\mu+2m)}\prod_{i=1}^{m-1}
    \frac{\bigl(\mu+2i+1\bigr)_{i+2}\,\bigl(\mu+2i+4\bigr)_{i-1}\,
      \bigl(\frac{\mu}{2}+2i+2\bigr)_{i-1}^2}{\bigl(\frac{\mu}{2}+i+2\bigr)_{i-1}^2}}_{r_2}.
\end{align*}
It is easy to check that $\ell_1=r_1$ by simplifying
$\prod\limits_{i=0}^{m-1}(i!)^2/\bigl((2i)!\bigr)^2 = \prod\limits_{i=1}^{m-1}1/(i+1)_i^2$.
Using induction on $m$ and applying \ref{itm:pochshift},
\ref{itm:pochconnect}, \ref{itm:pochinterlace}, we get $\ell_2=r_2$, which
implies that the identity holds.
\end{proof}

\begin{corollary}\lbl{cor:EDCorollary2}
Suppose $\mu$ is an indeterminate and $m$ is a positive integer. Then
\[
  \frac{\E22{\mu}{2m+1}}{\D12{\mu+3}{2m}} =
  \frac{-\mu(2m+1)(\mu+2m+3)}{(m+1)(\mu+2m+2)}.
\]
\end{corollary}
\begin{proof}
We use the formulas derived in \cite[Theorem~5]{KoutschanThanatipanonda13} for
the numerator and \thm{KTConj20} together with \lem{switch} for the
denominator as follows:
\begin{align*}
 \E22{\mu}{2m+1}=&
 \underbrace{\frac{(-1)^m \, 2^{2m-2}\prod\limits_{i=0}^m (i!)^2}
   {(m+1)! \,\prod\limits_{i=0}^m \bigl((2i)!\bigr)^2}}_{\ell_1} \cdot
 \underbrace{\vphantom{\frac{\prod\limits_{i=0}^m}{\prod\limits_{i=0}^m}}
   2^{4m-2}\,(\mu+3)\, \bigl(\tfrac{\mu}{2}\bigr)_{m+1}}_{\ell_2} \\ \displaybreak[1]
 &\times \underbrace{2^{m^2-3m+4}
   \prod\limits_{i=1}^{\lfloor\frac{m+2}{2}\rfloor}
   \bigl(\tfrac{\mu}{2}+3i-\tfrac{1}{2}\bigr)_{m-2i+2}^2
   \prod\limits_{i=1}^{\lfloor\frac{m+1}{2}\rfloor}
   \bigl(-\tfrac{\mu}{2}-3m+3i-3\bigr)_{m-2i+1}^2}_{\ell_3},
 \\[1ex]
 \D12{\mu+3}{2m}=&
 \underbrace{\frac{(-1)^{m-1}}{(2m+1)!\,\bigl(m+1\bigr)_m
     \prod\limits_{i=0}^{m-2} \bigl(i+1\bigr)_{i+1}^2}}_{r_1}
 \cdot \underbrace{\vphantom{\frac{(-1)^m}{\prod\limits_{i=0}^{m-2}}}
   \frac{\bigl(\mu+2\bigr)_{2m+1}\,\bigl(\mu+5\bigr)_{2m-1}}
        {\bigl(\frac{\mu}{2}+3\bigr)_{m-1}}}_{r_2} \\
 &\times \underbrace{
   \prod\limits_{i=0}^{m-2}\frac{\bigl(\mu+2i+9\bigr)_i^2\,
     \bigl(\frac{\mu}{2}+2i+6\bigr)_{i+1}^2}{\bigl(\frac{\mu}{2}+i+5\bigr)_{i}^2}}_{r_3}.
\end{align*}
Then it is easy to check that $\ell_1/r_1=-(2m+1)/(m+1)$. By
\ref{itm:pochinterlace}, we have
\[
  \frac{\ell_2}{r_2}=
  \frac{\mu}{\bigl(\mu+2m+2\bigr)_2\,\bigl(\frac{\mu}{2}+\frac{5}{2}\bigr)_{m-1}^2}.
\]
Finally, by induction on $m$ and applying \ref{itm:pochshift},
\ref{itm:pochconnect}, \ref{itm:pochinterlace}, we get
$\ell_3/r_3=4\bigl(\frac{\mu}{2}+\frac{5}{2}\bigr)_{m}^2$, which implies,
after some necessary cancellations, that the identity holds.
\end{proof}


\section{Triangle Relations}
\lbl{sec:triangle}

In this section, we present some relationships between $E$-determinants
(resp.\ $D$-determinants) that cannot be expressed by nice product formulas,
since they do not factor completely. Note that the term ``triangle'' refers to
the fact that we identify triples of determinants whose pairwise ratios are
products of linear factors, implying that these three determinants share the
same ``ugly'' factor. These relationships are depicted in \fig{triangle}, and
also in \fig{3Dfamilies} (in the colors \cEtriA, \cEtriBref, \cDtriAref, and
\cDtriB), where the triangles have been ``thinned'' for better
visibility. Otherwise, the notion of triangle has no other geometric meaning
here. Some of the relationships for the $D$-determinants have already been
stated in \cite[Corollaries~22 and~23]{KoutschanThanatipanonda19}, but we
recall them here for completeness. Also, we close a small gap by showing that
none of the determinants in these triangle relations vanish. Remarkably, the
proof uses a combinatorial argument, and it seems difficult to find a purely
algebraic proof.

\begin{figure}
\centering
\begin{tikzpicture}
  \node (E1A) at (2,3.46) {$\E{2r}1\mu{2m}$};
  \node (E1B) at (0,0) {$\E{2r}1\mu{2m+1}$};
  \node (E1C) at (4,0) {$\E{2r+1}1\mu{2m}$};
  \node at (2,1.15) {\cor{triangleE1}};
  \draw (E1A) to (E1B) to (E1C) to (E1A);
  \node (E-1A) at (10,3.46) {$\E{-1}{2r}\mu{2m}$};
  \node (E-1B) at (8,0) {$\E{-1}{2r-1}\mu{2m}$};
  \node (E-1C) at (12,0) {$\E{-1}{2r}\mu{2m-1}$};
  \node at (10,1.15) {\cor{triangleEneg1}};
  \draw (E-1A) to (E-1B) to (E-1C) to (E-1A);
  \node (D1A) at (2,-2.54) {$\D{2r-1}1\mu{2m-1}$};
  \node (D1B) at (0,-6) {$\D{2r-1}1\mu{2m}$};
  \node (D1C) at (4,-6) {$\D{2r}1\mu{2m-1}$};
  \node at (2,-4.85) {\cor{triangleD1}};
  \draw (D1A) to (D1B) to (D1C) to (D1A);
  \node (D-1A) at (10,-2.54) {$\D{-1}{2r+1}\mu{2m+1}$};
  \node (D-1B) at (8,-6) {$\D{-1}{2r+1}\mu{2m}$};
  \node (D-1C) at (12,-6) {$\D{-1}{2r}\mu{2m+1}$};
  \node at (10,-4.85) {\cor{triangleD-1}};
  \draw (D-1A) to (D-1B) to (D-1C) to (D-1A);
\end{tikzpicture}
\caption{A line connecting two determinants implies that their ratios have a
  nice closed form as given in the corresponding corollary.}
\lbl{fig:triangle}
\vspace*{\floatsep}
\includegraphics[width=0.44\textwidth]{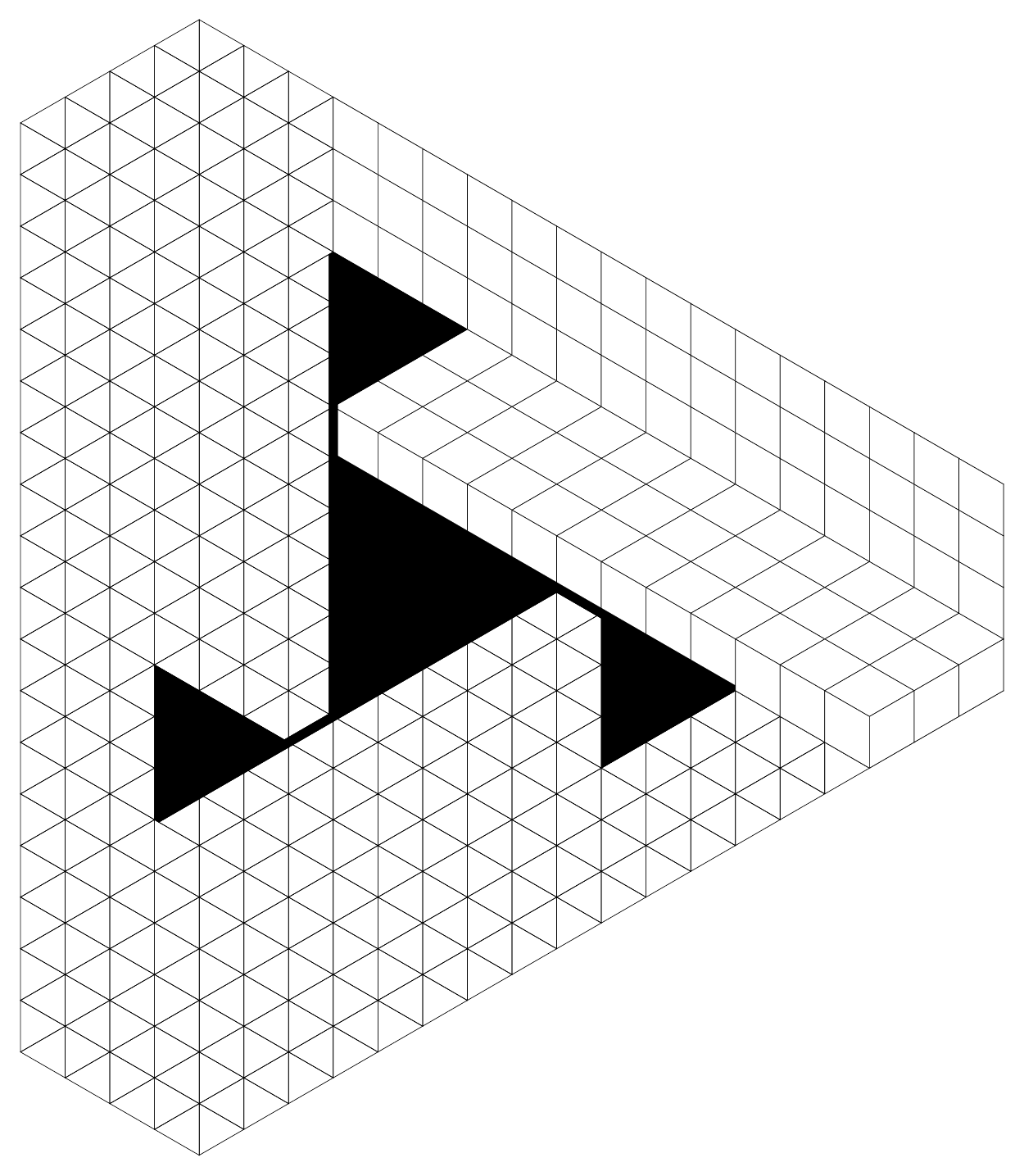}
\qquad\qquad
\raisebox{22pt}{\includegraphics[width=0.36\textwidth]{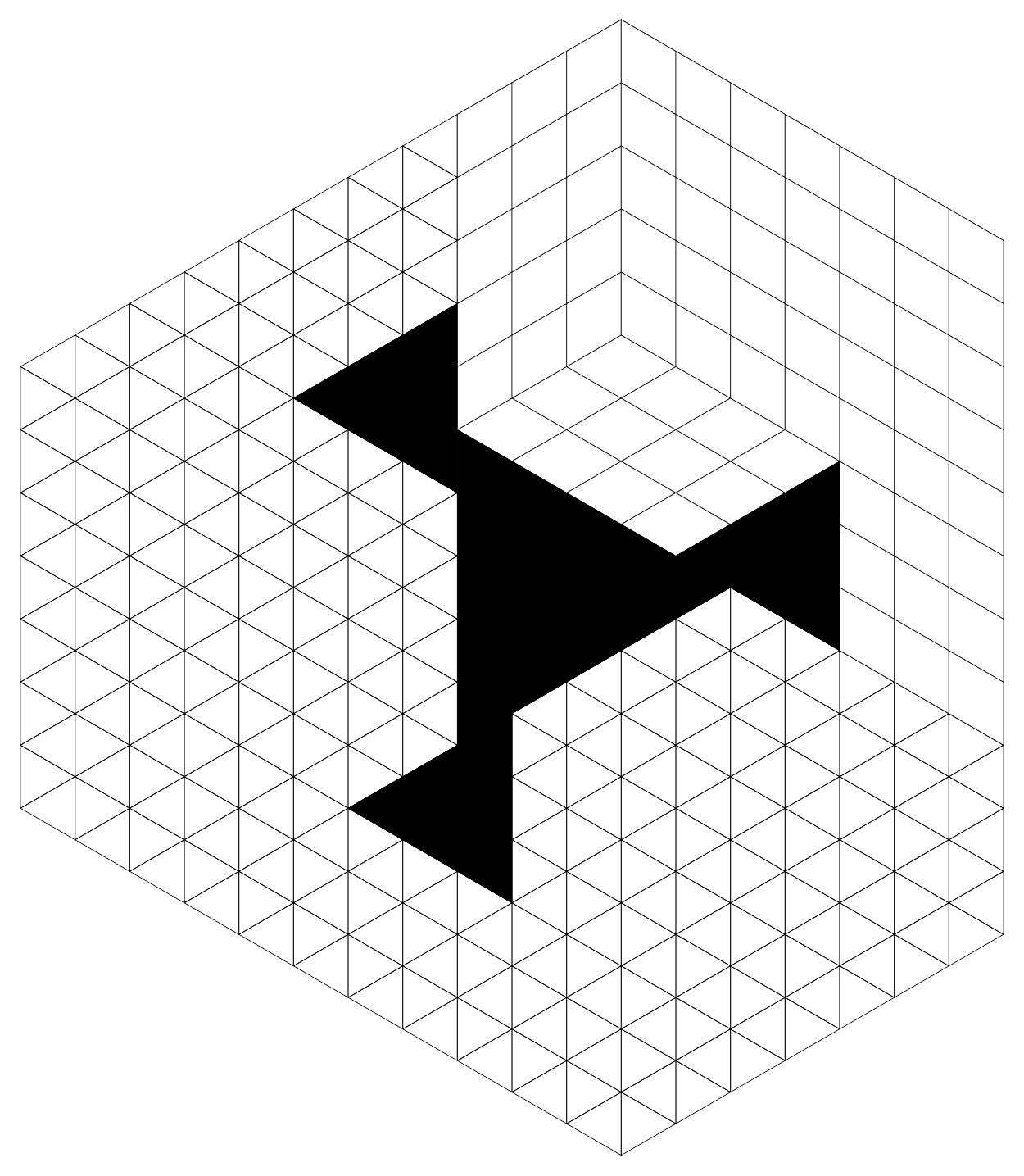}}
\caption{Hexagonal domains for the tiling problems counted by
  $\E4176$ (left) and $\E{-1}276$ (right). One particular rhombus tiling
  for each domain is sketched in the upper right third, and the tiling for
  the remaining two-thirds follow by cyclic symmetry.}
\lbl{fig:hexfortriangle}
\end{figure}

\begin{proposition}\lbl{prop:nonzero}
  Let $\mu$ be an indeterminate and let $n,r\in\Z^+$. If $n\geq 2r-1$,
  then $\E{2r}1\mu{n}$ and $\D{2r-1}1\mu{n}$ are nonzero, and if
  $n\geq 2r+1$, then $\E{-1}{2r}\mu{n}$ and $\D{-1}{2r+1}\mu{n}$ are
  nonzero (i.e., not identically zero as polynomials in~$\mu$).
\end{proposition}
\begin{proof}
  We prove the statement by appealing to the combinatorial interpretation of
  these determinants. First, we can see that $s-t$ is odd in the two $E$-determinants, 
  while for the two $D$-determinants it is even. By~\eqref{eq:sumofminors} and its analog
  \cite[(2.1)]{KoutschanThanatipanonda19}, all four determinants perform
  unweighted counts, i.e., they add up all rhombus tilings without signs.  Two examples of
  such hexagonal tiling regions are shown in \fig{hexfortriangle}.
  Note that $\E st\mu{n}$ and $\D st\mu{n}$ have the same tiling region, but
  differ only in the modus of counting (weighted vs.\ unweighted).  For each
  choice of the parameters $n,r,\mu$ there exist cyclically symmetric
  tilings (a ``canonical'' one for each type of region is shown in the
  figure), implying that these determinants cannot be (identically) zero.
\end{proof}

\begin{corollary}\lbl{cor:triangleE1}
  Let $\mu$ be an indeterminate, and let $m,r\in\Z$. If $m > r\geq 1$, then
  \begin{align*}
    \frac{E_{2r,1}^{\mu}(2m+1)}{E_{2r,1}^{\mu}(2m)}
    &= \frac{\bigl(\mu+2m+4r-1\bigr)_{m-r+1} \, \bigl(\frac{\mu}{2}+2m+r+1\bigr)_{m-r}}{
      \bigl(m-r+1\bigr)_{m-r+1} \, \bigl(\frac{\mu}{2}+m+2r\bigr)_{m-r}},
    \\[1ex]
    \frac{\E{2r}1{\mu}{2m+1}}{\E{2r+1}1{\mu}{2m}}
    &= \frac{(-1)^{m-r} \, \bigl(\frac{\mu}{2}+2m+r+1\bigr)_{m-r} \,
      \bigl(\frac{\mu}{2}+3r-\frac{1}{2}\bigr)_{m-r+1}}{
      \bigl(\frac{3}{2}\bigr)_{m-r} \, \bigl(m-r\bigr)_{m-r}},
    \\[1ex]
    \frac{\E{2r+1}1{\mu}{2m}}{\E{2r}1{\mu}{2m}}
    &= \frac{(-1)^{m-r} \, \bigl(\frac{1}{2}\bigr)_{m-r+1} \, \bigl(\mu+2m+4r-1\bigr)_{m-r+1}}{
      (2m-2r+1) \, \bigl(\frac{\mu}{2}+m+2r\bigr)_{m-r} \,
      \bigl(\frac{\mu}{2}+3r-\frac{1}{2}\bigr)_{m-r+1}}.
\end{align*}
\end{corollary}

\begin{proof}
Since the third identity is easily obtained as the quotient of the first
divided by the second, we focus on the first two identities.
We can use the Desnanot--Jacobi--Dodgson identity (see \sect{switch}) with
two different shifts of the first index:
\begin{align*}
  \E{2r-1}0{\mu}{2m+2} \E{2r}1{\mu}{2m} &=
  \E{2r-1}0{\mu}{2m+1} \E{2r}1{\mu}{2m+1} -
  \cancelto{0}{\E{2r}0{\mu}{2m+1}} \E{2r-1}1{\mu}{2m+1}, \\
  \E{2r}0{\mu}{2m+2} \E{2r+1}1{\mu}{2m} &= 
  \cancelto{0}{\E{2r}0{\mu}{2m+1}} \E{2r+1}1{\mu}{2m+1} - 
  \E{2r+1}0{\mu}{2m+1} \E{2r}1{\mu}{2m+1}.
\end{align*}
By the first identity of \lem{famA}
and~\cite[Theorem~19]{KoutschanThanatipanonda19}, it follows that
$E_{2r,0}^{\mu}$ vanishes at odd dimensions larger than $2r$, while
all other instances of $E_{s,0}^{\mu}$ are nonzero.  Together with
\prop{nonzero}, this implies that all members in the above two equations
(except the cancelled ones) are nonzero. This allows us to take quotients and
express our identities in terms of known determinants (using \lem{famA}):
\begin{alignat*}{2}
  \frac{\E{2r}1{\mu}{2m+1}}{\E{2r}1{\mu}{2m}} &=
  \frac{\E{2r-1}0{\mu}{2m+2}}{\E{2r-1}0{\mu}{2m+1}} &&=
  \frac{\D{2r}0{\mu-3}{2m+3}}{\D{2r}0{\mu-3}{2m+2}}, \\[1ex]
  \frac{\E{2r}1{\mu}{2m+1}}{\E{2r+1}1{\mu}{2m}} &=
  -\frac{\E{2r}0{\mu}{2m+2}}{\E{2r+1}0{\mu}{2m+1}} &&=
  -\frac{\D{2r-1}0{\mu+3}{2m+1}}{\D{2r}0{\mu+3}{2m}}.
\end{alignat*}
From \cite[Theorem~18]{KoutschanThanatipanonda19} we already know
$\D{2r}0{\mu}{n+1}/\D{2r}0{\mu}{n}$, so that the first quotient is
immediate. For the second identity, we combine Theorems~18 and~19
from~\cite{KoutschanThanatipanonda19} to find the ratio of
$\D{2r-1}0{\mu}{2m+1}$ and $\D{2r}0{\mu}{2m}$, and then perform some
simplifications on Pochhammer symbols by \ref{itm:pochinterlace},
\ref{itm:pochconnect} and \ref{itm:pochshift} to obtain the claimed formula.
\end{proof}

\begin{corollary}\lbl{cor:triangleEneg1}
Let $\mu$ be an indeterminate, and let $m,r\in\Z$. If $m-1 > r \geq 1$, then
\begin{align*}
    \frac{E_{-1,2r}^\mu(2m)}{E_{-1,2r}^\mu(2m-1)}
    &= \frac{\bigl(\mu+2m-3\bigr)_{2r+1} \, \bigl(\mu+2m+4r-2\bigr)_{m-r-1} \,
      \bigl(\frac{\mu}{2}+2m+r-1\bigr)_{m-r-1}}{2
      \bigl(2m-1\bigr)_{2r+1} \, \bigl(m-r-1\bigr)_{m-r-1} \,
      \bigl(\frac{\mu}{2}+m+2r-1\bigr)_{m-r-1}},
    \\[1ex]
    \frac{E_{-1,2r}^\mu(2m)}{E_{-1,2r-1}^\mu(2m)}
    &= \frac{(-1)^{m-r} \, \bigl(m-r\bigr)_{m-r} \, \bigl(\mu+2r-2\bigr)_{2m} \,
      \bigl(\frac{\mu}{2}+m+2r-\frac{3}{2}\bigr)_{m-r}}{
      \bigl(2r\bigr)_{2m} \, \bigl(\frac{\mu}{2}+3r-\frac{1}{2}\bigr)_{m-r} \,
      \bigl(\mu+3m+3r-3\bigr)_{m-r}},
    \\[1ex]
    \frac{E_{-1,2r-1}^\mu(2m)}{E_{-1,2r}^\mu(2m-1)}
    &= \frac{-(-4)^{m-r-1} \, \bigl(2r\bigr)_{2m-2r-1} \, \bigl(\frac{\mu}{2}+2m+r-2\bigr)_{m-r} \,
      \bigl(\frac{\mu}{2}+3r-\frac{1}{2}\bigr)_{m-r-1}}{\
      \bigl(m-r\bigr)_{m-r} \, \bigl(m-r-1\bigr)_{m-r-1} \, \bigl(\mu+2r-2\bigr)_{2m-2r-1}}.
\end{align*}
\end{corollary}

\begin{proof}
Similar to the proof of~\cor{triangleE1}, we apply the DJD identity to obtain
\begin{align*}
&E_{-1,2r}^{\mu}(2m)E_{0,2r+1}^{\mu}(2m-2)
=E_{-1,2r}^{\mu}(2m-1)E_{0,2r+1}^{\mu}(2m-1)-\cancelto{0}{E_{0,2r}^{\mu}(2m-1)}E_{-1,2r+1}^{\mu}(2m-1),\\[1ex]
&E_{-1,2r-1}^{\mu}(2m+1)\cancelto{0}{E_{0,2r}^{\mu}(2m-1)}
=E_{-1,2r-1}^{\mu}(2m)E_{0,2r}^{\mu}(2m)-E_{0,2r-1}^{\mu}(2m)E_{-1,2r}^{\mu}(2m),
\end{align*}
where $E_{0,2r}^{\mu}$ also vanishes at odd dimensions bigger than $2r$:
by \lem{switch}, \lem{famA}, and~\cite[Theorem~19]{KoutschanThanatipanonda19}, we obtain
$\E0{2r}{\mu}{2m-1}=(\ldots)\cdot\E{2r}0{\mu}{2m-1}=(\ldots)\cdot\D{2r-1}0{\mu+3}{2m-2}=0$.
Using a similar argument as in \cor{triangleE1}, we see that all other determinants are
nonzero, and hence we can express our identities in terms of known determinants:
\begin{alignat*}{2}
  \frac{E_{-1,2r}^\mu(2m)}{E_{-1,2r}^\mu(2m-1)}
  &= \frac{E_{0,2r+1}^\mu(2m-1)}{E_{0,2r+1}^\mu(2m-2)}
  &&= \frac{(\mu+2m-3)_{2r}}{(2m-1)_{2r}} \cdot \frac{D_{2r,0}^{\mu+3}(2m-2)}{D_{2r,0}^{\mu+3}(2m-3)},
  \\[1ex]
  \frac{E_{-1,2r}^\mu(2m)}{E_{-1,2r-1}^\mu(2m)}
  &= \frac{E_{0,2r}^\mu(2m)}{E_{0,2r-1}^\mu(2m)}
  &&= \frac{(\mu+2r-2)_{2m}}{(2r)_{2m}} \cdot \frac{D_{2r-1,0}^{\mu+3}(2m-1)}{D_{2r-2,0}^{\mu+3}(2m-1)}.
\end{alignat*}
Then applying~\cite[Theorem~18]{KoutschanThanatipanonda19}, the first identity is immediate. For the second one, combining Theorems~18 and~19 from~\cite{KoutschanThanatipanonda19} and 
performing \ref{itm:pochinterlace}, \ref{itm:pochconnect}, \ref{itm:pochshift} and \ref{itm:pochprodstep}, we can get the claimed formula. The third identity of the lemma follows from the quotient of the first divided by the second and some necessary calculations depending on  \ref{itm:pochinterlace}, \ref{itm:pochconnect} and \ref{itm:pochshift}.
\end{proof}

\begin{corollary}\lbl{cor:triangleD1}
Let $\mu$ be an indeterminate, and let $m,r\in\Z$. If $m > r \geq 1$, then
\begin{align*}
  \frac{\D{2r-1}1\mu{2m}}{\D{2r-1}1\mu{2m-1}}
  &= \frac{\bigl(\mu+2m+4r-4\bigr)_{m-r+1} \, \bigl(\frac{\mu}{2}+2m+r-\frac{1}{2}\bigr)_{m-r}}{
    \bigl(m-r+1\bigr)_{m-r+1} \, \bigl(\frac{\mu}{2}+m+2r-\frac{3}{2}\bigr)_{m-r}},
  \\[1ex]
  \frac{D_{2r,1}^\mu(2m-1)}{D_{2r-1,1}^\mu(2m)}
  &= \frac{(-1)^{m-r} \, \bigl(m-r\bigr)_{m-r} \, \bigl(m-r+1\bigr)_{m-r+1}}{
    2^{2m-2r} \, \bigl(\frac{\mu}{2}+2m+r-\frac{1}{2}\bigr)_{m-r} \,
    \bigl(\frac{\mu}{2}+3r-2\bigr)_{m-r+1}},
  \\[1ex]
  \frac{D_{2r,1}^\mu(2m-1)}{D_{2r-1,1}^\mu(2m-1)}
  &= \frac{(-1)^{m-r} \, \bigl(m-r\bigr)_{m-r} \, \bigl(\frac{\mu}{2}+m+2r-2\bigr)_{m-r}}{
    \bigl(\frac{\mu}{2}+3r-2\bigr)_{m-r+1} \, \bigl(\mu+3m+3r-3\bigr)_{m-r-1}}.
\end{align*}
\end{corollary}
\begin{proof}
The first identity is given by~\cite[Corollary~22]{KoutschanThanatipanonda19}.
For the second one, we similarly apply the DJD identity to obtain
\[
  D_{2r-1,0}^{\mu}(2m+1) D_{2r,1}^{\mu}(2m-1)
  =\cancelto{0}{D_{2r-1,0}^{\mu}(2m)} D_{2r,1}^{\mu}(2m)-D_{2r-1,1}^{\mu}(2m) D_{2r,0}^{\mu}(2m),
\]
where $D^{\mu}_{2r-1,0}$ vanishes at even dimensions no less than $2r$ by
\cite[Theorem~19]{KoutschanThanatipanonda19}. Together with \prop{nonzero},
it follows that all three determinants in this triangle relation are nonzero. Thus,
\[
  \frac{\D{2r}1\mu{2m-1}}{\D{2r-1}1\mu{2m}} = -\frac{\D{2r}0\mu{2m}}{\D{2r-1}0\mu{2m+1}}.
\]
This quotient in terms of known determinants \cite[Theorems~18
  and~19]{KoutschanThanatipanonda19} can be simplified as claimed by
\ref{itm:pochinterlace}, \ref{itm:pochconnect}, \ref{itm:pochshift} and
\ref{itm:pochprodstep}. Finally, we can obtain the third identity by combining
the first two and then performing \ref{itm:pochinterlace},
\ref{itm:pochconnect} and \ref{itm:pochshift}.
\end{proof}

\begin{corollary}\lbl{cor:triangleD-1}
Let $\mu$ be an indeterminate, and let $m,r\in\Z$. If $m > r \geq 0$, then
\begin{align*}
  \frac{\D{-1}{2r+1}\mu{2m+1}}{\D{-1}{2r+1}\mu{2m}}
  &= \frac{\bigl(\mu+2m-2\bigr)_{2r+2} \, \bigl(\mu+2m+4r+1\bigr)_{m-r-1} \,
    \bigl(\frac{\mu}{2}+2m+r+\frac{1}{2}\bigr)_{m-r-1}}{
    \bigl(2m\bigr)_{2r+2} \, \bigl(m-r\bigr)_{m-r-1} \,
    \bigl(\frac{\mu}{2}+m+2r+\frac{1}{2}\bigr)_{m-r-1}},
  \\[1ex]
  \frac{\D{-1}{2r+1}\mu{2m+1}}{\D{-1}{2r}\mu{2m+1}}
  &= \frac{(-1)^{m-r} (2)^{2m-2r-1} \, \bigl(\frac{1}{2}\bigr)_{m-r} \,
    \bigl(\frac{\mu}{2}+m+2r+1\bigr)_{m-r-1} \, \bigl(\mu+2r-1\bigr)_{2 m + 1}}{
    \bigl(2r+1\bigr)_{2m+1} \, \bigl(\frac{\mu}{2}+3r+1\bigr)_{m-r-1} \,
    \bigl(\mu+3m+3r\bigr)_{m-r}},
  \\[1ex]
  \frac{\D{-1}{2r}\mu{2m+1}}{\D{-1}{2r+1}\mu{2m}}
  &= \frac{(-1)^{m-r} \, \bigl(2r+1\bigr)_{2m-2r-1} \, \bigl(\frac{\mu}{2}+3r+1\bigr)_{m-r-1} \,
    \bigl(\frac{\mu}{2}+2m+r-\frac{1}{2}\bigr)_{m-r}}{
    \bigl(\frac{1}{2}\bigr)_{m-r} \, \bigl(m-r\bigr)_{m-r-1} \, \bigl(\mu+2r-1\bigr)_{2m-2r-1}}.
\end{align*}
\end{corollary}
\begin{proof}
The first identity is given by~\cite[Corollary~23]{KoutschanThanatipanonda19}.
Then we can use DJD to get
\[
  D_{-1,2r}^{\mu}(2m+2) \cancelto{0}{D_{0,2r+1}^{\mu}(2m)}
  =D_{-1,2r}^{\mu}(2m+1) D_{0,2r+1}^{\mu}(2m+1)-D_{-1,2r+1}^{\mu}(2m+1) D_{0,2r}^{\mu}(2m+1),
\]
where $D_{0,2r+1}$ also vanishes at even dimensions which are no less than
$2r$ by \lem{switch} and~\cite[Theorem~19]{KoutschanThanatipanonda19}. With
the help of \prop{nonzero}, we find that all other determinants in the above
DJD identity are nonzero. Then, by invoking \lem{switch}, we have
\[
  \frac{\D{-1}{2r+1}\mu{2m+1}}{\D{-1}{2r}\mu{2m+1}} =
  \frac{\bigl(\mu+2r-1\bigr)_{2m+1}}{\bigl(2r+1\bigr)_{2m+1}} \cdot
  \frac{\D{2r+1}0\mu{2m+1}}{\D{2r}0\mu{2m+1}}.
\]
The quotient, which is in terms of known determinants given in Theorem~18 and~19 in~\cite{KoutschanThanatipanonda19}, can be simplified as claimed in the second identity by \ref{itm:pochinterlace}, \ref{itm:pochconnect}, \ref{itm:pochshift} and \ref{itm:pochprodstep}. For the third identity, we combine the first two and then perform \ref{itm:pochinterlace}, \ref{itm:pochconnect} and \ref{itm:pochshift}.
\end{proof}


\section{Some Final Thoughts}
\lbl{sec:conclusion}

In this paper, we are able to tell a cohesive story about two related binomial
determinant families with signed Kronecker deltas located along a certain
diagonal in the corresponding matrices. In \fig{3Dfamilies}, we compile a
summary of this work. The reader has probably noticed the vastness of the
blank areas in the figure and may wonder if there are more results where the
determinant, viewed as a polynomial in $\mu$, factors into linear factors (we
refer to these expressions as ``nice''). For the case $s\geq 0$ and $t<0$, the
determinants evaluate to zero (see a related discussion for the $t=-1$ case in
\sect{gamma}). For the case $s<0$ and $t<0$, we have some zero determinants, a
few determinants having nice forms as discussed in
\cite[Corollary~15]{KoutschanThanatipanonda19}, and the rest are ugly. We made
a conscious decision not to include these results in order to simplify our
diagram.

In general, computer experiments for fixed $s,t$ and nontrivial $n$ led us to
rule out determinants that admit a form containing an irreducible factor of
degree greater than one. These experiments permitted us to narrow down our
search for nice expressions to within a strip of $\pm 1$ around the positive
axes. \lem{famA} and \cite[Theorems~18 and~19]{KoutschanThanatipanonda19},
together with \lem{switch}, gave us product formulas for both $D$ and $E$ on the
positive axes (with some being zero). From there, we could argue from the
perspective of the DJD identity (see \sect{switch}): as long as one of the six
determinants in the identity is zero, one of the three terms in the identity vanishes, and the remaining four determinants can be rearranged as
the equality of two ratios. This behavior manifests itself very clearly in
\sect{triangle}. In the case where the determinant that we want is paired with
the zero determinant (for example, $\E{2r-1}1{\mu}{2m-1}$ suffers from this
fate), we had to employ the holonomic ansatz to reveal other
relationships. DJD also failed to yield new results once we moved away from
this strip (e.g., where $|s|\geq 2$) because we were unable to leverage the
known results close to the axis and zero determinants. This behavior parallels
our observations from combinatorics: we have a simple interpretation if
min$(s,t)=0$, but once min$(s,t)>0$, a border line appears in our figures and
we no longer have an easy way to establish relationships between the families
or to take advantage of the symmetry of the figure to do the
counting. Therefore, while we cannot say for certain that there are no other
product formulas and ratios of the forms presented in this paper, we can
say that we searched in the places where we believe such forms are located
and anticipate that the expressions would only increase in complexity as one
moves further away from the axes.

A second remark is that we heavily relied on symbolic computation tools to
obtain our results in a reasonable amount of time. Such tools have enabled the
resolution of four out of the seven problems and conjectures discussed in
\cite[Section~5.5]{Krattenthaler05}, namely Problem~34 and
Conjectures~{35--37}. It is unknown whether or not the combinatorial
interpretation in \sect{comb} could help us prove some of our main lemmas and
theorems more easily or deduce more results.

The last thing that the engaged reader may be wondering is if the remaining
three problems and conjectures discussed in
\cite[Section~5.5]{Krattenthaler05} have also been resolved. We first note that these
three problems are of a different flavor than the ones in this paper, both in form
and in their combinatorics. In Problem~38, the goal was to compute a Pfaffian
whose entries are sums of signed binomial coefficients depending on some
entanglement of six different parameters. They arise from the
$(-1)$-enumeration of self-complementary plane partitions. This exercise was
resolved by Eisenk\"olbl in~2008 using combinatorial arguments 
\cite[Corollary]{Eisenkoelbl08}. In Conjecture~39, we see a
determinant that is a shuffling of two binomial determinants and counts
rhombus tilings of hexagons with a central triangular hole that is off center
by one unit. This was originally proposed in
\cite[Section~12]{CiucuEisenkoelblKrattenthalerZare02} and was resolved by
Rosengren in~2016 \cite[Theorem~2.1]{Rosengren16}, who used orthogonal
polynomials and analysis as the main tools. In Conjecture~40, we see another
shuffling of two different binomial determinants and Krattenthaler declared at
the time that it was one of the ``weirdest closed forms in enumeration that he
was aware of.'' It counts lozenge tilings of hexagons with cut off corners~\cite{CiucuKrattenthaler02}. This was resolved in a more general form by Ciucu
and Fischer in~2015 \cite[Theorem~2.3]{CiucuFischer15}. Their key approach is
also combinatorial. In the context of this paper, it might be relevant to
wonder whether or not a computer algebra approach to these three problems
would be efficient and would yield alternate proofs for these results.

\textbf{Acknowledgements.} C.~Koutschan and E.~Wong are both supported by the
Austrian Science Fund (FWF): F5011-N15. H.~Du is supported by the Austrian
Science Fund (FWF): F5011-N15, P31952, and P32301. These three authors would
like to acknowledge the Special Research Program ``Algorithmic and Enumerative
Combinatorics'' that encouraged and supported this research. The authors would
also like to thank Mihai Ciucu and Ilse Fischer for the update regarding
Conjecture~39 from \cite{Krattenthaler05}. E.~Wong would like to thank Manfred
Buchacher, Matteo Gallet and Ali Uncu for their friendship, support and advice
during the preparation of this manuscript, the
Symbolic Computation group at RICAM who gave ample freedom and space to conduct 
research, and especially RICAM's system administrators Florian Tischler and Wolfgang 
Forsthuber who provided technical support and infrastructure for the 
computations in this paper to complete. We would also like to thank the referees for their careful reading and helpful comments.

\bibliography{det} 
\bibliographystyle{plain}

\appendix

\section{}\lbl{ap:altproof}
Here is the alternative proof of \lem{biglemma2}.
\begin{proof}
We can see that both identities can be presented in a uniform way:
\[
  \lim\limits_{\eps\to 0}
  \left(\frac{\detnew A{s+\eps}{-1+\eps}{\mu}{n}}{\B{s-1+\eps}{-1+\eps}{\mu+3}{n-1}}\right)
  = \frac{2s(n-1)(\mu-3)(\mu+n+s-2)}{\mu(n+s)(\mu+n-3)(\mu+s-2)}
  =: R_{s,-1}^{\mu}(n),
\]
where $(A,B,s,n)=(D,E,2r,2m)$ or $(A,B,s,n)=(E,D,2r+1,2m+1)$. Like in the
proof of \lem{biglemma1}, we use an inductive argument to ensure that
$\lim_{\eps\to0}\bigl(\frac1\eps\detnew A{s+\eps}{-1+\eps}{\mu}{n}\bigr)$
exists and is nonzero. As a base case, we use
$\E{1+\eps}{-1+\eps}{\mu}{2m-2r+1}$ (see \lem{quoED1}), and as induction
hypothesis we assume from now on that
$\lim_{\eps\to0}\bigl(\frac1\eps\B{s-1+\eps}{-1+\eps}{\mu+3}{n-1}\bigr)$
exists and is nonzero.

Similar to the proof of  \lem{quoED1}, we do some basic row and column operations
for $\detnew A{s+\eps}{-1+\eps}{\mu}{n}$ by multiplying with the elementary
matrices $\matnot{L}_{n}$ from~\eqref{eq:LandR} and 
\begin{equation}\lbl{eq:newR}
\tilde{\matnot{R}}_{n}:=\begin{pmatrix} 
0&-1&0&0&0&\cdots \\ 
1&0&1&1&1&\cdots \\ 
0&0&1&1&1&\cdots \\ 
0&0&0&1&1&\cdots \\ 
\vdots & \vdots & \vdots &  & \ddots & \ddots \end{pmatrix},
\end{equation}
and then applying \lem{pascal} and the Taylor expansion~\eqref{eq:bctaylor2eps}
such that the transformed matrix $\matnot{L}_{n} \cdot \Mat{A}{s+\eps}{-1+\eps}{\mu}{n}\cdot \tilde{\matnot{R}}_{n}$ becomes
 \begin{equation}\lbl{eq:titleA}
  \left(
  \begin{array}{c:c:c}
  1+O(\eps) & -\frac{\eps}{\mu+s-2}+O(\eps^2) & \binom{\mu+s+j-3}{j-2}\pm \sum\limits_{k=1}^j\delta_{s,k-2} + O(\eps) \\
  & & \scriptstyle (3\;\leq\;j\;\leq\;n) \\
  \hdashline & & \\[-12pt]
  \begin{array}{c} \frac{\mu+s+i-5}{(\mu+s+i-4)_2} \cdot \eps +O(\eps^2) \\ \null \end{array} &
  \begin{array}{c}
    \frac{\eps}{(\mu+s+i-4)_2}+O(\eps^2) \\ \scriptstyle(2\;\leq\;i\;\leq\;n)
  \end{array} &
  \matnot{M}_{(n-1)\times(n-2)}+ O(\eps)
\end{array}
\right)
 \end{equation}
where $\pm$ is $+$ if $A=D$ and $-$ if $A=E$, and $\matnot{M}$ is the first $(n-2)$ columns of $\Mat{B}{s-1}0{\mu+3}{n-1}$.
Note that the determinantal value remained unaffected under this
transformation, and the $O(\eps)$ added to~$\matnot{M}$ means that it is added to
every entry. Since the determinant behaves like a linear function in the
columns of the matrix, the determinant of~\eqref{eq:titleA}
is equal to $\eps\cdot\tilde{A}+O(\eps^2)$, where
\[
  \tilde{A} := \det
\left(
\begin{array}{c:c:c}
  \mute{1} & \mute{-\frac{1}{\mu+s-2}} &
  \mute{\binom{\mu+s+j-3}{j-2}\pm \sum\limits_{k=1}^j\delta_{s,k-2}} \\
  & & \scriptstyle\mute{(3\;\leq\;j\;\leq\;2m+1)} \\
  \hdashline & & \\[-12pt]
  \begin{array}{c} \mute{0} \\ \null \end{array} &
  \begin{array}{c}
   \frac{1}{(\mu+s+i-4)_2} \\ \scriptstyle(2\;\leq\;i\;\leq\;n)
  \end{array} &
  \matnot{M}_{(n-1)\times(n-2)}
\end{array}
\right).
\] 
Denote by $\tilde{\matnot{A}}$ the bottom right $(n-1) \times (n-1)$ submatrix
of the above matrix, whose determinant also
equals~$\tilde{A}$. 
On the other hand, by the definition of the matrices $\Mat{D}st\mu{n}$ and $\Mat{E}st\mu{n}$ and the Taylor expansion~\eqref{eq:bctaylor2eps}, we have that
\begin{equation}\lbl{eq:titleB}
   \Mat{B}{s-1+\eps}{-1+\eps}{\mu+3}{n-1}=\left(
  \begin{array}{c:c}
  \begin{array}{c}
    \frac{\eps}{\mu+s+i-1}+O(\eps^2) \\ \scriptstyle(1\leq\;i\;\leq\;n-1)
  \end{array} &
  \matnot{M}_{(n-1)\times(n-2)}+ O(\eps)
 \end{array}
 \right).
\end{equation}
Then by linearity of the determinant in its columns, $\B{s-1+\eps}{-1+\eps}{\mu+3}{n-1}=\eps\cdot \tilde{B}+O(\eps^2)$, where $\tilde B$ is the determinant of $\tilde{\matnot{B}}=\left(
  \begin{array}{c:c}
  \begin{array}{c}
    \frac{1}{\mu+s+i-1} \\ \scriptstyle(1\leq\;i\;\leq\;n-1)
  \end{array} &
  \matnot{M}_{(n-1)\times(n-2)}
 \end{array}
 \right).$ Thus,
 \[ \lim\limits_{\eps\to 0}
  \left(\frac{\detnew A{s+\eps}{-1+\eps}{\mu}{n}}{\B{s-1+\eps}{-1+\eps}{\mu+3}{n-1}}\right)= 
  \lim\limits_{\eps\to 0}
  \left(\frac{\eps\cdot \tilde{A}+O(\eps^2)}{\eps\cdot \tilde{B}+O(\eps^2)}\right)=\frac{\tilde A}{\tilde B}.
\]
In order to compute the determinants $\tilde A$ and $\tilde B$, we choose to
expand about the first column of $\tilde{\matnot{A}}$ and
$\tilde{\matnot{B}}$, respectively, to get
\[
  \tilde A = \sum_{i=1}^{n-1}\frac{\Cof{i}{n-1}n}{(\mu+s+i-3)_2} \quad \mbox{and} \quad \tilde B = \sum_{i=1}^{n-1}\frac{\Cof{i}{n-1}n}{\mu+s+i-1},
\]
where $\Cof{i}{n-1}n$ is the $(i,n-1)$-cofactor of $\Mat{B}{s-1}0{\mu+3}{n-1}$. Since $\Cof{1}{n-1}n$ is equal to $(-1)^n B_{s,0}^{\mu+3}(n-2)$, which is nonzero by \lem{famA} and Propositions 8, 9 in~\cite{KoutschanThanatipanonda19}, we can define
\begin{equation}\lbl{eq:cquoED}
  c_{n,i}:=\frac{\Cof i{n-1}n}{\Cof1{n-1}{n}}.
\end{equation}
Let $\tilde b_{i,j}$ be the $(i,j)$-entry of $\Mat{B}{s-1}0{\mu+3}{n-1}$. Then for each fixed $n$ and $s$ with $n\geq s$, we have that $(c_{n,1}, \ldots, c_{n,n-1})$ satisfies the system of equations
\begin{equation}\lbl{eq:sys2}
\begin{cases}
  \; c_{n,1}=1, & \\[1ex]
  \; \displaystyle\sum_{i=1}^{n-1} c_{n,i}\cdot \tilde{b}_{i,j}=0, & 1\leq j \leq n-2.
\end{cases}
\end{equation}
Then the assertion will be confirmed provided that we can show that for all $n\geq s$:
\begin{equation}\lbl{eq:sys2eq3}
  \sum_{i=1}^{n-1}\frac{c_{n,i}}{(\mu+s+i-3)_2} =
  \sum_{i=1}^{n-1}\frac{c_{n,i}}{\mu+s+i-1} \cdot R_{s,-1}^\mu(n).
\end{equation}
Finally, we employ the holonomic framework to prove the following three identities:
\begin{align*}
  c_{n,1} &= 1, \\
  \sum_{i=1}^{n-1}\binom{\mu+i+j+s-2}{j-1}\cdot c_{n,i} &= \pm c_{n,j-s+1}, \qquad (1\leq j \leq n-2), \notag \\
  \sum_{i=1}^{n-1}\frac{c_{n,i}}{(\mu+s+i-3)_2} &=
  \sum_{i=1}^{n-1}\frac{c_{n,i}}{\mu+s+i-1} \cdot R_{s,-1}^\mu(n).
\end{align*}
where $c_{n,j-s+1}=0$ for $j< s$.
The computations for these identities turned out to be very similar to the
computations for the identities in \lem{biglemma1} so we will not repeat the
exposition. All of the computational details can be found in the accompanying
electronic material \cite{EM2}. However, we remark that the third identities
were much easier as there are no singularities in the certificates. Thus, the
annihilating ideal for the summation could be directly read off and certified
from the computation without further adjustments. Nevertheless, the overall computation time for these identities did not improve in comparison to the computation time for the identities in the proof of \lem{biglemma2} due to the appearance of an additional sum in the third identity. 

We can hence conclude that \eqref{eq:biglemeq3} and \eqref{eq:biglemeq4} hold,
which also completes our induction step.
\end{proof}

\end{document}